\documentclass[11pt,twoside]{article}
\usepackage{amssymb,amsmath,amsthm,amsfonts,mathrsfs,hyperref}
\usepackage{times}
\usepackage{enumerate}
\usepackage{cite,titletoc}
\usepackage{color}
\usepackage[toc,page,title,titletoc,header]{appendix}
\usepackage{multirow}

\allowdisplaybreaks

\def\cA{\mathcal A}
\def\cB{\mathcal B}

\def\cK{\mathcal K}
\def\cL{\mathcal L}

\def\cO{\mathcal O}

\def\cW{\mathcal W}

\def\sF{\mathscr F}
\def\N{\mathop{\mathbb N\kern 0pt}\nolimits}
\def\Z{\mathop{\mathbb Z\kern 0pt}\nolimits}
\def\Q{\mathop{\mathbb Q\kern 0pt}\nolimits}
\def\R{\mathop{\mathbb R\kern 0pt}\nolimits}
\def\T{\mathop{\mathbb T\kern 0pt}\nolimits}
\def\SS{\mathop{\mathbb S\kern 0pt}\nolimits}

\def\ds{\displaystyle}
\def\f{\frac}
\def\al{\alpha}
\def\supp{\mathop{\rm supp}\nolimits}

\def\p{\partial}

\def\ve{\varepsilon}
\def\dl{\delta}

\def\dive{\operatorname{div}}

\def\supp{\operatorname{supp}}
\def\ls{\lesssim}
\def\gt{\gtrsim}
\def\bn{\mathbf{n}}

\newcommand{\w}[1]{\langle {#1} \rangle}

\hypersetup{colorlinks=true,linkcolor=blue,citecolor=red,urlcolor=cyan}

\theoremstyle{plain}
\newtheorem{theorem}{Theorem}[section]

\newtheorem{lemma}[theorem]{Lemma}
\newtheorem{corollary}[theorem]{Corollary}

\newtheorem{remark}{Remark}[section]

\theoremstyle{definition}

\numberwithin{equation}{section}

\title{Almost global solutions of 2-D quadratic quasilinear wave equations in exterior domains without
the second null condition}

\author{Fei Hou$^{1}$ \quad Huicheng Yin$^{2,}$
    \footnote{Fei Hou (\texttt{fhou$@$nju.edu.cn}) and Huicheng Yin (\texttt{huicheng$@$nju.edu.cn}, \texttt{05407$@$njnu.edu.cn})
    are supported by the National key research and development program of China (No.2024YFA1013301).
    In addition, Fei Hou and Huicheng Yin are supported by the NSFC (No.~12571237, No.~12331007).}\\
    [12pt]{\small 1. School of Mathematics, Nanjing University, Nanjing, 210093, China}\\
    {\small 2. School of Mathematical Sciences and IMS, Nanjing Normal University, Nanjing, 210023, China}}

\begin{document}

\date{}
\maketitle
\thispagestyle{empty}

\begin{abstract}
For 3-D quadratic quasilinear wave equations in exterior domains
with homogeneous Dirichlet boundary conditions  and small initial data of size $\varepsilon>0$,
the authors of two papers [M. Keel, H. Smith, C.D. Sogge, Amer. J. Math. 122 (2000), no. 4, 805-842;
J. Amer. Math. Soc. 17 (2004), no. 1, 109-153] have shown that the lifespan ${\bar T}_{\varepsilon}$ of smooth solutions satisfies
${\bar T}_{\varepsilon}=\infty$ or ${\bar T}_{\varepsilon}\ge e^{C/\varepsilon}$,
depending on whether the null condition holds or not, respectively.
However, for the corresponding 2-D problem in exterior domains, it remains a difficult and unsolved problem
to determine the precise lifespan $T_{\ve}$ of smooth solutions.
In our previous papers, we have established that $T_{\ve}=\infty$ holds when both the first and second null conditions are satisfied;
and that $T_{\ve}\ge\f{C}{\ve^2}$  holds when the first null condition is violated.
In this paper, we will prove $T_{\ve}\ge e^{C/\varepsilon^2}$ when the first null condition holds but the second null condition
is violated. Combining these results, we have systematically established  sharp estimates for the lifespan $T_{\varepsilon}$
of 2-D  quadratic quasilinear wave equations in exterior domains with homogeneous Dirichlet boundary conditions.
Moreover, our series of results for quasilinear wave equations in 2-D exterior domains are consistent with the corresponding results
for the classical Cauchy problem.

\vskip 0.2 true cm


\qquad \quad lifespan, weighted energy estimate, Littlewood-Paley theory

\vskip 0.2 true cm
\noindent
\end{abstract}

\vskip 0.2 true cm

\addtocontents{toc}{\protect\thispagestyle{empty}}
\tableofcontents

\section{Introduction}

\subsection{Main results and remarks}

In this paper, we are concerned with the following IBVP (initial boundary value problem) of 2-D quadratic
quasilinear wave equations in exterior domains with homogeneous Dirichlet boundary conditions
\begin{equation}\label{QWE}
\left\{
\begin{aligned}
&\Box u=Q(\p u,\p^2u),\hspace{2.6cm} (t,x)\in[0,+\infty)\times\cK,\\
&u(t,x)=0,\hspace{3.7cm} (t,x)\in[0,+\infty)\times\p\cK,\\
&(u,\p_tu)(0,x)=(\ve u_0,\ve u_1)(x),\qquad x\in\cK,
\end{aligned}
\right.
\end{equation}
where $\p=(\p_0,\p_1,\p_2)=(\p_t,\p_{x_1},\p_{x_2})$, $\ds\Box=\p_t^2-\Delta$ with $\Delta=\p_1^2+\p_2^2$, $\cK=\R^2\setminus\cO$,
the obstacle $\cO\subset\R^2$ is compact and contains the origin, the boundary $\p\cK=\p\cO$ is smooth, $\ve>0$ is small,
$(u_0,u_1)\in C^\infty(\cK)$ and $\supp(u_0,u_1)\subset\{x\in\cK:|x|\le M_0\}$ with some fixed constant $M_0>1$.
In addition, the smooth nonlinearity is given by
$Q(\p u,\p^2u)=\ds\sum_{\alpha,\beta=0}^2Q^{\alpha\beta}(\p u)\p^2_{\alpha\beta}u$ with
\begin{equation}\label{nonlinear:def}
\begin{split}
Q^{\alpha\beta}(\p u)&=\sum_{\gamma=0}^2Q^{\alpha\beta\gamma}\p_{\gamma}u+\sum_{\gamma,\delta=0}^2Q^{\alpha\beta\gamma\delta}\p_{\gamma}u\p_{\delta}u+O(|\p u|^3),
\end{split}
\end{equation}
where the constants $Q^{\alpha\beta\gamma},Q^{\alpha\beta\gamma\delta}$ are naturally imposed to satisfy the symmetric conditions
\begin{equation}\label{sym:condition}
Q^{\alpha\beta\gamma}=Q^{\beta\alpha\gamma},
\quad Q^{\alpha\beta\gamma\delta}=Q^{\beta\alpha\gamma\delta}=Q^{\alpha\beta\delta\gamma}.
\end{equation}
Meanwhile, for $\xi=(\xi_0,\xi_1,\xi_2)\in\{\pm1\}\times\SS$, the first null condition is satisfied
\begin{equation}\label{1st:null:condition}
Q_{1st}(\xi)=\sum_{\alpha,\beta,\gamma=0}^2Q^{\alpha\beta\gamma}\xi_\alpha\xi_\beta\xi_\gamma\equiv0
\end{equation}
and the second null condition is violated
\begin{equation}\label{2nd:null:condition}
Q_{2nd}(\xi)=\sum_{\alpha,\beta,\gamma,\delta=0}^2Q^{\alpha\beta\gamma\delta}\xi_\alpha\xi_\beta\xi_\gamma\xi_\delta\not\equiv0.
\end{equation}
As usual, in order to find the smooth solutions of problem \eqref{QWE}, the compatibility conditions of the initial data $(\ve u_0,\ve u_1)$
with the boundary value $u|_{[0,+\infty)\times\p\cK}=0$
are necessary. That is,
let $J_ku=\{\p_x^{\alpha}u:0\le|\alpha|\le k\}$
and $\p_t^ku(0,x)=F_k(\ve J_ku_0,\ve J_{k-1}u_1)$ ($0\le k\le 2N$), where $F_k$ depends on $Q(\p u, \p^2u)$, $\ve J_ku_0$
and $\ve J_{k-1}u_1$
in terms of the equation and the initial data in problem \eqref{QWE},
then the compatibility conditions up to $(2N+1)-$order mean that
\begin{equation}\label{compat:condition}
\text{$F_k$ vanish on $\p\cK$ for all $0\le k\le 2N$.}
\end{equation}
Let $w_{hom}$ solve the homogeneous initial boundary value problem
\begin{equation}\label{hom:eqn}
\left\{
\begin{aligned}
&\Box w_{hom}=0,\qquad\quad (t,x)\in[0,+\infty)\times\cK,\\
&w_{hom}(t,x)=0,\qquad (t,x)\in[0,+\infty)\times\p\cK,\\
&(w_{hom},\p_tw_{hom})(0,x)=(\ve u_0,\ve u_1)(x),\qquad x\in\cK.
\end{aligned}
\right.
\end{equation}

Our main result can be stated as
\begin{theorem}\label{thm1}
Suppose that the obstacle $\cO$ is star-shaped and \eqref{nonlinear:def}-\eqref{compat:condition} hold.
Then there are two constants $\ve_0,\kappa>0$ such that when $\ve\le\ve_0$ and
\begin{equation}\label{initial:data}
\|u_0\|_{H^{2N+1}(\cK)}+\|u_1\|_{H^{2N}(\cK)}\le1 \quad\text{with $N\ge53$},
\end{equation}
problem \eqref{QWE} has a classical solution $u\in\bigcap\limits_{j=0}^{2N+1}C^{j}([0,e^{\kappa/\ve^2}), H^{2N+1-j}(\cK))$.
Moreover,
\addtocounter{equation}{1}
\begin{align}
\sum_{|a|\le2N-46}|\p Z^au|&\le C\ve\w{x}^{-1/2}\w{t-|x|}^{0.001-1},\tag{\theequation a}\label{thm1:decay:a}\\
\sum_{|a|\le2N-46}\sum_{i=1}^2|\bar\p_iZ^au|&\le C\ve\w{x}^{-1/2}\w{t+|x|}^{0.001-1},\tag{\theequation b}\label{thm1:decay:b}\\
\sum_{|a|\le2N-52}|Z^au|&\le C\ve\w{t+|x|}^{-1/2}\w{t-|x|}^{0.003-1/2},\tag{\theequation c}\label{thm1:decay:c}\\
\sum_{|a|\le2N-46}|\p Z^a(u-w_{hom})|&\le C\ve^2\w{t}^{-1},\quad|x|\le t/10,\tag{\theequation d}\label{thm1:decay:d}
\end{align}
where $Z=\{\p,\Omega\}$ with $\Omega=x_1\p_2-x_2\p_1$, $\w{x}=\sqrt{1+|x|^2}$, $|x|=r=\sqrt{x_1^2+x_2^2}$, $\bar\p_i=\p_i+\frac{x_i}{|x|}\p_t$ ($i=1,2$)
are the good derivatives (tangent to the outgoing light cone $|x|=t$).
Furthermore, the following time decay estimates of local energy hold
\addtocounter{equation}{1}
\begin{align}
\sum_{|a|\le2N-25}\|\p^au\|_{L^2(\cK_R)}\le C_R\ve(1+t)^{-1},\tag{\theequation a}\label{thm1:LEdecay:a}\\
\sum_{|a|\le2N-25}\|\p^a(u-w_{hom})\|_{L^2(\cK_R)}\le C_R\ve^2(1+t)^{-1},\tag{\theequation b}\label{thm1:LEdecay:b}
\end{align}
where $R>1$ is a fixed constant, $\cK_R=\cK\cap\{x: |x|\le R\}$ and  $C_R>0$ is a constant
depending on $R$.
\end{theorem}

\begin{remark}\label{YHC-1}
According to the expression of the quadratic nonlinearity $Q(\p u, \p^2u)$ in \eqref{nonlinear:def}, the 2-D quasilinear wave equation
in \eqref{QWE} can be divided into three cases:

\textbf{Case 1.} Both the two null conditions hold: $Q_{1st}(\xi)\equiv0$ and $Q_{2nd}(\xi)\equiv0$;

\textbf{Case 2.} The first null condition is violated: $Q_{1st}(\xi)\not\equiv0$;

\textbf{Case 3.} The first null condition holds but the second null condition is violated: $Q_{1st}(\xi)\equiv0$ and $Q_{2nd}(\xi)\not\equiv0$.

For \textbf{Case 1}, the global existence of the classical solution to problem \eqref{QWE} has been proved in \cite{HouYinYuan24}.

For \textbf{Case 2}, the lifespan $T_\ve$ of the classical solution to problem \eqref{QWE} was established in \cite{DingHouYin26}
as follows
\begin{equation}\label{life-lim-1st}
\begin{split}
&\liminf_{\ve\rightarrow0+}\ve\sqrt{T_\ve}\ge\tau_0=\Big(\max_{(\sigma, \theta)}\sum_{\alpha,\beta,\gamma=0}^2
Q^{\alpha\beta\gamma}\omega_\alpha\omega_\beta\omega_\gamma\p_\sigma^2F_0(\sigma,\theta)\Big)^{-1},\\
&\limsup_{\ve\rightarrow0+}\ve\sqrt{T_\ve}\le\tau_0\quad {\rm for~some~radial~problem},
\end{split}
\end{equation}
where $\omega=(\omega_0,\omega_1,\omega_2)=(-1,\omega_1,\omega_2)=(-1,\cos\theta,\sin\theta)$, $\theta\in [0, 2\pi]$, $x=r(\cos\theta,\sin\theta)$, $\sigma=r-t$ and
$\ds F_0(\sigma,\theta)=\lim_{r\rightarrow+\infty}\sqrt{r}w_{hom}(r-\sigma,r\cos\theta,r\sin\theta)$ with $w_{hom}$ being the solution of \eqref{hom:eqn}.

For \textbf{Case 3},  the lifespan $T_\ve$ of the solution to problem \eqref{QWE} is estimated in Theorem \ref{thm1} of
the paper.

It is pointed out that compared with the corresponding Cauchy problem of $\Box u=Q(\p u,\p^2u)$ with
$(u,\p_tu)(0,x)=(\ve u_0,\ve  u_1)(x)$ and $u_0,u_1\in C_0^{\infty}(\Bbb R^2)$, the estimates on $T_\ve$
in the above \textbf{Cases 1-3} are sharp (with respect to the related results of Cauchy problem, see \cite[Theorem 1$^\prime$]{Alinhac01a} for \textbf{Case 1}; \cite{Alinhac99b} and \cite[Theorem 6.5.8]{Hormander97book} for \textbf{Case 2};  \cite[Theorem 1]{Alinhac01a} and \cite{Alinhac01b} for \textbf{Case 3}, respectively).

Therefore, combining these results on $T_{\ve}$, we have systematically established  sharp estimates for $T_{\varepsilon}$
of 2-D  quadratic quasilinear wave equations in exterior domains with homogeneous Dirichlet boundary conditions.
Moreover, our series of results for quasilinear wave equations in 2-D exterior domains is consistent with the corresponding results
for the classical Cauchy problem.
\end{remark}

\begin{remark}\label{YHC-2}
Note that the second null condition \eqref{2nd:null:condition} is violated in Theorem \ref{thm1},
then the higher order small term $O(|\p u|^3)$ in \eqref{nonlinear:def} can be neglected.
For convenience, the term $O(|\p u|^3)$ in \eqref{nonlinear:def} will be omitted from now on.
In addition, the initial data $(u,\p_tu)(0,x)=(\ve u_0,\ve u_1)(x)$ in \eqref{QWE} can be extended to
more general case of $(u,\p_tu)(0,x)=(\tilde u_0(\ve,x),\tilde u_1(\ve,x))$ with $\tilde u_j(\ve,x)=\ve u_j(x)+\ve^2u_j^0(\ve,x)$ ($j=0,1$),
where $u_j(x)\in C_0^{\infty}(\cK)$, $u_j^0(\ve,x)\in C^{\infty}([0,1)\times\cK)$ and is supported compactly for $x$.
\end{remark}

\begin{remark}\label{YHC-02}
It is pointed out that the assumption on $N\ge53$ in \eqref{initial:data} comes from the requirement of $2N-52>N$
(see \eqref{thm1:decay:c} in Theorem \ref{thm1}). Of course, $N\ge 53$ is not sharp and we do not focus on the
precise lower bound of $N$.
\end{remark}

\begin{remark}\label{YHC-3}
In our forthcoming paper, we will establish a more accurate estimate for the lifespan $T_\ve$ of the solution to problem
\eqref{QWE}
\begin{equation}\label{life-lim-2nd}
\begin{split}
&\lim_{\ve\rightarrow0+}\ve^2\ln{T_\ve}=\tilde\tau_0=\Big(\max_{(\sigma,\theta)}\sum_{\alpha\beta\gamma,\delta=0}^2
Q^{\alpha\beta\gamma\delta}\omega_\alpha\omega_\beta\omega_\gamma\omega_\delta
\p_{\sigma}F_0(\sigma,\theta)\p_\sigma^2F_0(\sigma,\theta)\Big)^{-1}.
\end{split}
\end{equation}
Meanwhile, it will be shown that for the lifespan $T_\ve$ in \textbf{Case 2} of Remark \ref{YHC-1}
\begin{equation*}\label{life-lim-1st}
\begin{split}
&\lim_{\ve\rightarrow0+}\ve\sqrt{T_\ve}=\tau_0
\end{split}
\end{equation*}
and for the lifespan $\bar T_\ve$ corresponding to the 3-D exterior domain problem without the null condition
\begin{equation*}\label{life-lim-1st}
\begin{split}
&\lim_{\ve\rightarrow0+}\ve \ln{\bar T_\ve}=\bar \tau_0=\Big(\f12\max_{(\sigma, \omega)}\sum_{\alpha,\beta,\gamma=0}^3
Q^{\alpha\beta\gamma}\tilde\omega_\alpha\tilde\omega_\beta\tilde\omega_\gamma\p_\sigma^2\bar F_0(\sigma,\omega)\Big)^{-1},
\end{split}
\end{equation*}
where $\tilde\omega=(\omega_0,\omega)=(\omega_0,\omega_1,\omega_2,\omega_3)\in\{\pm1\}\times\Bbb S^2$,
$x=r\omega$ with $r=\sqrt{x_1^2+x_2^2+x_3^3}$, $\sigma=r-t$ and
$\ds \bar F_0(\sigma,\omega)=\lim_{r\rightarrow+\infty}r w_{hom}(r-\sigma,r\omega)$ with $w_{hom}$ being the solution
corresponding to the 3-D case of \eqref{hom:eqn}. Therefore, the precise position of shock formation for the 2-D
or 3-D exterior domain problem \eqref{QWE} can be determined. It is worth noting that although the shock formation
for the Cauchy problems of quasilinear wave equations or multidimensional hyperbolic conservation laws has been deeply studied
(see \cite{Ch1}, \cite{L-Y}, \cite{Luk}, \cite{M-Y} and \cite{Sp}), there are few results on the blowup mechanism of smooth solutions to the corresponding initial boundary
value problems in exterior domains.
\end{remark}

\begin{remark}\label{YHC-4}
The quasilinear wave equation in \eqref{QWE} not only admits a lot of interesting physical backgrounds
(for example, see \cite{HouYinYuan26b}) but also can be derived from some 2-D variational problem.
Indeed, for $z=(z^0,z^1,z^2)\in\R^{3}$,
let $\cL(z)$ be a smooth scalar Lagrangian density and assume
$\cL_{z^\alpha z^\beta}(0)=m^{\alpha\beta}$ with $(m^{\alpha\beta})_{\alpha,\beta=0}^3=\operatorname{diag}(1,-1,-1)$.
It follows from the first variation $\ds\sum_{\al=0}^2\p_\alpha\bigl(\cL_{z^\alpha}(\p u)\bigr)=0$ and Taylor expansion that
\begin{gather}\label{eq:1.7-YinH-0}
\Box u=\ds\sum_{\alpha,\beta,\gamma=0}^2g^{\alpha\beta\gamma}\p_\gamma  u\p_{\alpha\beta}^2  u
+\ds\sum_{\alpha,\beta,\gamma,\nu=0}^2h^{\alpha\beta\gamma\nu}\p_\gamma  u\p_\nu u\p_{\alpha\beta}^2 u+Q(\p u),
\end{gather}
where $g^{\alpha\beta\gamma}=-\cL_{z^\alpha z^\beta z^\gamma}(0)$, $h^{\alpha\beta\gamma\nu}$ are constants with $h^{\alpha\beta\gamma\nu}=
h^{\beta\alpha\gamma\nu}=h^{\alpha\beta\nu\gamma}$, and $Q(\p u)
=O(|\p u|^3|\p^2 u|)$ is a smooth function. If the first null condition holds but the second null condition
is violated in \eqref{eq:1.7-YinH-0},
then \eqref{eq:1.7-YinH-0}
fulfills the requirement on the quasilinear wave equation in \eqref{QWE}.
\end{remark}

\begin{remark}\label{YHC-004}
Consider the corresponding 2-D or 3-D Cauchy-Neumann problem of \eqref{QWE}
\begin{equation}\label{QWE-N}
\left\{
\begin{aligned}
&\Box u=Q(\p u,\p^2u),\hspace{2.6cm} (t,x)\in[0,+\infty)\times\cK,\\
&\frac{\p}{\p\bn}u(t,x)=0,\hspace{3.1cm} (t,x)\in[0,+\infty)\times\p\cK,\\
&(u,\p_tu)(0,x)=(\ve u_0,\ve u_1)(x),\qquad x\in\cK,
\end{aligned}
\right.
\end{equation}
where $\bn=(0, n_1,n_2)=(0, n_1(x),n_2(x))$ represents the unit outer normal direction of the boundary $[0,+\infty)\times\p\cK$.
In addition, the following admissible condition on $Q(\p u,\p^2u)$ for Neumann boundary
value problem \eqref{QWE-N} is imposed:
\begin{equation}\label{YHCC-2N}
\begin{split}
&\text{For any $C^1(\Bbb R^2)$ smooth functions $w$ and $\psi$ verifying $\frac{\p}{\p\bn}w
=\frac{\p}{\p\bn}\psi=0$ on $[0,+\infty)\times\p\cK$},\\
&\text{it holds that on $[0,+\infty)\times\p\cK$, $\sum_{\alpha,\beta,\mu=0}^2Q^{\alpha\beta\mu}n_{\alpha}\p_{\beta}w\p_{\mu}\psi\equiv0$
and $\sum_{\alpha,\beta,\mu,\nu=0}^2
Q^{\alpha\beta\mu\nu}n_{\alpha}\p_{\beta}w\p_{\mu}\psi\p_{\nu}\psi\equiv 0$.}
\end{split}
\end{equation}
Under the compatibility conditions of the initial data $(\ve u_0,\ve u_1)$
with the boundary value $\frac{\p}{\p\bn}u|_{[0,+\infty)\times\p\cK}=0$, there are some interesting
results on the lifespan of problem \eqref{QWE-N} (see \cite{Chen89}, \cite{Godin89}, \cite{HouYinYuan26b},
\cite{KKL13}, \cite{Li-Yin18},  \cite{Secchi-Shibata03}-\cite{Shib-Tsuti86} and \cite{Yuan22}). However, compared with the Cauchy-Dirichlet problem \eqref{QWE},
the results on the precise lifespan of the Cauchy-Neumann problem \eqref{QWE-N} are still rather incomplete.
\end{remark}

\subsection{Previous results}

The long time existence of small data smooth solutions to the quasilinear wave equations in the whole space $\Bbb R^n$
or in the exterior domains have been extensively studied.
We now give systematic reviews on the related results.
Consider the following Cauchy problem of the quasilinear wave equation
\begin{equation}\label{QWE:Cauchy}
\left\{
\begin{aligned}
&\Box w+\sum_{\alpha,\beta,\mu=0}^nQ^{\alpha\beta\mu}\p^2_{\alpha\beta}w\p_{\mu}w
+\sum_{\al,\beta,\mu,\nu=0}^nQ^{\alpha\beta\mu\nu}\p^2_{\alpha\beta}w\p_{\mu}w\p_{\nu}w=0 ,\quad t>0,\, x\in\R^n,\\
&(w,\p_tw)(0,x)=(\ve w^C_0,\ve w^C_1)(x),\qquad \qquad \qquad \qquad \qquad \qquad \qquad \qquad x\in\R^n
\end{aligned}
\right.
\end{equation}
and the initial boundary value problem in exterior domains
\begin{equation}\label{QWE:exterior}
\left\{
\begin{aligned}
&\Box w+\sum_{\alpha,\beta,\mu=0}^nQ^{\alpha\beta\mu}\p^2_{\al\beta}w\p_{\mu}w
+\sum_{\alpha,\beta,\mu,\nu=0}^nQ^{\alpha\beta\mu\nu}\p^2_{\alpha\beta}w\p_{\mu}w\p_{\nu}w=0 ,\quad t>0,\, x\in\cK,\\
&w=0,\qquad \qquad \qquad \qquad \qquad \qquad \qquad \qquad \qquad \qquad \qquad\qquad\qquad t>0, x\in\p\cK,\\
&(w,\p_tw)(0,x)=(\ve w^B_0,\ve w^B_1)(x),\qquad \qquad \qquad \qquad \qquad \qquad \qquad \qquad x\in\cK,
\end{aligned}
\right.
\end{equation}
where $n\ge 2$, $(w^C_0,w^C_1)\in C_0^{\infty}(\R^n)$, $(w^B_0,w^B_1)\in C_0^{\infty}(\cK)$, $\ve>0$ is small,
$\cK=\R^n\setminus\cO$, the obstacle $\cO\subset\R^n$ is compact and contains the origin, and the boundary $\p\cK=\p\cO$ is smooth.
Let $T_{\ve}^{C}$ and  $T_{\ve}^{D}$ denote the lifespan of smooth solutions to the Cauchy problem \eqref{QWE:Cauchy} and the Dirichlet boundary problem in \eqref{QWE:exterior}, respectively.
For $n=3$ in \eqref{QWE:Cauchy} and \eqref{QWE:exterior},
set $Q(\xi)=\ds\sum_{\alpha,\beta,\gamma=0}^3Q^{\alpha\beta\gamma}\xi_{\alpha}\xi_{\beta}\xi_{\gamma}$
for
$\xi=(\xi_0, \xi_1, \xi_2, \xi_3)$
\noindent $\in\{\pm1\}\times\SS^2$.
The previously known results on $T_{\ve}^{C} $ and $T_{\ve}^{D}$ can be briefly described in the following table:
\vskip 0.1 true cm

\begin{table}[!h]
\renewcommand\arraystretch{2.0}
\begin{tabular}{|c|c|c|c|}
\hline
Dimensions & Quadratic nonlinearity & Cauchy problem & Dirichlet problem\\
\hline
$n\ge4$ &  & $\ds T_{\ve}^C=\infty$ & $\ds T_{\ve}^D=\infty$ \\
\hline
\multirow{2}{*}{$n=3$} & $Q(\xi)\equiv0$ & $\ds T_{\ve}^C=\infty$ & $\ds T_{\ve}^D=\infty$\\
\cline{2-4}
 & $Q(\xi)\not\equiv0$  & $\ds T_{\ve}^C\ge e^{\frac{C}{\ve}}$ & $\ds T_{\ve}^D\ge e^{\frac{C}{\ve}}$
 \\
\hline
\multirow{3}{*}{$n=2$} & $Q_{1st}(\xi)\equiv0, Q_{2nd}(\xi)\equiv0$ & $\ds T_{\ve}^C=\infty$ & $\ds T_{\ve}^D=\infty$ \\
\cline{2-4}
& $Q_{1st}(\xi)\not\equiv0$  & $\ds T_{\ve}^C\ge C/\ve^2$ & $\ds T_{\ve}^D\ge C/\ve^2$ \\
\cline{2-4}
 & $Q_{1st}(\xi)\equiv0, Q_{2nd}(\xi)\not\equiv0$ & $\ds T_{\ve}^C\ge e^{C/\ve^2}$ & {\color{red}$\ds T_{\ve}^D\ge e^{C/\ve^2}$~Open}\\
\hline
\end{tabular}
\end{table}

Next we provide a detailed account of the conclusions in the table above.

$\bullet$  When $n\ge4$, the global existence of small data smooth solutions to both the problems \eqref{QWE:Cauchy} and \eqref{QWE:exterior} has been well known, i.e., $T_{\ve}^{C}, T_{\ve}^{D}=\infty$.
The Cauchy problem \eqref{QWE:Cauchy} is referred to \cite{Hormander97book,Klainerman80,KP83,Li-Chen} and see \cite{Chen89,Hayashi94,MetcalfeSogge06,Shib-Tsuti86} for the exterior domain problem \eqref{QWE:exterior}.

\vskip 0.1 true cm

$\bullet$ When $n=3$ and the null condition is satisfied (i.e., $Q(\xi)\equiv0$), the global existence of small data smooth solutions for both the Cauchy problem \eqref{QWE:Cauchy} and the initial boundary value problem \eqref{QWE:exterior} can be found in \cite{Christodoulou86,Klainerman} and \cite{Godin95,KK08,KSS00,MetcalfeSogge05}, respectively.
While the null condition is violated (that means, $Q(\xi)\not\equiv0$), the sharp lower bounds $T^C_{\ve}, T_{\ve}^D\ge e^{C/\ve}$ have been established for problem \eqref{QWE:Cauchy} (see \cite{JohnKlainerman84}) and problem \eqref{QWE:exterior} (one can be referred to \cite{Godin89,KSS04,MetcalfeSogge06}), respectively.

\vskip 0.1 true cm

$\bullet$ When $n=2$ and all $Q^{\alpha\beta\mu}$ vanish (i.e., $Q^{\alpha\beta\mu}\equiv0$ for all $\alpha,\beta,\mu=0,1,2$),
the global smooth solutions with the second null condition ($Q_{2nd}(\xi)\equiv0$) and the lower bound of the lifespan without the second null condition ($Q_{2nd}(\xi)\not\equiv0$) for both the problems \eqref{QWE:Cauchy} and \eqref{QWE:exterior} have been shown,
see \cite{Katayama95,Kubo13,HouYinYuan26a} for the global existence, see \cite{Godin93,Kovalyov87,Hoshiga95}
for $T_{\ve}^{C}\ge e^{C/\ve^2}$ and \cite{Kubo14} for $T_{\ve}^{D}\ge e^{C/\ve^2}$, respectively.

\vskip 0.2 true cm

$\bullet$ When $n=2$ and $Q^{\alpha\beta\mu}\neq0$ for some $\alpha,\beta,\mu=0,1,2$, as is well known from \cite{Alinhac99b,Alinhac01a,Alinhac01b,Hormander97book} that $T^C_{\ve}=\infty$ for \textbf{Case 1} (i.e., $Q_{1st}(\xi)\equiv0$ and $Q_{2nd}(\xi)\equiv0$), the lifespan verifies $\ds\lim_{\ve\rightarrow 0+}\ve\sqrt{T_\ve^C}=\tau_0>0$ for \textbf{Case 2} (namely, $Q_{1st}(\xi)\not\equiv0$ in \eqref{1st:null:condition}),
$\ds\lim_{\ve\rightarrow 0+}\ve^2\ln{T_\ve^C}=\tilde\tau_0>0$ for \textbf{Case 3} (that is, $Q_{1st}(\xi)\equiv0$ but $Q_{2nd}(\xi)\not\equiv0$ in \eqref{2nd:null:condition}), respectively.
For the exterior domain problem \eqref{QWE:exterior}, $T^D_{\ve}=\infty$ is shown in \cite{HouYinYuan24} for \textbf{Case 1}.
Recently, the authors in \cite{LRX25} proved $T^D_\ve\ge\frac{C}{\ve^2|\ln \ve|^3}$ for \textbf{Case 2}, which contains $|\ln \ve|^{-3}$ loss than the corresponding Cauchy problem.
In \cite{DingHouYin26}, $\ds\liminf_{\ve\rightarrow 0+}\ve\sqrt{T_\ve^D}\ge\tau_0$ is established for \textbf{Case 2} and $\ds\limsup_{\ve\rightarrow 0+}\ve\sqrt{T_\ve^D}\le\tau_0$ is also shown for some radial problem.

\vskip 0.1 true cm

This paper focuses on the remaining \textbf{Case 3} of the 2-D quasilinear wave equations in exterior domains, namely, we will
prove the sharp lifespan $\ds T_{\ve}^D\ge e^{C/\ve^2}$ for problem \eqref{QWE} under the first null condition and without the second null condition.

\subsection{Sketches of the proof}

Some essential difficulties for studying the 2-D exterior domain problem of quasilinear wave equations have been illustrated in \cite{HouYinYuan24}.
Note that the Lorentz boost fields $t\p_i+x_i\p_t$ $(i=1,2)$, which play key roles in the related Cauchy problem,
are no longer applicable for the exterior domain problems since these fields contain the large  normal component factor $t$,
see \cite[page 110]{KSS04}.
On the other hand, it is also difficult to apply the scaling vector field $t\p_t+r\p_r$ for the 2-D exterior domain problem
since the local energy of 2-D wave equation decays in the rate $(1+t)^{-1}$, while the local energy of the related 3-D problem decays exponentially in time, one can be referred to \cite{Morawetz75,Vainberg75}.
It is emphasized that the exponential decay rate of local energy for the 3-D problem is essential in proving the
key Lemma 2.8 of \cite{MetcalfeSogge05}.
Moreover, as discussed in \cite[Page 485]{Hepd-Metc23}, the crucial KSS estimates in 3-D problem \cite{KSS04,MetcalfeSogge05,MetcalfeSogge06} fail in the corresponding 2-D problem.

When both the first and second null conditions hold, the global existence of smooth solutions to the 2-D quadratic quasilinear wave equations in exterior domains has been proved in \cite{HouYinYuan24} and the following crucial pointwise estimates and local energy decay estimates have also been established
\begin{equation*}\label{hyy24:decay}
\begin{split}
\sum_{|a|\le N+1}|\p Z^au|&\le C\ve\w{x}^{-1/2}\w{t-|x|}^{-1},\\
\sum_{|a|\le N+1}|\bar\p Z^au|&\le C\ve\w{x}^{-1/2}\w{t+|x|}^{0.001-1},\\
\end{split}
\end{equation*}
\begin{equation}\label{hyy24:decay}
\begin{split}
\sum_{|a|\le N+1}|Z^au|&\le C\ve\w{t+|x|}^{0.001-1/2}\w{t-|x|}^{0.001-1/2},\\
\sum_{|a|\le2N-41}\|\p^au\|_{L^2(\cK_R)}&\le C_R\ve(1+t)^{-1}.
\end{split}
\end{equation}

In the paper, in order to prove Theorem \ref{thm1} for $t\le e^{\kappa/\ve^2}$, we intend to prove some estimates analogous
to \eqref{hyy24:decay}
through overcoming the difficulties mentioned above. However, due to the absence of the second null condition,
only a slightly weaker estimate  \eqref{thm1:decay:a} can be obtained than that in the first line of \eqref{hyy24:decay},
which will be enough to show Theorem \ref{thm1}. In this process, we identify the following favorable structure for the cubic nonlinearity
in \eqref{QWE}
\begin{equation}\label{intro:cubic:decom}
\begin{split}
&\sum_{\alpha,\beta,\gamma,\delta=0}^2Q^{\alpha\beta\gamma\delta}\p^2_{\alpha\beta}u\p_{\gamma}u\p_{\delta}u
=\frac13\underbrace{\p_t[\sum_{\alpha,\beta,\gamma,\delta=0}^2Q^{\alpha\beta\gamma\delta}\omega_{\alpha}\omega_{\beta}\omega_{\gamma}
\omega_{\delta}(\p_tu)^3]}
_{\rm \bf the ~divergence~part}\\
&\qquad+\sum_{\alpha,\beta,\gamma,\delta=0}^2Q^{\alpha\beta\gamma\delta}\underbrace{(\p^2_{\alpha\beta}u\p_{\gamma}u\p_{\delta}u
-\omega_{\alpha}\omega_{\beta}\omega_{\gamma}\omega_{\delta}\p_t^2u(\p_tu)^2)}_{\rm \bf the ~second ~null~condition~part},
\end{split}
\end{equation}
which is composed of a divergence part and a second null condition part.

Based on \eqref{intro:cubic:decom} and the proof procedure of \eqref{hyy24:decay} for both null conditions,
we may roughly focus on the following special cubic nonlinear problem which does not contain the quadratic nonlinearity
with the first null condition
and the cubic nonlinearity with the second null condition
\begin{equation}\label{cubic:Cauchy:model}
\Box w=\p_t[(\p_tw)^3].
\end{equation}
To solve \eqref{cubic:Cauchy:model}, at first, through investigating the solution of \eqref{hom:eqn} by some delicate analyses, we
can obtain the following new estimate (see Section \ref{sect3} below)
\begin{equation}\label{intr:hom:decay}
|\p Z^aw_{hom}(t,x)|\le C\ve\w{x}^{-1/2}\w{t-|x|}^{-1},\quad x\in\cK.
\end{equation}
This, together with the weighted $L^\infty-L^\infty$ estimate for the solution $w-w_{hom}$ of the related nonlinear wave equation with
divergence form in Lemma \ref{lem:pw:div}, yields that for $t\le e^{\kappa/\ve^2}$,
\begin{equation*}
|\p w(t,x)|\le C\ve^3\ln(2+t+|x|)\w{x}^{-1/2}\w{t-|x|}^{-(1-\eta)},\quad \eta\in(0,1/2).
\end{equation*}
Thus, \eqref{thm1:decay:a} can be attained for $t\in[0,e^{\kappa/\varepsilon^2}]$, but which does not admit the same good spacetime decay rate
as in the first line of \eqref{hyy24:decay} due to the appearance of the large factor $\ln(2+t+|x|)$.

Next, we explain how to control the energy bounds of problem \eqref{QWE} for $t\in [0,e^{\kappa/\varepsilon^2}]$, which
is one of the key points in proving Theorem \ref{thm1}.
Due to the essential differences from \cite{HouYinYuan24}, the desired energy estimate can not be established only by  \eqref{thm1:decay:a}-\eqref{thm1:decay:c}.
Indeed, in the region $|x|\le t/10$, it follows from \eqref{thm1:decay:a} and the standard energy
estimate for the quadratic nonlinear wave equation that
\begin{equation*}
\frac{d}{dt}E(t)\le \frac{C\ve E(t)}{(1+t)^{1-\delta_0}},
\end{equation*}
where $\delta_0=10^{-3}$ and $\ds E(t)=\sum_{|a|\le N_0}\|\p\p^au(t)\|_{L^2}$.
Together with  Gronwall's Lemma, this yields
\begin{equation*}
E(t)\le E(0)e^{C\ve(1+t)^{\delta_0}},
\end{equation*}
which means that $E(t)$ grows exponentially in time and can not be uniformly bounded  for $t\in [0,e^{\kappa/\varepsilon^2}]$.
On the other hand, if one intends to use  the estimate analogous to the first line of \eqref{hyy24:decay},
then it holds that in the region $|x|\le t/10$,
\begin{equation}\label{intr:energy1}
\begin{split}
\frac{d}{dt}E(t)&\le \frac{C\ve E(t)}{1+t}\quad\text{and thus $E(t)\le E(0)e^{C\ve\ln(2+t)}=E(0)(2+t)^{C\ve}$,}
\end{split}
\end{equation}
which also means that $E(t)$ is not bounded for $t\in[0,e^{\kappa/\ve^2}]$.
If one still continues to apply the energy estimate \eqref{intr:energy1} as in \cite{HouYinYuan24}, then it follows from the Sobolev embedding that
\begin{equation}\label{intr:x-decay}
|\p Z^au(t,x)|\le C\ve\w{x}^{-1/2}(2+t)^{C\ve}.
\end{equation}
Therefore, the estimate \eqref{thm1:decay:a} can not be derived by the weighted $L^\infty-L^\infty$ estimate in Lemma \ref{lem:pw:div} with \eqref{intr:hom:decay} and \eqref{intr:x-decay}.
To overcome this essential difficulty, a better decay estimate (including the small factor $\ln(2+t)^{-2}$ for large $t$) will be established  
for the solution $w_{hom}$ of \eqref{hom:eqn}
(see  Section \ref{sect3} below)
\begin{equation}\label{hom:decay:ac}
|\p Z^aw_{hom}(t,x)|\ls\ve(1+t)^{-1}\ln(2+t)^{-2},\quad |x|\le t/4.
\end{equation}
If the better estimate \eqref{thm1:decay:d} is shown, then it follows from
some related methods in \cite{HouYinYuan24} and \eqref{hom:decay:ac} that
for $t\in [0,e^{\kappa/\ve^2}]$,
$$\frac{d}{dt}E(t)\le CE(t)(\ve(1+t)^{-1}\ln(2+t)^{-2}+\ve^2(1+t)^{-1})$$
and
$$E(t)\le E(0)e^{C(\ve+\ve^2\ln(2+t))}\le E(0)e^{C(\ve+\kappa)},$$
which leads to the uniform bound of the energy $E(t)$.

To prove the key estimate \eqref{thm1:decay:d}, we will refine the estimate in Lemma \ref{lem:pw:div} in the region $|x|\le t/10$
by the Fourier method and the Littlewood-Paley theory (see Lemma \ref{lem:dt:impr}). Namely,
for the solution  $w$ of the nonhomogeneous linear Cauchy problem
\begin{equation*}
\left\{
\begin{aligned}
&\Box w=F(t,x),\qquad \qquad \qquad (t,x)\in[0,+\infty)\times\R^2,\\
&(w,\p_tw)(0,x)=(0,0),\qquad\quad x\in\R^2,
\end{aligned}
\right.
\end{equation*}
where $\supp_yF(s,y)\subset\{y: \frac{9s}{10}\le|y|\le s+M_0\}$ and $M_0>0$, then
it holds that for any constant $\mu>0$ and $|x|\le\frac{t}{10}$,
\begin{equation}\label{dt:impr-0}
|\p_tw(t,x)|\ls\f{1}{\w{t}}\sup_{(s,y)\in[0,t]\times\R^2}\w{y}^{1/2}\w{s+|y|}(\min\{\w{y},\w{s-|y|}\})^{1+\mu}|(1-\Delta)^2F(s,y)|,
\end{equation}
which eliminates the large factor $\ln(2+t+|x|)$ on the right hand side of \eqref{dt:impr-0} appeared in Lemma \ref{lem:pw:div}.
Combining \eqref{dt:impr-0} with  the energy estimate of $v=u-w_{hom}$ (see Section \ref{sect5})
and the Sobolev embedding, one has the following important pointwise estimate (see \eqref{BA3:impr:C2})
\begin{equation}\label{YHC-5}
\sum_{|a|\le2N-15}|\p Z^av(t,x)|\le C\ve\w{x}^{-1/2}.
\end{equation}
Based on \eqref{YHC-5}, through applying some related methods for treating the exterior domain problem of
the quadratic quasilineat wave equation with null forms in \cite{HouYinYuan24} (see Lemmas \ref{lem:Q0i:structure}-\ref{lem:eqn:null:div} below) and utilizing the basic structure expression \eqref{intro:cubic:decom}, the pointwise estimates \eqref{thm1:decay:a}-\eqref{thm1:decay:d} will be obtained by the weighted $L^\infty-L^\infty$ estimates of the IBVP established in Section \ref{sect4} following the steps below
\begin{equation*}
\begin{split}
\sum_{|a|\le2N-15}|Z^av|\le C\ve\w{t+|x|}^{0+} & \Rightarrow \sum_{|a|\le2N-17}|\bar\p Z^av|\le C\ve\w{x}^{-1+},|x|\ge1+t/2, \\
& \Downarrow \\
\sum_{|a|\le2N-21}|Z^av|\le C\ve\w{t+|x|}^{-1/2}, & \sum_{|a|\le2N-21}|\p Z^av|\le C\ve\w{x}^{-1/2}\w{t-|x|}^{-1/2}\w{t}^{0+}, \\
& \Downarrow \\
\sum_{|a|\le2N-27}|Z^av|\le C\ve\w{t+|x|}^{-1/2}\w{t-|x|}^{-0.4}, & \sum_{|a|\le2N-23}|\bar\p Z^av|\le C\ve\w{x}^{-3/2+},|x|\ge1+t/2,\\
\sum_{|a|\le2N-25}\|\p^av\|_{L^2(\cK_R)}&\le C\ve^2(1+t)^{-1},\\
& \Downarrow \\
\sum_{|a|\le2N-36}|\p Z^av|&\le C\ve\w{x}^{-1/2}\w{t-|x|}^{-1+}\w{t}^{0+},\\
& \Downarrow \\
\sum_{|a|\le2N-46}|\p Z^av|\le C\ve^2\w{t}^{-1},|x|\le t/10,&
\sum_{|a|\le2N-46}|\p Z^av|\le C\ve\w{x}^{-1/2}\w{t-|x|}^{-1+}, \\
& \Downarrow \\
\sum_{|a|\le2N-46}|\bar\p Z^av|\le C\ve\w{x}^{-1/2}\w{t+|x|}^{-1+} & \Rightarrow \sum_{|a|\le2N-52}|Z^av|\le C\ve\w{t+|x|}^{-1/2}\w{t-|x|}^{-1/2+},
\end{split}
\end{equation*}
whose details are given in Section \ref{sect6}.
Due to these obtained conclusions, Theorem \ref{thm1} is eventually proved.

\vskip 0.1 true cm

The paper is organized as follows.
In Section \ref{sect2}, several basic lemmas, some important results and good properties of null structures are stated.
In Section \ref{sect2.5}, the key Lemma \ref{lem:pw:div} is proved in details.
Two crucial estimates are derived in Section \ref{sect3}.
In Section \ref{sect4}, some new pointwise estimates for the exterior domain problem are established.
In Section \ref{sect5}, the energy estimates are shown to remain uniformly bounded for $t\in[0,e^{\kappa/\varepsilon^2}]$.
In Section \ref{sect6}, the key pointwise estimates are established step by step and the proof of Theorem \ref{thm1} is finished.

\vskip 0.2 true cm

\noindent \textbf{Notations:}
\begin{itemize}
  \item $\cK=\R^2\setminus\cO$, the obstacle $\cO\subset\R^2$ is compact and contains the origin,
  $\cK_R=\cK\cap\{x: |x|\le R\}$, $R>1$ is a fixed constant which may be changed in different places.
  \item Without loss of generality, assume $\p\cK\subset\{x:c_0<|x|<1/2\}$ with $c_0>0$ being a constant.
  \item The cutoff function $\chi_{[a,b]}(x)=\chi_{[a,b]}(|x|)$, $\chi_{[a,b]}(s)\in C^\infty(\R)$ with $a,b\in\R$, $0<a<b$, $0\le\chi_{[a,b]}(s)\le1$ and
  \begin{equation*}
    \chi_{[a,b]}(s)=\left\{
    \begin{aligned}
    0,\qquad s\le a,\\
    1,\qquad s\ge b.
    \end{aligned}
    \right.
  \end{equation*}
  \item $\w{x}=\sqrt{1+|x|^2}$.
  \item $\N_0=\{0,1,2,\cdots\}$ and $\N=\{1,2,\cdots\}$.
  \item $\p_0=\p_t$, $\p_1=\p_{x_1}$, $\p_2=\p_{x_2}$, $\p_x=\nabla_x=(\p_1,\p_2)$, $\p=(\p_t,\p_1,\p_2)$.
  \item For $|x|>0$, define $\bar\p_i=\p_i+\frac{x_i}{|x|}\p_t$ ($i=1,2$) and $\bar\p=(\bar\p_1,\bar\p_2)$.
  \item $r=|x|$, $x=r(\cos\theta,\sin\theta)$ with $\theta\in [0, 2\pi]$, $\Omega_{ij}=x_i\p_j-x_j\p_i$ with $i,j=1,2$, $\Omega=\p_{\theta}=\Omega_{12}$,  $Z=\{Z_1,Z_2,Z_3,Z_4\}=\{\p_t,\p_1,\p_2,\Omega\}$,
  $\tilde Z=\chi_{[1/2,1]}(x)Z$.
  \item $\p_x^a=\p_1^{a_1}\p_2^{a_2}$ for $a\in\N_0^2$, $\p^a=\p_t^{a_1}\p_1^{a_2}\p_2^{a_3}$ for $a\in\N_0^3$
      and $Z^a=Z_1^{a_1}Z_2^{a_2}Z_2^{a_3}Z_4^{a_4}$ for $a\in\N_0^4$.
  \item The commutator $[X,Y]=XY-YX$.
  \item For $f\ge0$ and $g\ge0$, $f\ls g$ or $g\gt f$ denotes $f\le Cg$ for a generic constant $C>0$ independent of $\ve>0$,
  and $f\approx g$ means $f\ls g$ and $g\ls f$.
  \item Denote $\|u(t)\|=\|u(t,\cdot)\|$ with norm $\|\cdot\|=\|\cdot\|_{L^2(\cK)},\|\cdot\|_{L^2(\cK_R)},\|\cdot\|_{H^k(\cK)}$.
  \item $|\bar\p f|=(|\bar\p_1f|^2+|\bar\p_2f|^2)^{1/2}$.
  \item $\cW_{\mu,\nu}(t,x)=\w{t+|x|}^\mu(\min\{\w{x},\w{t-|x|}\})^\nu$ for $\mu, \nu\in\Bbb R$.
  \item $\ds|Z^{\le j}f|=\big(\sum_{0\le|a|\le j}|Z^af|^2\big)^\frac12$ and $Z^{\le1}f=(f,Zf)$.
  \item $\delta_0=10^{-3}$.
\end{itemize}

\section{Preliminaries}\label{sect2}

\subsection{Several useful lemmas}

In this subsection, some lemmas on the local energy decay estimate, the elliptic estimate, Gronwall's lemma and
the Sobolev embedding are listed.

\begin{lemma}
Suppose that the obstacle $\cO$ is star-shaped and $\cK=\R^2\setminus\cO$.
Let $w$ be the solution of the IBVP
\begin{equation}\label{loc:decay:eqn}
\left\{
\begin{aligned}
&\Box w=G(t,x),\qquad (t,x)\in[0,\infty)\times\cK,\\
&w=0,\qquad\qquad\quad (t,x)\in[0,\infty)\times\p\cK,\\
&(w,\p_tw)(0,x)=(w_0,w_1)(x),\quad x\in\cK,
\end{aligned}
\right.
\end{equation}
where $\supp_x(w_0(x),w_1(x),G(t,x))\subset\{x: |x|\le R_1\}$ for $R_1>1$.
Then for fixed constant $R>1$, $m\in\N$ and $\eta\in(0,1]$, it holds
\begin{equation}\label{loc:decay}
\begin{split}
\sum_{|a|\le m}\|\p^aw\|_{L^2(\cK_R)}
\ls\w{t}^{-1}\ln^{-2}(2+t)(\|w_0\|_{H^{m}(\cK)}+\|w_1\|_{H^{m-1}(\cK)})\\
+\w{t}^{-\eta}\sum_{|a|\le m-1}\sup_{0\le s\le t}\w{s}^{\eta}\|\p^aG(s)\|_{L^2(\cK)}.
\end{split}
\end{equation}
\end{lemma}
\begin{proof}
One can see the end of Section 3 in \cite{Vainberg75} or Lemma 3.2 of \cite{Kubo14}.
\end{proof}

\begin{lemma}[Lemma 3.1 of \cite{Kubo14}]
Let $w\in H^j(\cK)$ and $\ds w|_{\p\cK}=0$ with integer $j\ge2$.
Then for any fixed constant $R>1$ and multi-index $a\in\N_0^2$ with $2\le|a|\le j$, one has
\begin{equation}\label{ellip}
\|\p_x^aw\|_{L^2(\cK)}\ls\|\Delta w\|_{H^{|a|-2}(\cK)}+\|w\|_{H^{|a|-1}(\cK_{R+1})}.
\end{equation}
\end{lemma}

\begin{lemma}\label{lem-Gronwall}
For any positive function $f(t),g(t),h(t)\ge0$ satisfying
\begin{equation*}
f(t)\le g(t)+\int_0^tf(s)h(s)ds,
\end{equation*}
then it holds
\begin{equation*}
f(t)\le e^{\int_0^th(s)ds}\sup_{0\le s\le t}g(s).
\end{equation*}
\end{lemma}
\begin{proof}
The proof follows from Gronwall's inequality and direct computations.
\end{proof}

\begin{lemma}[Lemma 3.6 of \cite{Kubo14}]
Given any function $f(x)\in C_0^2(\overline{\cK})$, one has that for all $x\in\cK$,
\begin{equation}\label{Sobo:ineq}
\w{x}^{1/2}|f(x)|\ls\sum_{|a|\le2}\|Z^af\|_{L^2(\cK)}.
\end{equation}
\end{lemma}

\subsection{Some pointwise estimates of Cauchy problem}

In order to prove Theorem \ref{thm1}, the following spacetime weighted estimates for the solution of 2-D linear wave equation
will be crucial.

\begin{lemma}\label{lem:pw:div}
Let $w$ be the solution of the Cauchy problem
\begin{equation*}
\left\{
\begin{aligned}
&\Box w=\sum_{\alpha=0}^2\p_\alpha H^\alpha(t,x),\qquad(t,x)\in[0,\infty)\times\R^2,\\
&(w,\p_tw)(0,x)=(0,0),\qquad x\in\R^2.
\end{aligned}
\right.
\end{equation*}
Then for $H^0(t,x)\equiv0$ and $\eta\in(0,1/2)$, one has
\begin{equation}\label{pw:ivp:div}
\w{x}^{1/2}\w{t-|x|}^{1-\eta}|w|\ls\ln(2+t+|x|)\sum_{i=1}^2\sum_{|a|\le1}\sup_{(s,y)\in[0,t]\times\R^2}\w{y}^{1/2}
\cW_{1,1}(s,y)|Z^aH^i(s,y)|,
\end{equation}
where $\ds\cW_{\mu,\nu}(t,x)=\w{t+|x|}^\mu(\min\{\w{x},\w{t-|x|}\})^\nu$;
For general $H^0(t,x)\not\equiv 0$ and $\eta\in(0,1/2)$, it holds that
\begin{equation}\label{pw:ivp:div'}
\w{t+|x|}^{1/2}\w{t-|x|}^{1/2-\eta}|w|\ls\ln(2+t+|x|)\sum_{\alpha=0}^2\sum_{|a|\le1}\sup_{(s,y)\in[0,t]\times\R^2}\w{y}^{1/2}
\cW_{1,1}(s,y)|Z^aH^{\alpha}(s,y)|.
\end{equation}
\end{lemma}
\begin{proof}
Its proof will be postponed to Section \ref{sect2.5}.
\end{proof}

As a corollary of Lemma \ref{lem:pw:div}, we have
\begin{corollary}
Let $w$ be the solution of the Cauchy problem
\begin{equation*}
\left\{
\begin{aligned}
&\Box w=F(t,x),\qquad(t,x)\in[0,\infty)\times\R^2,\\
&(w,\p_tw)(0,x)=(0,0),\qquad x\in\R^2.
\end{aligned}
\right.
\end{equation*}
Then for $\eta\in(0,1/2)$, it holds that
\begin{equation}\label{dtpw:ivp}
\w{x}^{1/2}\w{t-|x|}^{1-\eta}|\p_t^2w|\ls\ln(2+t+|x|)\sum_{|a|\le2}\sup_{(s,y)\in[0,t]\times\R^2}\w{y}^{1/2}\cW_{1,1}(s,y)|Z^aF(s,y)|.
\end{equation}
\end{corollary}
\begin{proof}
Note that
\begin{equation}\label{YHC-6}
\begin{split}
w(t,x)&=\int_0^t\frac{\sin(t-s)|\nabla|}{|\nabla|}F(s)ds,\\
\p^2_tw(t,x)&=F(t,x)+\Delta w(t,x)=F(t,x)+\int_0^t\frac{\sin(t-s)|\nabla|}{|\nabla|}\Delta F(s)ds.
\end{split}
\end{equation}
In addition, \eqref{pw:ivp:div} can be written as
\begin{equation}\label{YHC-7}
\begin{split}
&\w{x}^{1/2}\w{t-|x|}^{1-\eta}\Big|\int_0^t\frac{\sin(t-s)|\nabla|}{|\nabla|}\nabla G(s)ds\Big|\\
&\ls\ln(2+t+|x|)\sum_{|a|\le1}\sup_{(s,y)\in[0,t]\times\R^2}\w{y}^{1/2}\cW_{1,1}(s,y)|Z^aG(s,y)|.
\end{split}
\end{equation}
Combining \eqref{YHC-7} with \eqref{YHC-6} yields \eqref{dtpw:ivp}.
\end{proof}

\begin{lemma}\label{lem:pw:initial}
Let $w$ be the solution of the Cauchy problem
\begin{equation*}
\left\{
\begin{aligned}
&\Box w=H(t,x),\qquad(t,x)\in[0,\infty)\times\R^2,\\
&(w,\p_tw)(0,x)=(w_0,w_1),\quad x\in\R^2.
\end{aligned}
\right.
\end{equation*}
Then for $\mu,\nu\in(0,1/2)$, we have
\begin{align}
&\w{t+|x|}^{1/2}\w{t-|x|}^\mu|w|\ls\cA_{3,0}[w_0,w_1]+\sup_{(s,y)\in[0,t]\times\R^2}\w{y}^{1/2}\cW_{1+\mu,1+\nu}(s,y)|H(s,y)|,\label{pw:ivp1}\\
&\w{t+|x|}^{1/2}\w{t-|x|}^{1/2}|w|\ls\cA_{3,0}[w_0,w_1]+\sup_{(s,y)\in[0,t]\times\R^2}\w{y}^{1/2}\cW_{3/2+\mu,1+\nu}(s,y)|H(s,y)|,\label{pw:ivp2}\\
&\w{t+|x|}^{1/2}\w{t-|x|}^{1/2}|w|\ls\cA_{3,0}[w_0,w_1]+\ln(2+|t-|x||)\sup_{(s,y)\in[0,t]\times\R^2}\w{y}^{1/2}\cW_{3/2,1+\nu}(s,y)|H(s,y)|\label{pw:ivp3}
\end{align}
and
\begin{equation}\label{dpw:ivp}
\begin{split}
&\w{x}^{1/2}\w{t-|x|}^{1+\mu}|\p w|\ls\cA_{4,1}[w_0,w_1]+\sum_{|a|\le1}\sup_{(s,y)\in[0,t]\times\R^2}\w{y}^{1/2}\cW_{1+\mu+\nu,1}(s,y)|Z^aH(s,y)|,
\end{split}
\end{equation}
where $\ds\cA_{\kappa,s}[f,g]=\sum_{\tilde\Gamma\in\{\p_1,\p_2,\Omega\}}
(\sum_{|a|\le s+1}\|\w{z}^\kappa\tilde\Gamma^af(z)\|_{L^\infty}+\sum_{|a|\le s}\|\w{z}^\kappa\tilde\Gamma^ag(z)\|_{L^\infty})$.
\end{lemma}
\begin{proof}
See (3.4), (3.5), (3.7) of \cite{HouYinYuan24} (or (4.2), (4.3), (4.4) and (A.13) of \cite{Kubo19}) for \eqref{pw:ivp1}, \eqref{pw:ivp2}, \eqref{dpw:ivp}, respectively.
The estimate \eqref{pw:ivp3} can be obtained from (4.2) and (A.13) of \cite{Kubo19}.
\end{proof}

\begin{lemma}[Lemma 2.8 of \cite{HouYinYuan26b}]\label{lem:pw:compact}
Let $w$ be the solution of the Cauchy problem
\begin{equation*}
\left\{
\begin{aligned}
&\Box w=H(t,x),\qquad\qquad(t,x)\in[0,\infty)\times\R^2,\\
&(w,\p_tw)(0,x)=(0,0),\quad x\in\R^2,
\end{aligned}
\right.
\end{equation*}
where $\supp_x H(t,x)\subset\{x: |x|\le R\}$.
Then we have that for $\mu,\nu>0$ and $\mu+\nu<1/2$,
\begin{equation}\label{pw:compact}
\w{t+|x|}^{1/2}\w{t-|x|}^{\mu}|w|\ls\sum_{|a|\le2}\sup_{(s,y)\in[0,t]\times\R^2}\w{s}^{1/2+\mu+\nu}|\p^aH(s,y)|
\end{equation}
and
\begin{equation}\label{dpw:compact}
\w{x}^{1/2}\w{t-|x|}|\p w|
\ls\sup_{(s,y)\in[0,t]\times\R^2}\w{s}|(1-\Delta)^3H(s,y)|.
\end{equation}
\end{lemma}

\subsection{Null conditions}\label{sect2-3}

Set
\begin{equation}\label{null:form:def}
\begin{split}
Q_0(f,g)&=\p_tf\p_tg-\sum_{i=1}^2\p_if\p_ig,\\
Q_{\mu\nu}(f,g)&=\p_{\mu}f\p_{\nu}g-\p_{\nu}f\p_{\mu}g,\quad \mu,\nu=0,1,2.
\end{split}
\end{equation}

\begin{lemma}\label{lem:null:structure}
Suppose that the constants $N^{\alpha\beta\gamma}$ satisfy
\begin{equation*}
\sum_{\alpha,\beta,\gamma=0}^2N^{\alpha\beta\gamma}\xi_\alpha\xi_\beta\xi_\gamma\equiv0
\quad\text{for $\xi=(\xi_0,\xi_1,\xi_2)\in\{\pm1\}\times\SS$.}
\end{equation*}
Then for smooth functions $f,g$ and $h$, it holds that
\begin{equation}\label{null:structure}
\begin{split}
\Big|\sum_{\alpha,\beta,\gamma=0}^2N^{\alpha\beta\gamma}
\p^2_{\alpha\beta}f\p_{\gamma}g\Big|&\ls|\bar\p\p f||\p g|+|\p^2f||\bar\p g|,\\
\Big|\sum_{\alpha,\beta,\gamma=0}^2N^{\alpha\beta\gamma}
\p_{\alpha}f\p_{\beta}g\p_{\gamma}h\Big|&\ls|\bar\p f||\p g||\p h|+|\p f||\bar\p g||\p h|+|\p f||\p g||\bar\p h|
\end{split}
\end{equation}
and
\begin{equation}\label{Q12:null:structure}
\begin{split}
Q_{12}(f,g)&=\frac{1}{r}(\p_rf\p_\theta g-\p_\theta f\p_rg)=\sum_{i=1,2}\p_i(\frac{x_if\Omega g}{|x|^2})
-\p_\theta(\frac{f\p_rg}{|x|}),\quad x\neq0,\\
|Q_{12}(f,g)|&\ls(1+|x|)^{-1}(|Zf||\p g|+|\p f||Zg|).
\end{split}
\end{equation}
\end{lemma}
\begin{proof}
Since the proof of \eqref{null:structure} is rather analogous to that in Section 9.1 of \cite{Alinhac:book}
and \cite[Lemma 2.2]{HouYin20jde}, we omit it here.
In addition, \eqref{Q12:null:structure} comes from (2.5) of \cite{HouYinYuan24} directly.
\end{proof}

\begin{lemma}\label{lem:eqn:high}
Let $u$ be the solution of \eqref{QWE} and $Q(\p u,\p^2u)$ be defined by \eqref{nonlinear:def}.
Suppose that the symmetric conditions \eqref{sym:condition} and the first null condition \eqref{1st:null:condition} hold.
Then for any multi-index $a$, one has
\begin{equation}\label{eqn:high}
\begin{split}
&\Box Z^au=Z^aQ(\p u,\p^2u)\\
&=\sum_{\alpha,\beta,\gamma=0}^2\sum_{b+c\le a}Q^{\alpha\beta\gamma}_{abc}\p^2_{\alpha\beta}Z^bu\p_{\gamma}Z^cu
+\sum_{\alpha,\beta,\gamma,\delta=0}^2\sum_{b+c+d\le a}Q^{\alpha\beta\gamma\delta}_{abcd}\p^2_{\alpha\beta}Z^bu\p_{\gamma}Z^cu\p_{\delta}Z^du,
\end{split}
\end{equation}
where $Q_{abc}^{\alpha\beta\gamma}$ and $Q_{abcd}^{\alpha\beta\gamma\delta}$ are constants, in particular $Q_{aa0}^{\alpha\beta\gamma}=Q^{\alpha\beta\gamma}$
and $Q_{aa00}^{\alpha\beta\gamma\delta}=Q^{\alpha\beta\gamma\delta}$ hold.
Furthermore, for any $(\xi_0,\xi_1,\xi_2)\in\{\pm1\}\times\SS$, one has
\begin{equation}\label{null:high}
\sum_{\alpha,\beta,\gamma=0}^2Q_{abc}^{\alpha\beta\gamma}\xi_\alpha\xi_\beta\xi_\gamma\equiv0.
\end{equation}
\end{lemma}
\begin{proof}
See Lemma 6.6.5 of \cite{Hormander97book}.
\end{proof}

\begin{lemma}[Structures of $Q_{0i}$ null forms, Lemma 2.3 of \cite{HouYinYuan24}]\label{lem:Q0i:structure}
For any functions $f,g$, and $\alpha=0,1,2$, we have that for $x\neq0$, $i=1,2$,
\begin{equation}\label{Q0i:structure1}
\begin{split}
Q_{0i}(\p_{\alpha}f,g)&=\sum_{\mu=0}^2\Big\{C_{i\alpha}^{12,\mu}Q_{12}(\p_{\mu}f,g)
+C_{i\alpha}^{0,\mu}(\omega)Q_0(\p_{\mu}f,g)+C_{i\alpha}^{\mu}(\omega)\Box f\p_{\mu}g\Big\}\\
&\quad+\sum_{\mu,\nu=0}^2\frac{C_{i\alpha}^{\mu\nu,1}(\omega)\Omega\p_{\mu}f\p_{\nu}g
+C_{i\alpha}^{\mu\nu,2}(\omega)\Omega g\p^2_{\mu\nu}f}{|x|},
\end{split}
\end{equation}
where $C_{i\alpha}^{12,\mu}$ are constants, $C_{i\alpha}^{0,\mu}(\omega)$,
$C_{i\alpha}^{\mu}(\omega)$, $C_{i\alpha}^{\mu\nu,1}(\omega)$ and $C_{i\alpha}^{\mu\nu,2}(\omega)$
are polynomials of $\omega$.
Furthermore, it holds that
\begin{equation}\label{Q0i:structure2}
\begin{split}
Q_{0i}(\p_{\alpha}f,g)&=\sum_{\mu=0}^2\Big\{C_{i\alpha}^{12,\mu}Q_{12}(\p_{\mu}f,g)
+C_{i\alpha}^{0,\mu}(\omega)Q_0(\p_{\mu}f,g)+C_{i\alpha}^{\mu}(\omega)\Box f\p_{\mu}g\Big\}\\
&\quad+\sum_{\mu=0}^2\p_{\mu}(\frac{C(\omega)g\Omega\p f}{|x|})+\p_\theta(\frac{C(\omega)g\p^2f}{|x|})
+\frac{C(\omega)g\Omega\p f}{|x|^2},
\end{split}
\end{equation}
where $C(\omega)$ stands for the generic polynomial of $\omega$.
\end{lemma}

\begin{lemma}\label{lem:eqn:null:div}
Let $u$ be the solution of \eqref{QWE} and $Q(\p u,\p^2u)$ be defined by \eqref{nonlinear:def}.
Suppose that the symmetric conditions \eqref{sym:condition} and the first null condition \eqref{1st:null:condition} hold.
Then for any multi-index $a$, one has
\begin{equation}\label{eqn:null}
\begin{split}
&\Box Z^au=\sum_{b+c\le a}\Big\{\sum_{\mu,\nu=0}^2|x|^{-1}(C(\omega)\p^2_{\mu\nu}Z^bu\Omega Z^cu
+C(\omega)\Omega\p_{\mu}Z^bu\p_{\nu}Z^cu)\\
&\quad +\sum_{\mu=0}^2\Big[C_{abc}^{12,\mu}Q_{12}(\p_{\mu}Z^bu,Z^cu)
+C_{abc}^{0,\mu}(\omega)Q_0(\p_{\mu}Z^bu,Z^cu)+C(\omega)\Box Z^bu\p_{\mu}Z^cu\Big]\Big\}\\
&\quad +Z^a\big[\sum_{\alpha,\beta,\gamma,\delta=0}^2Q^{\alpha\beta\gamma\delta}\p^2_{\alpha\beta}u\p_{\gamma}u\p_{\delta}u\big],
\end{split}
\end{equation}
where $C_{abc}^{12,\mu}$ are constants and $C_{abc}^{0,\mu}(\omega)$ are polynomials of $\omega$.
Moreover, one has
\begin{equation}\label{eqn:null:div}
\begin{split}
&\Box Z^au=\sum_{b+c\le a}\Big\{\p_\theta(\frac{C(\omega)Z^cu\p^2Z^bu}{|x|})
+\sum_{\mu=0}^2\p_{\mu}(\frac{C(\omega)Z^cuZ\p Z^bu}{|x|})+\frac{C(\omega)Z^cuZ\p Z^bu}{|x|^2}\\
&\quad +\sum_{\mu=0}^2[C_{abc}^{0,\mu}(\omega)Q_0(\p_{\mu}Z^bu,Z^cu)+C(\omega)\p_{\mu}Z^cu\Box Z^bu]\Big\}
+Z^a\big[\sum_{\alpha,\beta,\gamma,\delta=0}^2Q^{\alpha\beta\gamma\delta}\p^2_{\alpha\beta}u\p_{\gamma}u\p_{\delta}u\big].
\end{split}
\end{equation}
\end{lemma}
\begin{proof}
Note that the null condition \eqref{null:high} ensures (see \cite{Klainerman} or \cite[Page 79]{MetcalfeSogge05})
\begin{equation}\label{eqn:null:pf}
\sum_{\alpha,\beta,\gamma=0}^2Q^{\alpha\beta\gamma}_{abc}\p^2_{\alpha\beta}Z^bu\p_{\gamma}Z^cu
=\sum_{\mu=0}^2[C_{abc}^{0,\mu}Q_0(\p_{\mu}Z^bu,Z^cu)+\sum_{\alpha,\beta=0}^2C_{abc}^{\alpha\beta,\mu}Q_{\alpha\beta}(\p_{\mu}Z^bu,Z^cu)],
\end{equation}
where $C_{abc}^{0,\mu}$ and $C_{abc}^{\alpha\beta,\mu}$ are constants.
The equations \eqref{eqn:null} and \eqref{eqn:null:div} can be obtained by \eqref{eqn:high}, \eqref{eqn:null:pf} with \eqref{Q0i:structure1} and by \eqref{Q12:null:structure}, \eqref{Q0i:structure2}, respectively.
\end{proof}

\section{Proof of Lemma \ref{lem:pw:div}}\label{sect2.5}
In order to prove Lemma \ref{lem:pw:div}, we first establish the following basic estimate.

\begin{lemma}\label{YHC-06}
For $\mu\in(0,1/2)$, it holds that
\begin{equation}\label{app:ivp}
\w{t+|x|}^{1/2}\w{t-|x|}^{\mu}|w(t,x)|\ls\ln(2+t+|x|)\sum_{|a|\le1}\sup_{(s,y)\in[0,t]\times\R^2}\w{y}^{1/2}\cW_{1,1}(s,y)|Z^aH(s,y)|,
\end{equation}
provided that $w(t,x)$ is the solution of
\begin{equation*}
\left\{
\begin{aligned}
&\Box w=\p_tH(t,x),\qquad(t,x)\in[0,\infty)\times\R^2,\\
&(w,\p_tw)(0,x)=(0,0),\qquad x\in\R^2.
\end{aligned}
\right.
\end{equation*}
\end{lemma}
\begin{proof}
By Poisson's formula, one has
\begin{equation}\label{app:ivp:pf1}
w(t,x)=\frac{1}{2\pi}\int_{|x-y|<t}\int_{0}^{t-|x-y|}\frac{\p_sH(s,y)dsdy}{\sqrt{(t-s)^2-|x-y|^2}}.
\end{equation}
Set $\rho=|x-y|$, $\tau=t-s-\rho$ and $K(\rho,\tau)=\frac{1}{\sqrt{\tau(2\rho+\tau)}}$.
We decompose  \eqref{app:ivp:pf1} into the following parts
\begin{equation}\label{app:ivp:pf2}
\begin{split}
&2\pi w(t,x)=\int_{\rho<t}\int_{0}^{t-\rho}(1-\chi_{[1/2,1]}(\tau))K(\rho,\tau)\p_sH(s,y)dsdy\\
&\quad+\int_{\rho<t}\int_{0}^{t-\rho}\chi_{[1/2,1]}(\tau)K(\rho,\tau)\p_sH(s,y)dsdy\\
&=\int_{\rho<t}\int_{0}^{t-\rho}(1-\chi_{[1/2,1]}(\tau))K(\rho,\tau)\p_sH(s,y)dsdy\\
&\quad-\int_{\rho<t}\chi_{[1/2,1]}(t-\rho)K(\rho,t-\rho)H(0,y)dy
+\int_{\rho<t}\int_{0}^{t-\rho}\p_{\tau}(\chi_{[1/2,1]}(\tau)K(\rho,\tau))H(s,y)dsdy.
\end{split}
\end{equation}
On the other hand, it is easy to know
$|\p_{\tau}K(\rho,\tau)|=\frac{\rho+\tau}{\tau(2\rho+\tau)}K(\rho,\tau)\le\frac{K(\rho,\tau)}{\tau}$.
Then we have
\begin{equation}\label{app:ivp:pf3}
|(1-\chi_{[1/2,1]}(\tau))K(\rho,\tau)|+|\p_\tau(\chi_{[1/2,1]}(\tau)K(\rho,\tau))|\ls\frac{K(\rho,\tau)}{1+\tau}.
\end{equation}
Substituting \eqref{app:ivp:pf3} into \eqref{app:ivp:pf2} yields
\begin{equation}\label{app:ivp:pf4}
\begin{split}
&|w(t,x)|\ls M_t\big(X(t,x)+\int_0^tk(\tau)Y(t-\tau,x)d\tau\big),
\end{split}
\end{equation}
where
\begin{equation*}\label{app:ivp:pf4-0}
\begin{split}
&M_t=\sum_{|a|\le1}\sup_{(s,y)\in[0,t]\times\R^2}W(s,y)|Z^aH(s,y)|,\quad
X(t,x)=\int_{|x-y|<t}\frac{(1+|y|)^{-5/2}dy}{\sqrt{t^2-|x-y|^2}},\\
&k(\tau)=\frac{1}{\sqrt{\tau}(1+\tau)},\quad Y(t,x)=\int_{|x-y|<t}\frac{|x-y|^{-1/2}dy}{W(t-|x-y|,y)}\quad
\text{and  $W(s,y)=\w{y}^{1/2}\cW_{1,1}(s,y)$}.
\end{split}
\end{equation*}
Next we estimate $X(t,x)$ and $Y(t,x)$. For notational convenience,
set $r=|x|$, $p=t+r$, $q=|t-r|$, $\lambda=|y|$ and $\zeta=t-r$.

\vskip 0.2 true cm

{\bf $\bullet$ Estimate of $X(t,x)$}
\vskip 0.1 true cm

We assert that
\begin{equation}\label{app:ivp:pf5}
X(t,x)\ls\w{t+r}^{-1/2}\w{t-r}^{-1/2}.
\end{equation}
If $p\le8$, then \eqref{app:ivp:pf5} can be obtained by
\begin{equation}\label{app:ivp:pf6}
\int_{|x-y|<t}\frac{dy}{\sqrt{t^2-|x-y|^2}}=2\pi t.
\end{equation}
In what follows, $p\ge8$ is assumed.  Denote
\begin{equation}\label{app:ivp:pf7}
\begin{split}
&X_R(t,x)=\int_{|x-y|<t,R/2\le|y|\le R}\frac{dy}{\sqrt{t^2-|x-y|^2}},\\
&\tilde X_R(t,x)=\int_{|x-y|<t,|y|\le R}\frac{dy}{\sqrt{t^2-|x-y|^2}}.
\end{split}
\end{equation}
If $R=2^n\ge p/4$, then it follows from \eqref{app:ivp:pf6} that
\begin{equation}\label{app:ivp:pf8}
\sum_{R=2^n\ge p/4}R^{-5/2}X_R(t,x)\ls p^{-5/2}t\ls p^{-3/2}.
\end{equation}
Note that $r=|x|\le|x-y|+|y|<t+R$, otherwise, $\tilde X_R(t,x)$ vanishes.

If $R=2^n\in[1,p/4]$, then $t>\frac{p-R}{2}\ge\frac{3p}{8}\ge\frac{3R}{2}$ holds.
By the rotation of coordinates, one can assume $x=(r,0)$, $y=(z,h)$ and $z_\pm(h)=r\pm\sqrt{t^2-h^2}$.
Note that $|y|\le R$ holds in the integral domain of $\tilde X_R(t,x)$. This implies that $|h|\le R\le 2t/3$, $|z|\le R\le p/4$, $z_+(h)-z\ge\frac{\sqrt{5}-2}{8}p$ and
$t^2-|x-y|^2=(z-z_-(h))(z_+(h)-z)$.
Therefore, one can achieve that for $R=2^n\in[1,p/4]$,
\begin{equation}\label{app:ivp:pf9}
X_R(t,x)\ls\tilde X_R(t,x)\ls p^{-1/2}\int_{-R}^R\int_{-R\le z\le R,z\ge z_-(h)}(z-z_-(h))^{-1/2}dzdh\ls R^{3/2}p^{-1/2}.
\end{equation}
For $q=|t-r|\le2$, \eqref{app:ivp:pf9} leads to
\begin{equation}\label{app:ivp:pf10}
\tilde X_1(t,x)+\sum_{n\ge0}\sum_{R=2^n\le p/4}R^{-5/2}X_R(t,x)\ls p^{-1/2}.
\end{equation}
For $q\ge2$ and $t\ge r=|x|$, when $R=2^n\in[q/4,p/4]$, \eqref{app:ivp:pf9} gives
\begin{equation}\label{app:ivp:pf11}
\sum_{n\ge0}\sum_{R=2^n\in[q/4,p/4]}R^{-5/2}X_R(t,x)\ls\sum_{n\ge0}\sum_{R=2^n\ge q/4}R^{-1}p^{-1/2}\ls p^{-1/2}(1+q)^{-1};
\end{equation}
when $R=2^n\le q/4=|t-r|/4$, it is easy to get $t^2-|x-y|^2\ge t^2-(r+R)^2\gt pq$ and
\begin{equation}\label{app:ivp:pf12}
\begin{split}
&\tilde X_R(t,x)\ls R^2(pq)^{-1/2},\\
&\tilde X_1(t,x)+\sum_{n\ge0}\sum_{R=2^n\le q/4}R^{-5/2}X_R(t,x)\ls(pq)^{-1/2}.
\end{split}
\end{equation}
For $q\ge2$ and $t\le r$, the analogous estimate for  \eqref{app:ivp:pf11} can be obtained.
In this case, one has $R\ge q=r-t$ in the integral domain of $\tilde X_R(t,x)$ in \eqref{app:ivp:pf7}, otherwise,
$|y|\le R\le|x|-t\le|x-y|+|y|-t$ holds and $\tilde X_R(t,x)$ vanishes.

Therefore, the estimate \eqref{app:ivp:pf5} can be obtained by combining \eqref{app:ivp:pf7}-\eqref{app:ivp:pf12} with the dyadic decomposition in space $\Bbb R^2$.

\vskip 0.2 true cm
{\bf $\bullet$ Estimate of $Y(t,x)$}
\vskip 0.1 true cm

We assert that
\begin{equation}\label{app:ivp:pf13}
Y(t,x)\ls\frac{\ln(2+t+|x|)\ln(2+|t-|x||)}{\w{t+|x|}^{1/2}\w{t-|x|}^{1/2}}.
\end{equation}
For fixed $x\neq0$, the Jacobian determinant of the coordinate transformation $y=(y_1,y_2)\rightarrow(\lambda,\rho)$ is
\begin{equation}\label{app:ivp:pf14}
\Big|\frac{\p(\lambda,\rho)}{\p(y_1,y_2)}\Big|=\frac{\sqrt{(\lambda+\rho+r)(\lambda+\rho-r)(r+\lambda-\rho)(r-\lambda+\rho)}}{2\lambda\rho},
\end{equation}
where $\lambda=|y|$, $\rho=|x-y|$ and $r=|x|$.

Set
$$\alpha=\lambda+t-\rho\quad \text{and}\quad \beta=\lambda+\rho-t.$$
From $\rho<t$, one has
\begin{equation}\label{app:ivp:pf15}
\alpha\in(q,p),\quad \beta\in(-\zeta,\alpha),\quad \lambda=\frac{1}{2}(\alpha+\beta),
\quad \rho=\frac{1}{2}(2t-\alpha+\beta),\quad \Big|\frac{\p(\lambda,\rho)}{\p(\alpha,\beta)}\Big|=1/2.
\end{equation}
Substituting \eqref{app:ivp:pf14} and \eqref{app:ivp:pf15} into $Y(t,x)$ together with the facts of $\rho\le(p+\beta)/2$ and $\lambda(1+\lambda)^{-1/2}\ls\lambda^{1/2}$ yields
\begin{equation}\label{app:ivp:pf16}
\begin{split}
Y(t,x)&\ls\int_q^p\frac{Z(\alpha,\zeta)d\alpha}{(1+\alpha)\sqrt{(\alpha-\zeta)(p-\alpha)}}\quad\text{with
$Z(\alpha,\zeta)=\int_{-\zeta}^{\alpha}\frac{\lambda^{1/2}d\beta}{\sqrt{\beta+\zeta}\min\{1+\lambda,1+|\beta|\}}$},
\end{split}
\end{equation}
where $\zeta=t-r$.

We now deal with $Z(\alpha,\zeta)$. Note that
\begin{equation}\label{app:ivp:pf17}
\begin{split}
Z(\alpha,\zeta)&\ls Z_1(\alpha,\zeta)+\alpha^{1/2}Z_2(\alpha,\zeta),
\end{split}
\end{equation}
where $Z_1(\alpha,\zeta)=\int_{-\zeta}^{\alpha}\frac{\lambda^{1/2}d\beta}{\sqrt{\beta+\zeta}(1+\lambda)}$
and $Z_2(\alpha,\zeta)=\int_{-\zeta}^{\alpha}\frac{d\beta}{\sqrt{\beta+\zeta}(1+|\beta|)}$.

For $Z_1(\alpha,\zeta)$, set $z=\beta+\zeta$, one can arrive at
\begin{equation}\label{app:ivp:pf18}
Z_1(\alpha,\zeta)\ls\int_{0}^{\alpha+\zeta}\frac{dz}{\sqrt{z(1+\alpha-\zeta+z)}}\ls\ln(2+\alpha),
\end{equation}
where the last inequality comes from the transformation of $z=(1+\alpha-\zeta)\sinh^2\psi$.

For $Z_2(\alpha,\zeta)$, we now show
\begin{equation}\label{app:ivp:pf19}
Z_2(\alpha,\zeta)\ls\ln(2+q)(1+q)^{-1/2}\quad\text{with $q=|\zeta|$.}
\end{equation}
Indeed, for $\zeta=-q\le0$, one has
\begin{equation}\label{app:ivp:pf20}
Z_2(\alpha,\zeta)\ls(1+q)^{-1}\int_{q}^{1+2q}\frac{d\beta}{\sqrt{\beta-q}}+\int_{1+2q}^{\infty}\frac{d\beta}{(1+\beta)^{3/2}}\ls(1+q)^{-1/2};
\end{equation}
for $\zeta=q>0$ and $q\ge1$, analogously to \eqref{app:ivp:pf20}, we have
\begin{equation}\label{app:ivp:pf21}
\begin{split}
Z_2(\alpha,\zeta)&\ls q^{-1}\int_{-q}^{-q/2}\frac{d\beta}{\sqrt{\beta+q}}+q^{-1/2}\int_{-q/2}^{q}\frac{d\beta}{1+|\beta|}
+\int_{q}^{\infty}\frac{d\beta}{(1+\beta)^{3/2}}\\
&\ls(1+q)^{-1/2}\ln(2+q);
\end{split}
\end{equation}
when $q\le1$, it holds that
\begin{equation}\label{app:ivp:pf22}
Z_2(\alpha,\zeta)\ls\int_{-q}^{2-q}\frac{d\beta}{\sqrt{\beta+q}}+\int_{2-q}^{\infty}\frac{d\beta}{(1+\beta)^{3/2}}\ls1.
\end{equation}
Collecting \eqref{app:ivp:pf20}-\eqref{app:ivp:pf22} yields \eqref{app:ivp:pf19}.

Substituting \eqref{app:ivp:pf18}, \eqref{app:ivp:pf19} into \eqref{app:ivp:pf17}, we have
\begin{equation}\label{app:ivp:pf23}
Z(\alpha,\zeta)\ls\ln(2+\alpha)+\alpha^{1/2}\ln(2+q)(1+q)^{-1/2}.
\end{equation}
Based on \eqref{app:ivp:pf16} and \eqref{app:ivp:pf23}, we next estimate $Y(t,x)$. Note that
\begin{equation}\label{app:ivp:pf24}
\begin{split}
Y(t,x)&\ls\ln(2+q)(1+q)^{-1/2}\int_q^p\frac{d\alpha}{\sqrt{(1+\alpha)(\alpha-q)(p-\alpha)}}\\
&\quad+\ln(2+p)\int_q^p\frac{d\alpha}{(1+\alpha)\sqrt{(\alpha-q)(p-\alpha)}}.
\end{split}
\end{equation}
Let $\alpha=q+(p-q)\sin^2\psi$. Then the integral in the first line of \eqref{app:ivp:pf24} can be reduced to
\begin{equation}\label{app:ivp:pf25}
\int_q^p\frac{d\alpha}{\sqrt{(1+\alpha)(\alpha-q)(p-\alpha)}}
=2\int_0^{\pi/2}\frac{d\psi}{\sqrt{1+q+(p-q)\sin^2\psi}}.
\end{equation}
If $p-q\le1+q$, then one has
\begin{equation}\label{app:ivp:pf26}
\int_0^{\pi/2}\frac{d\psi}{\sqrt{1+q+(p-q)\sin^2\psi}}\ls(1+q)^{-1/2}\ls(1+p)^{-1/2}.
\end{equation}
If $p-q\ge1+q$, then it follows from $\sin\psi\ge2\psi/\pi$ that
\begin{equation}\label{app:ivp:pf27}
\begin{split}
\int_0^{\pi/2}\frac{d\psi}{\sqrt{1+q+(p-q)\sin^2\psi}}&\ls\int_0^{\pi/2}\frac{d\psi}{\sqrt{1+q+4(p-q)\psi^2/\pi^2}}\\
&\ls(p-q)^{-1/2}{\rm arsinh}\big(\frac{p-q}{1+q}\big)^{1/2}\ls\ln(2+p)(1+p)^{-1/2}.
\end{split}
\end{equation}
Similarly, the integral in the second line of \eqref{app:ivp:pf24} can be treated as follows
\begin{equation}\label{app:ivp:pf28}
\begin{split}
\int_q^p\frac{d\alpha}{(1+\alpha)\sqrt{(\alpha-q)(p-\alpha)}}
&=2\int_0^{\pi/2}\frac{d\psi}{1+q+(p-q)\sin^2\psi}\\
&=2\int_0^{\infty}\frac{dz}{1+q+(1+p)z^2}\ls(1+p)^{-1/2}(1+q)^{-1/2},
\end{split}
\end{equation}
where the transformation $z=\tan\psi$ is used in the last line of \eqref{app:ivp:pf28}.

Substituting \eqref{app:ivp:pf25}-\eqref{app:ivp:pf28} into \eqref{app:ivp:pf24} yields \eqref{app:ivp:pf13} for $x\neq0$.
In addition, in the case of $x=0$,
\eqref{app:ivp:pf13} can be analogously obtained.

We start to prove \eqref{app:ivp}.
It follows from \eqref{app:ivp:pf4}, \eqref{app:ivp:pf5} and \eqref{app:ivp:pf13} that for $\mu\in(0,1/2)$,
\begin{equation}\label{app:ivp:pf29}
|w(t,x)|\ls M_t\Big(\w{t+r}^{-1/2}\w{t-r}^{-1/2}+\ln(2+t+r)\int_0^t\frac{k(\tau)d\tau}{\w{t-\tau+r}^{1/2}\w{t-\tau-r}^{\mu}}\Big).
\end{equation}
For $r=|x|>t$, it is easy to know
\begin{equation*}
t-\tau+r\gt t+r, \quad |t-\tau-r|=r+\tau-t\ge r-t>0.
\end{equation*}
This, together with the fact of $\ds\int_0^{\infty}k(\tau)\ls1$, implies
\begin{equation}\label{app:ivp:pf30}
\int_0^t\frac{k(\tau)d\tau}{\w{t-\tau+r}^{1/2}\w{t-\tau-r}^{\mu}}\ls\w{t+r}^{-1/2}\w{t-r}^{-\mu}.
\end{equation}
For $r=|x|\le t$, analogously to \eqref{app:ivp:pf21}, we have that for $t-r\ge1$,
\begin{equation}\label{app:ivp:pf31}
\begin{split}
&\int_0^{\infty}k(\tau)(1+|t-\tau-r|)^{-\mu}d\tau\\
&\ls\w{t-r}^{-\mu}\int_0^{(t-r)/2}k(\tau)d\tau+\w{t-r}^{-3/2}\int_{(t-r)/2}^{2(t-r)}(1+|t-\tau-r|)^{-\mu}d\tau\\
&\quad+\int_{2(t-r)}^{\infty}\tau^{-3/2-\mu}d\tau\\
&\ls\w{t-r}^{-\mu};
\end{split}
\end{equation}
when $t-r\le1$, \eqref{app:ivp:pf31} follows from $\ds\int_0^{\infty}k(\tau)\ls1$.
Therefore, we can arrive at
\begin{equation}\label{app:ivp:pf32}
\begin{split}
&\int_0^t\frac{k(\tau)d\tau}{\w{t-\tau+r}^{1/2}\w{t-\tau-r}^{\mu}}\\
&\ls\w{t+r}^{-1/2}\int_0^{(1+t+r)/2}k(\tau)(1+|t-\tau-r|)^{-\mu}d\tau\\
&\quad+\w{t+r}^{-3/2}\int_{(1+t+r)/2}^t(1+t-\tau+r)^{-1/2}d\tau\\
&\ls\w{t+r}^{-1/2}\w{t-r}^{-\mu},
\end{split}
\end{equation}
where  \eqref{app:ivp:pf31} is used.
Therefore, \eqref{app:ivp} is proved by \eqref{app:ivp:pf29}, \eqref{app:ivp:pf30} and \eqref{app:ivp:pf32}.
\end{proof}
Based on Lemma \ref{YHC-06}, we now prove Lemma \ref{lem:pw:div}.

\begin{proof}[Proof of Lemma \ref{lem:pw:div}]
At first, we prove \eqref{pw:ivp:div} for $H^0(t,x)\equiv0$.
Choosing $\mu=\eta>0$, $\nu=1-\eta$ and $\kappa=1$ in (B.6) of \cite{Kubo14}, one then has
\begin{equation}\label{pw:ivp:div1}
\w{x}^{1/2}\w{t-|x|}^{1-\eta}|v|\ls\ln(2+t+|x|)\sum_{i=1}^2
\sum_{|a|\le1}\sup_{(s,y)\in[0,t]\times\R^2}\w{y}^{3/2}\w{s+|y|}|Z^aH^i(s,y)|.
\end{equation}
Recall (B.7) of \cite{Kubo14}
\begin{equation}\label{pw:ivp:div2}
\w{x}^{1/2}\w{t-|x|}^{1-\eta}|v|\ls\ln(2+t+|x|)\sum_{i=1}^2
\sum_{|a|\le1}\sup_{(s,y)\in[0,t]\times\R^2}\w{y}^{1/2}\w{s+|y|}\w{s-|y|}|Z^aH^i(s,y)|.
\end{equation}
Combining \eqref{pw:ivp:div1} with \eqref{pw:ivp:div2} yields \eqref{pw:ivp:div}.
It is pointed out that the estimates (B.6) and (B.7) of \cite{Kubo14} are only proved for $H^0(t,x)\equiv0$.

When $H^0(t,x)\not\equiv0$, the estimate \eqref{pw:ivp:div'} can be achieved by \eqref{pw:ivp:div} and \eqref{app:ivp} directly.
\end{proof}

\section{Some key estimates}\label{sect3}

We now recall the Littlewood-Paley decomposition.
Choosing a smooth cutoff function $\psi: \R\rightarrow[0,1]$, which equals 1 on $[-5/4,5/4]$ and vanishes outside $[-8/5,8/5]$,
one sets
\begin{equation*}
\psi_k(\xi)=\psi(|\xi|/2^k)-\psi(|\xi|/2^{k-1}),\quad k\in\Z,\quad\xi\in\R^2.
\end{equation*}
Let $P_k$ be the Littlewood-Paley projection onto frequency $2^k$
\begin{equation*}
\sF(P_kf)(\xi)=\psi_k(\xi)\sF f(\xi),\quad k\in\Z.
\end{equation*}

\subsection{Estimates of the homogeneous initial boundary value problem}

\begin{lemma}\label{YHC-8}
Let $w_{hom}$ be the solution to the homogeneous initial boundary value problem \eqref{hom:eqn},
then for any multi-index $a$, one has that
\begin{equation}\label{hom:pw1}
\begin{split}
|\p Z^aw_{hom}(t,x)|&\ls\ve\w{x}^{-1/2}\w{t-|x|}^{-1},\\
|Z^aw_{hom}(t,x)|&\ls\ve\w{t+|x|}^{-1/2}\w{t-|x|}^{\dl_0-1/2},\\
|\bar\p Z^aw_{hom}(t,x)|&\ls\ve\w{x}^{-1/2}\w{t+|x|}^{-1},
\end{split}
\end{equation}
where $\dl_0=10^{-3}$.
Moreover, for $|x|\le t/4$, it holds that
\begin{equation}\label{hom:pw2}
|\p Z^aw_{hom}(t,x)|\ls\ve\w{t}^{-1}\ln(2+t)^{-2}.
\end{equation}
\end{lemma}
\begin{proof}
It follows from the local energy decay estimate \eqref{loc:decay} that for any multi-index $b$,
\begin{equation}\label{hom:pw:pf1}
\|\p^bw_{hom}(t,x)\|_{L^2(\cK_R)}\ls\ve\w{t}^{-1}\ln(2+t)^{-2},
\end{equation}
where $R>0$ is a large constant.
Thus, the estimates \eqref{hom:pw1} and \eqref{hom:pw2} in the region $|x|\le2$ can be achieved.

In what follows, we focus on the case of $|x|\ge2$ and $t\ge1$.
Note that $w^c(t,x)=\chi_{[1,2]}(x)w_{hom}(t,x)$ is the solution of the Cauchy problem
\begin{equation}\label{hom:pw:pf1-1}
\Box w^c=-[\Delta,\chi_{[1,2]}(x)]w_{hom},\quad(w^c,\p_tw^c)(0,x)=(\ve\chi_{[1,2]}(x)u_0(x),\ve\chi_{[1,2]}(x)u_1(x)).
\end{equation}
The first two inequalities in \eqref{hom:pw1} follow from Lemmas \ref{lem:pw:initial}-\ref{lem:pw:compact} with $\mu=1/2-\delta_0$, $\nu=\delta_0/2$ and \eqref{hom:pw:pf1}.

The third inequality of \eqref{hom:pw1} in the region $|x|\le t/2$ can be obtained by the first inequality.

In the region $|x|\ge t/2$, it is easy to get that (one can also see Section 4.6 of \cite{Kubo13} or Section 7 of \cite{HouYinYuan24})
\begin{equation}\label{hom:pw:pf2}
\begin{split}
&\p_+(r^{1/2}w)(t,r\frac{x}{|x|})-\p_+(r^{1/2}w)(0,(t+r)\frac{x}{|x|})\\
&\quad=\int_0^t\{(r+t-s)^{1/2}\Box w+(r+t-s)^{-3/2}(w/4+\Omega^2w)\}(s,(t+r-s)\frac{x}{|x|})ds,
\end{split}
\end{equation}
where $\bar\p_1=\frac{x_1}{r}\p_+-\frac{x_2}{r^2}\Omega$, $\bar\p_2=\frac{x_2}{r}\p_++\frac{x_1}{r^2}\Omega$
and $\p_{\pm}=\p_t\pm\p_r$.

Set $w=Z^aw_{hom}$ in \eqref{hom:pw:pf2}. Then one obtains
\begin{equation*}
\begin{split}
|\p_+(r^{1/2}Z^aw_{hom})(t,x)|&\ls\ve\w{x}^{-1}+\ve\w{r+t}^{-2}\int_0^t(1+|r+t-2s|)^{\delta_0-1/2}ds\\
&\ls\ve\w{t+|x|}^{-1},
\end{split}
\end{equation*}
where the second inequality in \eqref{hom:pw1} is used.
Thus, the third inequality of \eqref{hom:pw1} is achieved.

Finally, we focus on the proof of \eqref{hom:pw2}.
The solution to problem \eqref{hom:pw:pf1-1} can be decomposed into $w^c=w^c_1+w^c_2$, where $w^c_1$ and $w^c_2$ are the solutions
of the Cauchy problems
\begin{equation*}
\begin{split}
&\Box w^c_1=H=-[\Delta,\chi_{[1,2]}(x)]w_{hom}, \quad (w^c_1,\p_tw^c_1)(0,x)=(0,0),\\
&\Box w^c_2=0,\quad (w^c_2,\p_tw^c_2)(0,x)=(\ve\chi_{[1,2]}(x)u_0(x),\ve\chi_{[1,2]}(x)u_1(x)),
\end{split}
\end{equation*}
respectively.
The estimate of $\p_{\alpha}w^c_2$ ($\alpha=0,1,2$) can be obtained by \eqref{dpw:ivp}.
Note that $\p_tw^c_1$ and $\p_jw^c_1$ ($j=1,2$) are expressed as
\begin{equation*}
\begin{split}
\p_tw^c_1(t,x)&=\int_0^t\big(\cos(t-s)|\nabla|\big)H(s)ds
=\int_0^t\frac{e^{i(t-s)|\nabla|}+e^{-i(t-s)|\nabla|}}{2}H(s)ds\\
&=\frac{1}{2}\int_0^t\int_{\R^2}\sum_{k\in\Z}(K^{0+}_k+K^{0-}_k)(t-s,x-y)(1-\Delta)^3H(s,y)dyds\\
&\text{with~~}K^{0\pm}_k(t,x)=\frac{1}{(2\pi)^2}\int_{\R^2}e^{i(\pm t|\xi|+x\cdot\xi)}\frac{\psi_k(\xi)}{(1+|\xi|^2)^3}d\xi
\end{split}
\end{equation*}
and
\begin{equation}\label{hom:pw:pf3}
\begin{split}
\p_jw^c_1(t,x)&=\int_0^t\big(\sin(t-s)|\nabla|\big)R_jH(s)ds\\
&=\frac{1}{2i}\int_0^t\int_{\R^2}\sum_{k\in\Z}(K^{j+}_k-K^{j-}_k)(t-s,x-y)P_kR_j(1-\Delta)^3H(s,y)dyds\\
&\text{with~~}K^{j\pm}_k(t,x)=\frac{1}{(2\pi)^2}\int_{\R^2}e^{i(\pm t|\xi|+x\cdot\xi)}\frac{\psi_{[k-1,k+1]}(\xi)}{(1+|\xi|^2)^3}d\xi,
\end{split}
\end{equation}
where $R_j=\p_j|\nabla|^{-1}$ ($j=1,2$) is the Riesz transformation and the fact of $\psi_k(\xi)=\psi_{[k-1,k+1]}(\xi)\psi_k(\xi)$
has been used.
In the following, it suffices to deal with the estimates of $\p_jw_1^c$ since the estimate of $\p_tw_1^c$ can be treated similarly.
Recall (4.28) and (4.30) of \cite{HouYinYuan24}
\begin{equation}\label{hom:pw:pf4}
\begin{split}
|P_kR_j(1-\Delta)^3H(s,x)|&\ls2^{2k}(1+2^{k-k_+}|x|)^{-5}\|(1-\Delta)^3H(s,\cdot)\|_{L^\infty(\R^2)},\\
|K^{j\pm}_k(t,x)|&\ls2^{2k-6k_+}(1+2^k|x|)^{-1/2}(1+2^k|t-|x||)^{-2}.
\end{split}
\end{equation}
Substituting \eqref{hom:pw:pf1} and \eqref{hom:pw:pf4} into \eqref{hom:pw:pf3} yields
\begin{equation}\label{hom:pw:pf5}
\begin{split}
&|\p_jw^c_1|\ls\ve\sum_{k\in\Z}2^{4k-6k_+}I^k,\\
&I^k=\int_0^t\int_{\R^2}\frac{dyds}{(1+s)\ln^2(2+s)(1+2^k|t-s-|x-y||)^2(1+2^{k-k_+}|y|)^5}.
\end{split}
\end{equation}
The integral $I^k$ can be separated into the following three parts
\begin{equation}\label{hom:pw:pf6}
\begin{split}
I^k&=I^k_1+I^k_2+I^k_3,\\
I^k_1&=\int_{|x-y|\ge t/2}\int_0^t\frac{dsdy}{(1+s)\ln^2(2+s)(1+2^k|t-s-|x-y||)^2(1+2^{k-k_+}|y|)^5},\\
I^k_2&=\int_{|x-y|\le t/2}\int_0^{\frac{t-|x-y|}{2}}\frac{dsdy}{(1+s)\ln^2(2+s)(1+2^k(t-s-|x-y|))^2(1+2^{k-k_+}|y|)^5},\\
I^k_3&=\int_{|x-y|\le t/2}\int_{\frac{t-|x-y|}{2}}^t\frac{dsdy}{(1+s)\ln^2(2+s)(1+2^k|t-s-|x-y||)^2(1+2^{k-k_+}|y|)^5}.
\end{split}
\end{equation}
For the integral $I^k_1$, one can see that $|y|\ge|x-y|-|x|\ge t/2-t/4=t/4$ and
\begin{equation}\label{hom:pw:pf7}
\begin{split}
I^k_1&\ls\w{t}^{-3/2}2^{3(k_+-k)/2}\int_0^t\int\frac{dsdy}{(1+s)\ln^2(2+s)(1+2^{k-k_+}|y|)^{7/2}}\\
&\ls\w{t}^{-3/2}2^{7(k_+-k)/2}.
\end{split}
\end{equation}
For the integral $I^k_2$, one has that $t-s-|x-y|\ge\frac{t-|x-y|}{2}\ge t/4$ and
\begin{equation}\label{hom:pw:pf8}
\begin{split}
I^k_2&\ls\w{t}^{-3/2}2^{3(k_+-k)/2}\int_0^t\int\frac{dsdy}{(1+s)\ln^2(2+s)(1+2^{k-k_+}|y|)^5}\\
&\ls\w{t}^{-3/2}2^{7(k_+-k)/2}.
\end{split}
\end{equation}
For the integral $I^k_3$, due to $s\ge\frac{t-|x-y|}{2}\ge t/4$, it holds that
\begin{equation}\label{hom:pw:pf9}
\begin{split}
I^k_3&\ls\w{t}^{-1}\ln^{-2}(2+t)\int_0^t\int\frac{dsdy}{(1+2^k|t-s-|x-y||)^2(1+2^{k-k_+}|y|)^5}\\
&\ls\w{t}^{-1}\ln^{-2}(2+t)2^{3(k_+-k)}.
\end{split}
\end{equation}
Collecting \eqref{hom:pw:pf5}-\eqref{hom:pw:pf9} leads to
\begin{equation*}
|\p_jw^c_1|\ls\ve\w{t}^{-1}\ln^{-2}(2+t)\sum_{k\in\Z}2^{k/2-5k_+/2}\ls\ve\w{t}^{-1}\ln^{-2}(2+t).
\end{equation*}
This completes the proof of Lemma \ref{YHC-8}.
\end{proof}

\subsection{Some estimates of Cauchy problem}

\begin{lemma}\label{lem:dt:impr}
Let $w$ be the solution of the linear Cauchy problem
\begin{equation*}
\left\{
\begin{aligned}
&\Box w=F(t,x),\qquad \qquad \qquad (t,x)\in[0,+\infty)\times\R^2,\\
&(w,\p_tw)(0,x)=(0,0),\qquad\quad x\in\R^2,
\end{aligned}
\right.
\end{equation*}
provided that $\supp_yF(s,y)\subset\{y: \frac{9s}{10}\le|y|\le s+M_0\}$ and $M_0>0$.
Then for any constant $R>0$ and $\mu>0$, it holds that
\begin{equation}\label{dt:impr}
\w{t}|\p_tw(t,x)|\ls\sup_{(s,y)\in[0,t]\times\R^2}\w{y}^{1/2}\cW_{1,1+\mu}|(1-\Delta)^2F(s,y)|.
\end{equation}
where $x\in\{x\in\R^2:|x|\le\frac{t}{10}~{\rm or}~|x|\le R\}$.
\end{lemma}
\begin{proof}
Without loss of generality, we can assume $t\ge20M_0+10R$ since \eqref{dt:impr} obviously holds for the case of $t\le20M_0+10R$.
Due to $\{|x|\le R\}\subset\{|x|\le t/10\}$, we only focus on the case of $|x|\le t/10$.
Note that
\begin{equation*}
\begin{split}
\p_tw(t,x)&=\frac{1}{2}\int_0^t\int_{\R^2}\sum_{k\in\Z}(\tilde K^{0+}_k+\tilde K^{0-}_k)(t-s,x-y)(1-\Delta)^2F(s,y)dyds\\
&\text{with~~}\tilde K^{0\pm}_k(t,x)=\frac{1}{(2\pi)^2}\int_{\R^2}e^{i(\pm t|\xi|+x\cdot\xi)}\frac{\psi_k(\xi)}{(1+|\xi|^2)^2}d\xi
\end{split}
\end{equation*}
and
\begin{equation}\label{dt:impr:pf1}
\begin{split}
&|\tilde K^{0\pm}_k(t,x)|\ls2^{2k-4k_+}(1+2^k|x|)^{-1/2}(1+2^k|t-|x||)^{-2},\\
&|\p_tw|\ls\sum_{k\in\Z}2^{2k-4k_+}J^k\sup_{(s,y)\in[0,t]\times\R^2}\w{y}^{1/2}\cW_{1,1+\mu}|(1-\Delta)^2F(s,y)|,\\
\end{split}
\end{equation}
where
$$J^k=\int_0^t\int_{\frac{9s}{10}\le|y|\le s+M_0}\frac{(1+2^k|x-y|)^{-1/2}dyds}{(1+s+|y|)^{3/2}(2M_0+s-|y|)^{1+\mu}(1+2^k|t-s-|x-y||)^2}.$$

The integral $J^k$ can be separated into three parts
\begin{equation}\label{dt:impr:pf2}
\begin{split}
&J^k=J^k_1+J^k_2+J^k_3,\\
&J^k_1=\int_0^{2t/5}\int_{\frac{9s}{10}\le|y|\le s+M_0}\frac{(1+2^k|x-y|)^{-1/2}dyds}{(1+s+|y|)^{3/2}(2M_0+s-|y|)^{1+\mu}(1+2^k|t-s-|x-y||)^2},\\
&J^k_2=\int_{3t/5}^t\int_{\frac{9s}{10}\le|y|\le s+M_0}\frac{(1+2^k|x-y|)^{-1/2}dyds}{(1+s+|y|)^{3/2}(2M_0+s-|y|)^{1+\mu}(1+2^k|t-s-|x-y||)^2},\\
&J^k_3=\int_{2t/5}^{3t/5}\int_{\frac{9s}{10}\le|y|\le s+M_0}\frac{(1+2^k|x-y|)^{-1/2}dyds}{(1+s+|y|)^{3/2}(2M_0+s-|y|)^{1+\mu}(1+2^k|t-s-|x-y||)^2}.
\end{split}
\end{equation}
For the integral $J^k_1$, one sees that $t-s-|x-y|\ge t-s-|x|-|y|\ge t-2s-|x|-M_0\ge \frac{t}{10}-M_0\ge\frac{t}{20}$ holds,
while for the integral $J^k_2$, one has $|t-s-|x-y||=s+|x-y|-t\ge s+|y|-|x|-t\ge t(\frac{57}{50}-\frac{11}{10})=\frac{t}{25}$.
Therefore,
\begin{equation}\label{dt:impr:pf3}
\begin{split}
J^k_1+J^k_2&\ls2^{3(k_+-k)/2}(1+t)^{-3/2}\int_0^t\int_{|y|\le s+M_0}\frac{dyds}{(1+s+|y|)^{3/2}(2M_0+s-|y|)^{1+\mu}}\\
&\ls2^{3(k_+-k)/2}(1+t)^{-1}.
\end{split}
\end{equation}
For the integral $J^k_3$, $|y|\ge\frac{9t}{25}\ge3|x|$ and $|x-y|\ge\frac{2|y|}{3}\ge t/5$ hold.
Thus, we have
\begin{equation}\label{dt:impr:pf4}
\begin{split}
&J^k_3\ls(1+t)^{-2}2^{(k_+-k)/2}J^k_{31},\\
\end{split}
\end{equation}
where
$$J^k_{31}=\int_{2t/5}^{3t/5}\int_{\frac{9s}{10}\le|y|\le s+M_0}\frac{dyds}{(2M_0+s-|y|)^{1+\mu}(1+2^k|t-s-|x-y||)^2}.$$
For $x=0$, one easily gets
\begin{equation}\label{dt:impr:pf5}
J^k_{31}\ls 2^{5(k_+-k)/4}(1+t).
\end{equation}
For $x\neq0$, let $\phi\in[0,\pi]$ be the angle between $x$ and $y$.
Taking the coordinate transformation $y=(y_1,y_2)\to (\tilde r, \tilde \lambda)$:
\begin{equation*}
\tilde r=|y|,\quad \tilde \lambda=\frac{|y|-|x-y|}{|x|}\in[-1,1],
\end{equation*}
whose Jacobian is
\begin{equation*}
\Big|\frac{\p(\tilde r,\tilde \lambda)}{\p(y_1,y_2)}\Big|=\frac{|\sin\phi|}{|x-y|}.
\end{equation*}
If $\phi\ge\pi/6$, then it holds that
\begin{equation}\label{dt:impr:pf6}
\begin{split}
1-{\tilde \lambda}^2
&=\frac{|x|^2-(|y|-|x-y|)^2}{|x|^2}=\frac{(|x|+|x-y|-|y|)(|x|+|y|-|x-y|)}{|x|^2}\\
&=\frac{(|x|+|x-y|-|y|)(|x|^2+|y|^2+2|x||y|-|x-y|^2)}{|x|^2(|x|+|x-y|+|y|)}\\
&=\frac{2|x||y|(|x|+|x-y|-|y|)(1+\cos\phi)}{|x|^2(|x|+|x-y|+|y|)}\\
&\ls\cos^2\frac{\phi}{2}=\frac{\sin^2\phi}{4\sin^2\frac{\phi}{2}}\ls|x-y|^2\Big|\frac{\p(\tilde r,\tilde \lambda)}{\p(y_1,y_2)}\Big|^2.
\end{split}
\end{equation}
If $\phi\le\pi/6$, then $|y|\ge3|x|$ implies $\psi>\pi/2$, where $\psi$ is the angle between $x$ and $x-y$ ( or between $-x$ and $y-x$).
On the other hand, it follows from  $2|y|/3\le|x-y|\le4|y|/3$ that $\frac{\sin\phi}{\sin\psi}=\frac{|y-x|}{|y|}\in[2/3,4/3]$.
Analogously to \eqref{dt:impr:pf6}, we can achieve
\begin{equation*}
\begin{split}
1-{\tilde \lambda}^2
&=\frac{(|x|+|x-y|-|x-(x-y)|)(|x|+|y|-|x-y|)}{|x|^2}\\
&=\frac{2|x||x-y|(1+\cos\psi)(|x|+|y|-|x-y|)}{|x|^2(|x|+|x-y|+|y|)}\\
&\ls\cos^2\frac{\psi}{2}=\frac{\sin^2\psi}{4\sin^2\frac{\psi}{2}}\ls|x-y|^2\Big|\frac{\p(\tilde r,\tilde \lambda)}{\p(y_1,y_2)}\Big|^2.
\end{split}
\end{equation*}
Therefore,
\begin{equation*}
\big|\frac{\p(\tilde r,\tilde \lambda)}{\p(y_1,y_2)}\big|^{-1}\ls\frac{|x-y|}{\sqrt{1-{\tilde \lambda}^2}}\ls\frac{1+t}{\sqrt{1-{\tilde \lambda}^2}}.
\end{equation*}
In the new coordinate $(\tilde r,\tilde \lambda)$, the integral $J^k_{31}$ can be treated as
\begin{equation*}
J^k_{31}\ls(1+t)\int_0^t\int_0^{s+M_0}\int_{-1}^1\frac{d\tilde \lambda d\tilde rds}{\sqrt{1-{\tilde \lambda}^2}(2M_0+s-\tilde r)^{1+\mu}(1+2^k|t-s-\tilde r+|x||\tilde \lambda|)^2}.
\end{equation*}
Set $\zeta_\pm=s\pm \tilde r$. Then has
\begin{equation}\label{dt:impr:pf7}
\begin{split}
J^k_{31}&\ls2^{5(k_+-k)/4}(1+t)\int_{-1}^1\int_{-M_0}^t\int_0^{2t+M_0}\frac{d\zeta_+d\zeta_-d\tilde \lambda}{\sqrt{1-{\tilde \lambda}^2}(2M_0+\zeta_-)^{1+\mu}
(1+|t-\zeta_++|x||\tilde \lambda|)^{5/4}}\\
&\ls2^{5(k_+-k)/4}(1+t)\int_{-1}^1\int_0^t\frac{d\zeta_-d\tilde \lambda}{\sqrt{1-{\tilde \lambda}^2}(2M_0+\zeta_-)^{1+\mu}}\\
&\ls2^{5(k_+-k)/4}(1+t).
\end{split}
\end{equation}
Collecting \eqref{dt:impr:pf1}-\eqref{dt:impr:pf5} and \eqref{dt:impr:pf7} yields \eqref{dt:impr}.
\end{proof}

Based on Lemma \ref{lem:dt:impr} and Lemma \ref{lem:pw:div},  we have
\begin{corollary}\label{coro:dt:impr}
Let $w$ be the solution to the linear initial boundary value problem
\begin{equation*}
\left\{
\begin{aligned}
&\Box w=\p_tF(t,x),\qquad \qquad (t,x)\in[0,+\infty)\times\cK,\\
&w(t,x)=0,\qquad \qquad \qquad (t,x)\in[0,+\infty)\times\p\cK,\\
&(w,\p_tw)(0,x)=(0,0),\qquad x\in\cK,
\end{aligned}
\right.
\end{equation*}
provided that $\supp_yF(s,y)\subset\{y: \frac{9s}{10}\le|y|\le s+M_0\}$ and $M_0>0$. Then for any constant $R>0$ and $\mu>0$, it holds that
\begin{equation}\label{YHC-9}
\w{t}|w(t,x)|\ls\sup_{y}|F(0,y)|+\sup_{(s,y)\in[0,t]\times\R^2}\w{y}^{1/2}\cW_{1,1+\mu}|(1-\Delta)^2F(s,y)|,
\end{equation}
where $x\in\{x\in\R^2:|x|\le\frac{t}{10}~{\rm or}~|x|\le R\}$.
\end{corollary}
\begin{proof}
The solution can be decomposed into $w=\p_tw_0+w_r$, where
\begin{equation*}
\Box w_0=F(t,x),\quad (w_0,\p_tw_0)(0,x)=(0,0); \quad \Box w_r=0,\quad (w_r,\p_tw_r)(0,x)=(0,-F(0,x)).
\end{equation*}
This, together with \eqref{pw:ivp2} and \eqref{dt:impr}, yields the proof of \eqref{YHC-9}.
\end{proof}

\section{Pointwise estimates of the linear initial boundary value problem}\label{sect4}

This section aims to establish some important pointwise estimates by use of Lemma \ref{lem:dt:impr}.
\begin{lemma}
Suppose that the obstacle $\cO\subset\Bbb R^2$ is star-shaped, $\cK=\R^2\setminus\cO$ and let $w$ be the solution of the IBVP
\begin{equation}\label{dt:loc:eqn}
\left\{
\begin{aligned}
&\Box w=F_0(t,x)+F_1(t,x),\qquad(t,x)\in[0,\infty)\times\cK,\\
&w=0,\qquad\qquad\qquad\qquad \qquad (t,x)\in[0,\infty)\times\p\cK,\\
&(w,\p_tw)(0,x)=(w_0,w_1)(x),\quad x\in\cK,
\end{aligned}
\right.
\end{equation}
where $(w_0,w_1)$ has compact support, $\supp_x F_0(t,x)\subset\{x: 9t/10\le|x|\le t+M_0\}$ and $\supp_x F_1(t,x)\subset\{x: |x|\le t+M_0\}$.
Then for $\mu\in(0,1/2)$, one has
\begin{equation}\label{dt:loc}
\begin{split}
&\w{t}\|\p^{\le1}\p_tw(t,x)\|_{L^2(\cK_R)}
\ls\|(w_0,w_1)\|_{H^2(\cK)}+\sum_{l=0,1}\sum_{|a|\le1}\sup_{s\in[0,t]}\w{s}\|\p^aF_l(s)\|_{L^2(\cK_3)}\\
&\quad+\sum_{|a|\le5}\sup_{(s,y)\in[0,t]\times\overline{\R^2\setminus\cK_2}}\w{y}^{1/2}\cW_{1,1+\mu}(s,y)|\p^aF_0(s,y)|\\
&\quad+\sum_{|a|\le2}\sup_{(s,y)\in[0,t]\times\overline{\R^2\setminus\cK_2}}\w{y}^{1/2}\cW_{1+\mu,1}(s,y)|Z^aF_1(s,y)|.
\end{split}
\end{equation}
\end{lemma}
\begin{remark}
The support property of $F_0(t,x)$ in \eqref{dt:loc:eqn} is essential in the proof of Theorem \ref{thm1} without the second null condition.
The estimate \eqref{dt:loc} actually improves the estimate (5.25) of \cite{HouYinYuan24}
\begin{equation*}
\begin{split}
&\w{t}\|\p_tw(t,x)\|_{L^\infty(\cK_R)}
\ls\|(w_0,w_1)\|_{H^3(\cK)}+\sum_{|a|\le2}\sup_{s\in[0,t]}\w{s}\|\p^aF(s)\|_{L^2(\cK_3)}\\
&\quad+\sum_{|a|\le2}\sup_{y\in\cK}|\p^aF(0,y)|+\sum_{|a|\le3}\sup_{(s,y)\in[0,t]\times\cK}\w{y}^{1/2}\cW_{1+\mu+\nu,1}(s,y)|Z^aF(s,y)|,
\end{split}
\end{equation*}
where $w$ is the solution of $\Box w=F$,  and $\mu,\nu>0$ are constants.

\end{remark}
\begin{proof}
Although the proof is somewhat similar to that of (5.25) in \cite{HouYinYuan24} through utilizing Lemma \ref{lem:dt:impr},
we still give the details for readers' convenience.
Let $w_1^b,w_2^b$ be the solutions to
\begin{equation}\label{dt:loc:pf1}
\left\{
\begin{aligned}
&\Box w_1^b=(1-\chi_{[2,3]}(x))(F_0(t,x)+F_1(t,x)),\\
&\Box w_2^b=\chi_{[2,3]}(x)(F_0(t,x)+F_1(t,x)),\\
&w_1^b|_{\p\cK}=w_2^b|_{\p\cK}=0,\\
&(w_1^b,\p_tw_1^b)(0,x)=(w_0,w_1)(x),\quad (w_2^b,\p_tw_2^b)(0,x)=(0,0),
\end{aligned}
\right.
\end{equation}
respectively.
It follows from the uniqueness of smooth solution to the IBVP with \eqref{dt:loc:eqn} and \eqref{dt:loc:pf1} that $w=w_1^b+w_2^b$.
By \eqref{loc:decay} with $m=2$ and $\eta=1$, one has
\begin{equation}\label{dt:loc:pf2}
\w{t}\|\p^{\le1}\p_tw_1^b(t,x)\|_{L^2(\cK_R)}
\ls\|(w_0,w_1)\|_{H^2(\cK)}+\sum_{l=0,1}\sum_{|a|\le1}\sup_{s\in[0,t]}\w{s}\|\p^aF_l(s)\|_{L^2(\cK_3)}.
\end{equation}
Next, we turn to the estimate of $w_2^b$.
Let $w_2^c$ be the solution of the Cauchy problem
\begin{equation}\label{dt:loc:pf3}
\Box w_2^c=\left\{
\begin{aligned}
&\chi_{[2,3]}(x)(F_0(t,x)+F_1(t,x)),\qquad && x\in\R^2\setminus\cO,\\
&0,\qquad && x\in\overline{\cO}
\end{aligned}
\right.
\end{equation}
with $(w_2^c, \p_tw_2^c)(0,x)=(0,0)$.
It is easy to find that $v_2^b=\p_tw_2^b-\chi_{[1,2]}(x)\p_tw_2^c$ defined on $[0,\infty)\times\cK$ is the solution of the IBVP
\begin{equation*}
\Box v_2^b=[\Delta,\chi_{[1,2]}]\p_tw_2^c,\quad v_2^b|_{\p\cK}=0
\end{equation*}
with $(v_2^b, \p_tv_2^b)(0,x)=(0,0)$.
Since $[\Delta,\chi_{[1,2]}]\p_tw_2^c$ has a compact support in the space variable, the estimate of $v_2^b$ is analogous to that of $w_1^b$ in \eqref{dt:loc:pf2} and it holds that
\begin{equation}\label{dt:loc:pf4}
\begin{split}
\w{t}\|\p^{\le1}v_2^b(t,x)\|_{L^2(\cK_R)}
&\ls\sup_{s\in[0,t]}\w{s}\|[\Delta,\chi_{[1,2]}]\p_tw_2^c(s)\|_{L^2(\cK)}\\
&\ls\sup_{s\in[0,t]}\w{s}\|\p_x^{\le1}\p_tw_2^c(s)\|_{L^\infty(|x|\le2)}.
\end{split}
\end{equation}
On the other hand, the solution of problem \eqref{dt:loc:pf3} can be decomposed into $w_2^c=w_{2,1}^c+w_{2,2}^c$,
where $w_{2,l}^c$ are the solutions of the problems $\Box w_{2,l}^c=\chi_{[2,3]}(x)F_l(t,x)$ with  $(w_{2,l}^c, \p_tw_{2,l}^c)(0,x)=(0,0)$ ($l=0,1$), respectively.
Applying Lemma \ref{lem:dt:impr} and \eqref{dpw:ivp} to $w_{2,0}^c$ and $w_{2,1}^c$, respectively, yields
\begin{equation}\label{dt:loc:pf5}
\begin{split}
\w{t}\|\p_x^{\le1}\p_tw_{2,0}^c(t,x)\|_{L^\infty(|x|\le t/10)}
&\ls\sum_{|a|\le5}\sup_{(s,y)\in[0,t]\times\overline{\R^2\setminus\cK_2}}\w{y}^{1/2}\cW_{1,1+\dl_0}(s,y)|\p^aF_0(s,y)|,\\
\w{x}^{1/2}\w{t-|x}|\p_x^{\le1}\p_tw_{2,1}^c(t,x)|
&\ls\sum_{|a|\le2}\sup_{(s,y)\in[0,t]\times\overline{\R^2\setminus\cK_2}}\w{y}^{1/2}\cW_{1+\dl_0,1}(s,y)|Z^aF_1(s,y)|.
\end{split}
\end{equation}
For $t\ge\max\{10R,20\}$, \eqref{dt:loc} can be achieved by \eqref{dt:loc:pf2}, \eqref{dt:loc:pf4} and \eqref{dt:loc:pf5}.
The proof of \eqref{dt:loc} for $t\le\max\{10R,20\}$ is much easier.
Therefore, the proof of \eqref{dt:loc} is completed.
\end{proof}

\begin{lemma}
Suppose that the obstacle $\cO\subset\Bbb R^2$ is star-shaped, $\cK=\R^2\setminus\cO$ and let $w$ be the solution of the IBVP
\begin{equation*}
\left\{
\begin{aligned}
&\Box w=\p_tF_0(t,x)+F_1(t,x),\qquad(t,x)\in[0,\infty)\times\cK,\\
&w=0,\qquad \qquad \qquad \qquad \qquad \quad (t,x)\in[0,\infty)\times\p\cK,\\
&(w,\p_tw)(0,x)=(w_0,w_1)(x),\qquad x\in\cK,
\end{aligned}
\right.
\end{equation*}
where $(w_0,w_1)$ has compact support and $\supp_x(F_0,F_1)(t,x)\subset\{x: |x|\le t+M_0\}$.
Then one has that for any $\mu,\nu>0$ and $\mu+\nu<1/2$,
\begin{equation}\label{pw:ibvp}
\begin{split}
&\w{t+|x|}^{1/2}\w{t-|x|}^\mu|w(t,x)|
\ls\sum_{|a|\le4}\sup_{s\in[0,t]}\w{s}^{1/2+\mu+\nu}\|\p^a(\p_tF_0,F_1)(s)\|_{L^2(\cK_3)}\\
&\quad+\|(w_0,w_1)\|_{H^5(\cK)}+\ln(2+t+|x|)\sum_{|a|\le6}\sup_{(s,y)\in[0,t]\times\overline{\R^2\setminus\cK_2}}\w{y}^{1/2}
\cW_{1,1}(s,y)|Z^aF_0(s,y)|\\
&\quad+\sum_{|a|\le5}\sup_{y\in\cK}|\p^a(F_0,F_1)(0,y)|+\sum_{|a|\le5}\sup_{(s,y)\in[0,t]\times\overline{\R^2\setminus\cK_2}}
\w{y}^{1/2}\cW_{1+\mu+\nu,1}(s,y)|\p^aF_1(s,y)|.
\end{split}
\end{equation}
Furthermore, if $F_0(t,x)\equiv0$, then it holds that
\begin{equation}\label{dpw:ibvp}
\begin{split}
&\w{x}^{1/2}\w{t-|x|}^{1/2-\mu}|\p w(t,x)|
\ls\|(w_0,w_1)\|_{H^4(\cK)}+\sum_{|a|\le3}\sup_{s\in[0,t]}\w{s}^{1/2}\|\p^aF_1(s)\|_{L^2(\cK_4)}\\
&\quad+\sum_{|a|\le3}\sup_{y\in\cK}|\p^aF_1(0,y)|
+\sum_{|a|\le4}\sup_{(s,y)\in[0,t]\times\overline{\R^2\setminus\cK_2}}\w{y}^{1/2}\cW_{1+\mu+\nu,1}(s,y)|Z^aF_1(s,y)|.
\end{split}
\end{equation}
\end{lemma}
\begin{proof}
Since the estimate \eqref{dpw:ibvp} is just (5.1) of \cite{HouYinYuan24}, we omit the related proof.
Next, we turn to prove \eqref{pw:ibvp}.
Let $w_1^b,w_2^b$ be the solutions to
\begin{equation}\label{pw:ibvp:pf1}
\left\{
\begin{aligned}
&\Box w_1^b=(1-\chi_{[2,3]}(x))(\p_tF_0+F_1),\\
&\Box w_2^b=\chi_{[2,3]}(x)(\p_tF_0+F_1),\\
&w_1^b|_{\p\cK}=w_2^b|_{\p\cK}=0,\\
&(w_1^b,\p_tw_1^b)(0,x)=(w_0,w_1)(x),\quad (w_2^b,\p_tw_2^b)(0,x)=(0,0),
\end{aligned}
\right.
\end{equation}
respectively.
Thus $w=w_1^b+w_2^b$ holds. By \eqref{loc:decay} with $m=|a|+2$, one has
\begin{equation}\label{pw:ibvp:pf2}
\begin{split}
&\w{t}^{\eta}\|\p^aw_1^b(t,x)\|_{L^\infty(\cK_R)}\ls\w{t}^{\eta}\|\p^aw_1^b(t,x)\|_{H^2(\cK_{R+1})}\\
&\ls\|(w_0,w_1)\|_{H^{|a|+2}(\cK)}+\sum_{|b|\le|a|+1}\sup_{s\in[0,t]}\w{s}^{\eta}\|\p^b(\p_tF_0,F_1)(s)\|_{L^2(\cK_3)},
\end{split}
\end{equation}
where $\eta\in(0,1)$.

We now focus on the estimate of $w_1^b$ in the region $|x|\ge R$.
Let $w_1^c=\chi_{[3,4]}(x)w_1^b$ be the solution of the Cauchy problem
\begin{equation*}
\left\{
\begin{aligned}
&\Box w_1^c=-[\Delta,\chi_{[3,4]}(x)]w_1^b,\\
&(w_1^c,\p_tw_1^c)(0,x)=(\chi_{[3,4]}w_0,\chi_{[3,4]}w_1)(x).
\end{aligned}
\right.
\end{equation*}
Note that $[\Delta,\chi_{[3,4]}(x)]w_1^b$ is supported in $\{x\in\R^2,|x|\le4\}$. Then it follows from \eqref{pw:ivp1}
and \eqref{pw:compact} that
\begin{equation}\label{pw:ibvp:pf3}
\begin{split}
&\quad \w{t+|x|}^{1/2}\w{t-|x|}^{\mu}|w_1^c(t,x)|\\
&\ls\|(w_0,w_1)\|_{W^{1,\infty}(\cK)}+\sum_{|a|\le2}\sup_{s\in[0,t]}\w{s}^{1/2+\mu+\nu}\|\p^a[\Delta,\chi_{[3,4]}]w_1^b(s)\|_{L^\infty}\\
&\ls\|(w_0,w_1)\|_{H^3(\cK)}+\sum_{|a|\le3}\sup_{s\in[0,t]}\w{s}^{1/2+\mu+\nu}\|\p^aw_1^b(s)\|_{L^\infty(\cK_4)}\\
&\ls\|(w_0,w_1)\|_{H^5(\cK)}+\sum_{|a|\le4}\sup_{s\in[0,t]}\w{s}^{1/2+\mu+\nu}\|\p^a(\p_tF_0,F_1)(s)\|_{L^2(\cK_3)},
\end{split}
\end{equation}
where \eqref{pw:ibvp:pf2} has been used for $\eta=1/2+\mu+\nu$.

Collecting \eqref{pw:ibvp:pf2} and \eqref{pw:ibvp:pf3} with $w_1^c=\chi_{[3,4]}(x)w_1^b$ and $R>4$ yields
\begin{equation}\label{pw:ibvp:pf4}
\begin{split}
&\quad \w{t+|x|}^{1/2}\w{t-|x|}^{\mu}|w_1^b(t,x)|\\
&\ls\|(w_0,w_1)\|_{H^5(\cK)}+\sum_{|a|\le4}\sup_{s\in[0,t]}\w{s}^{1/2+\mu+\nu}\|\p^a(\p_tF_0,F_1)(s)\|_{L^2(\cK_3)}.
\end{split}
\end{equation}
Next, we deal with $w_2^b$.
Let $w_{2,0}^c$ and $w_{2,1}^c$ be the Cauchy problems of
\begin{equation*}
\Box w_{2,0}^c=\p_t(\chi_{[2,3]}(x)F_0),\quad \Box w_{2,1}^c=\chi_{[2,3]}(x)F_1,
\quad (w_{2,l}^c,\p_tw_{2,l}^c)(0,x)=(0,0),\quad l=0,1,
\end{equation*}
respectively.
Then $v_2^b=w_2^b-\chi_{[1,2]}(x)(w_{2,0}^c+w_{2,1}^c)$ is the solution of the IBVP
\begin{equation*}
\Box v_2^b=[\Delta,\chi_{[1,2]}(x)](w_{2,0}^c+w_{2,1}^c),\quad v_2^b|_{\p\cK}=0,\quad (v_2^b,\p_tv_2^b)(0,x)=(0,0).
\end{equation*}
Since the estimate of $v_2^b$ is similar to that of $w_1^b$ in \eqref{pw:ibvp:pf4}, we then have
\begin{equation}\label{pw:ibvp:pf5}
\begin{split}
&\quad \w{t+|x|}^{1/2}\w{t-|x|}^{\mu}|v_2^b(t,x)|\\
&\ls\sum_{|a|\le4}\sup_{s\in[0,t]}\w{s}^{1/2+\mu+\nu/2}\|\p^a[\Delta,\chi_{[1,2]}](w_{2,0}^c+w_{2,1}^c)(s)\|_{L^2}\\
&\ls\sum_{|a|\le5}\sup_{s\in[0,t]}\w{s}^{1/2+\mu+\nu/2}\|\p^a(w_{2,0}^c+w_{2,1}^c)(s,y)\|_{L^\infty(|y|\le2)}.
\end{split}
\end{equation}
Applying \eqref{pw:ivp:div'} with $\eta=1/2-\mu-\nu$ and \eqref{pw:ivp1} to $\p^aw_{2,0}^c$ and $\p^aw_{2,1}^c$, respectively,
yields that for any $|a|\le5$,
\begin{equation}\label{pw:ibvp:pf6}
\begin{split}
&\quad \w{t+|x|}^{1/2}\w{t-|x|}^{\mu+\nu}|\p^aw_{2,0}^c(t,x)|\\
&\ls\sum_{|b|\le5}\sup_{y\in\cK}|\p^bF_0(0,y)|
+\ln(2+t+|x|)\sum_{|b|\le6}\sup_{(s,y)\in[0,t]\times\overline{\R^2\setminus\cK_2}}\w{y}^{1/2}\cW_{1,1}(s,y)|Z^b F_0(s,y)|,\\
&\quad \w{t+|x|}^{1/2}\w{t-|x|}^{\mu+\nu/2}|\p^aw_{2,1}^c(t,x)|\\
&\ls\sum_{|b|\le4}\sup_{y\in\cK}|\p^bF_1(0,y)|
+\sum_{|b|\le5}\sup_{(s,y)\in[0,t]\times\overline{\R^2\setminus\cK_2}}\w{y}^{1/2}\cW_{1+\mu+\nu,1}(s,y)|\p^b F_1(s,y)|,
\end{split}
\end{equation}
where $\ds\sup_{y\in\cK}|\p^bF_0(0,y)|$ and $\ds\sup_{y\in\cK}|\p^bF_1(0,y)|$ come from the initial data.
Hence, \eqref{pw:ibvp} can be obtained by \eqref{pw:ibvp:pf4}-\eqref{pw:ibvp:pf6} with $w=w_1^b+w_2^b$ and $w_2^b=v_2^b+\chi_{[1,2]}(x)(w_{2,0}^c+w_{2,1}^c)$.
\end{proof}

\begin{lemma}
Suppose that the obstacle $\cO\subset\Bbb R^2$ is star-shaped, $\cK=\R^2\setminus\cO$ and let $w$ be the solution of the IBVP
\begin{equation*}
\left\{
\begin{aligned}
&\Box w=\p_tF^0(t,x)+\sum_{\Xi=0,1,2,\theta}\p_{\Xi}G^{\Xi}(t,x)+G^r_1(t,x)+G^r_2(t,x),\quad(t,x)\in[0,\infty)\times\cK,\\
&w=0,\qquad(t,x)\in[0,\infty)\times\p\cK,\\
&(w,\p_tw)(0,x)=(w_0,w_1)(x),\quad x\in\cK,
\end{aligned}
\right.
\end{equation*}
where $(w_0,w_1)$ has compact support, $\supp_x F^0(t,x)\subset\{x: \frac{9t}{10}\le|x|\le t+M_0\}$ and $G^r_1,G^r_2,G^{\Xi}$ ($\Xi=0,1,2,\theta$) are supported in $\{x: |x|\le t+M_0\}$.
Then one has that for any $\eta,\mu,\nu>0$, $\eta<1/2$ and $\mu+\nu<1/2$,
\begin{equation*}\label{dpw:div}
\begin{split}
&\quad\w{t}\|\p w(t,x)\|_{L^\infty(|x|\le t/10)}+m^{\#}\w{x}^{1/2}\w{t-|x|}^{1-\eta}|\p w(t,x)|\\
&\ls\|(w_0,w_1)\|_{H^9(\cK)}+\sum_{|a|\le8}\sum_{\Xi=0,1,2,\theta}\sup_{y\in\cK}|\p^a(\p^{\le1}F^0,\p^{\le1}G^{\Xi},G^r_1,G^r_2)(0,y)|\\
&\quad+\sum_{|a|\le8}\sum_{\Xi=0,1,2,\theta}\sup_{s\in[0,t]}\w{s}\|\p^a(\p_tF^0,\p^{\le1}G^{\Xi},G^r_1,G^r_2)(s)\|_{L^2(\cK_4)}\\
&\quad+\sum_{i=1,2}\sup_{s\in[0,t]}\w{s}^{3/2+\mu}\|G^i(s)\|_{L^\infty(\cK_3)}
+\sup_{(s,y)\in[0,t]\times\overline{\R^2\setminus\cK_2}}\w{y}^{1/2}\cW_{3/2+\mu,1}(s,y)|G^r_1(s,y)|\\
&\quad+\ln(2+t+|x|)\sup_{(s,y)\in[0,t]\times\overline{\R^2\setminus\cK_2}}\w{y}^{1/2}\cW_{3/2,1+\nu}(s,y)|G^r_2(s,y)|\\
\end{split}
\end{equation*}
\begin{equation}\label{dpw:div}
\begin{split}
&\quad+\sum_{|a|\le9}\sum_{\Xi=0,1,2,\theta}\sup_{(s,y)\in[0,t]\times\overline{\R^2\setminus\cK_2}}
\w{y}^{1/2}\cW_{1+\mu+\nu,1}(s,y)|Z^a(Z^{\le1}G^{\Xi},G^r_1,G^r_2)(s,y)|\\
&\quad+(1+m^{\#}\ln(2+t+|x|))\sum_{|a|\le13}\sup_{(s,y)\in[0,t]\times\overline{\R^2\setminus\cK_2}}\w{y}^{1/2}\cW_{1,1+\mu}(s,y)|Z^aF^0(s,y)|,
\end{split}
\end{equation}
where $m^{\#}=0,1$.
\end{lemma}
\begin{remark}
Note that when $m^{\#}=0$, there is no $\ln(2+t+|x|)$ factor in the last line of \eqref{dpw:div} for the estimate of $\|\p w(t,x)\|_{L^\infty(|x|\le t/10)}$.
\end{remark}
\begin{proof}
First, recall (5.4) of \cite{HouYinYuan24}
\begin{equation}\label{dpw:div:pf1}
\begin{split}
&\w{x}^{1/2}\w{t-|x|}|\p w^{div}(t,x)|
\ls\|(w_0,w_1)\|_{H^9(\cK)}+\sum_{|a|\le8}\sum_{\Xi=0,1,2,\theta}\sup_{y\in\cK}|\p^a(\p^{\le1}G^{\Xi},G^r_1)(0,y)|\\
&\quad+\sum_{|a|\le8}\sum_{\Xi=0,1,2,\theta}\sup_{s\in[0,t]}\w{s}\|\p^a(\p^{\le1}G^{\Xi},G^r_1)(s)\|_{L^2(\cK_4)}
+\sum_{i=1,2}\sup_{s\in[0,t]}\w{s}^{3/2+\mu}\|G^i(s)\|_{L^\infty(\cK_3)}\\
&\quad+\sup_{(s,y)\in[0,t]\times\overline{\R^2\setminus\cK_2}}\w{y}^{1/2}\cW_{3/2+\mu,1}(s,y)|G^r_1(s,y)|\\
&\quad+\sum_{|a|\le9}\sum_{\Xi=0,1,2,\theta}\sup_{(s,y)\in[0,t]\times\overline{\R^2\setminus\cK_2}}
\w{y}^{1/2}\cW_{1+\mu+\nu,1}(s,y)|Z^a(G^r_1,Z^{\le1}G^{\Xi})(s,y)|,
\end{split}
\end{equation}
provided that $w^{div}$ is the solution of
\begin{equation*}
\left\{
\begin{aligned}
&\Box w^{div}=\sum_{\Xi=0,1,2,\theta}\p_{\Xi}G^{\Xi}(t,x)+G^r_1(t,x),\qquad(t,x)\in[0,\infty)\times\cK,\\
&w^{div}=0,\qquad \qquad \qquad \qquad \qquad \qquad \qquad \quad (t,x)\in[0,\infty)\times\p\cK,\\
&(w^{div},\p_tw^{div})(0,x)=(w_0,w_1)(x),\quad\qquad \qquad x\in\cK.
\end{aligned}
\right.
\end{equation*}
Thus, it suffices to treat the estimates of $\tilde w=w-w^{div}$, where  $\tilde w$ is the solution of
\begin{equation*}
\Box\tilde w=\p_tF^0(t,x)+G^r_2(t,x),\quad \tilde w|_{\p\cK}=0,\quad (\tilde w,\p_t\tilde w)(0,x)=(0,0).
\end{equation*}
Let $w_1^b$ and $w_2^b$ be the solutions of
\begin{equation}\label{dpw:div:pf2}
\left\{
\begin{aligned}
&\Box w_1^b=(1-\chi_{[2,3]}(x))(\p_tF^0+G^r_2),\\
&\Box w_2^b=\chi_{[2,3]}(x)(\p_tF^0+G^r_2),\\
&w_1^b|_{\p\cK}=w_2^b|_{\p\cK}=0,\\
&(w_1^b,\p_tw_1^b)(0,x)=(w_2^b,\p_tw_2^b)(0,x)=(0,0),
\end{aligned}
\right.
\end{equation}
respectively.
By \eqref{loc:decay} with $m=|a|+2$, one has
\begin{equation}\label{dpw:div:pf3}
\begin{split}
\w{t}\|\p^aw_1^b(t,x)\|_{L^\infty(\cK_R)}&\ls\w{t}\|\p^aw_1^b(t,x)\|_{H^2(\cK_{R+1})}\\
&\ls\sum_{|b|\le|a|+1}\sup_{s\in[0,t]}\w{s}\|\p^b(\p_tF^0,G^r_2)(s)\|_{L^2(\cK_3)}.
\end{split}
\end{equation}
Next, we turn to the estimate of $w_1^b$ in the region $|x|\ge R$.
It is easy to find that $w_1^c=\chi_{[3,4]}(x)w_1^b$ is the solution of the Cauchy problem
\begin{equation*}
\Box w_1^c=-[\Delta,\chi_{[3,4]}(x)]w_1^b,\quad (w_1^c,\p_tw_1^c)(0,x)=(0,0).
\end{equation*}
By applying \eqref{dpw:compact} to $\Box w_1^c$ together with \eqref{dpw:div:pf3}, one obtains
\begin{equation}\label{dpw:div:pf4}
\begin{split}
\w{x}^{1/2}\w{t-|x|}|\p w_1^c(t,x)|
&\ls\sum_{|b|\le6}\sup_{s\in[0,t]}\w{s}\|\p^b[\Delta,\chi_{[3,4]}]w_1^b(s)\|_{L^\infty}\\
&\ls\sum_{|b|\le7}\sup_{s\in[0,t]}\w{s}\|\p^bw_1^b(s)\|_{L^\infty(\cK_4)}\\
&\ls\sum_{|b|\le8}\sup_{s\in[0,t]}\w{s}\|\p^b(\p_tF^0,G^r_2)(s)\|_{L^2(\cK_3)}.
\end{split}
\end{equation}
Note that $\p w_1^b=\p w_1^c$ holds in the region $|x|\ge4$. Then collecting \eqref{dpw:div:pf3}-\eqref{dpw:div:pf4} with $R>4$ leads to
\begin{equation}\label{dpw:div:pf5}
\w{x}^{1/2}\w{t-|x|}|\p w_1^b(t,x)|\ls\sum_{|b|\le8}\sup_{s\in[0,t]}\w{s}\|\p^b(\p_tF^0,G^r_2)(s)\|_{L^2(\cK_3)}.
\end{equation}
We next deal with $\p_tw_2^b$ and $\p_iw_2^b$ ($i=1,2$), separately.
Let $w_2^c$ and $w_r^c$ be the solutions of the Cauchy problems
\begin{equation}\label{dpw:div:pf6}
\Box w_2^c=\chi_{[2,3]}(x)F^0(t,x),\quad \Box w_r^c=\chi_{[2,3]}(x)G^r_2(t,x)
\end{equation}
with $(w_2^c, \p_tw_2^c)(0,x)=(0,0)$ and $(w_r^c, \p_tw_r^c)(0,x)=(0,0)$, respectively.
Then, by \eqref{dpw:div:pf2} and \eqref{dpw:div:pf6}, $v_2^b=w_2^b-\chi_{[1,2]}(x)(\p_tw_2^c+w_r^c)$ is the solution of the IBVP
\begin{equation*}
\Box v_2^b=[\Delta,\chi_{[1,2]}(x)](\p_tw_2^c+w_r^c),\quad v_2^b|_{\p\cK}=0,\quad (v_2^b,\p_tv_2^b)(0,x)=(0,-\chi_{[2,3]}(x)F^0(0,x)).
\end{equation*}
Since the estimate of $v_2^b$ is analogous to that of $w_1^b$ in \eqref{dpw:div:pf5}, then it holds that
\begin{equation}\label{dpw:div:pf7}
\begin{split}
\w{x}^{1/2}\w{t-|x|}|\p v_2^b(t,x)|\ls\|F^0(0,y)\|_{H^8(\cK)}+\sum_{|b|\le9}\sup_{s\in[0,t]}\w{s}\|\p^b(\p_tw_2^c,w_r^c)(s,y)\|_{L^\infty(|y|\le2)}.
\end{split}
\end{equation}
For the estimates of $w_2^c$ and $w_r^c$, applying Lemma \ref{lem:dt:impr}, \eqref{dtpw:ivp}, \eqref{pw:ivp3}, \eqref{dpw:ivp} to \eqref{dpw:div:pf6}, respectively, yields
\begin{equation}\label{dpw:div:pf8}
\begin{split}
\sum_{|b|\le9}\w{t}|\p^b\p_tw_2^c(t,x)|
&\ls\sum_{|a|\le13}\sup_{(s,y)\in[0,t]\times\overline{\R^2\setminus\cK_2}}\w{y}^{1/2}\cW_{1,1+\mu}(s,y)|Z^aF^0(s,y)|\\
&\quad+\sum_{|a|\le8}\sup_{y\in\cK}|\p^aF^0(0,y)|,\qquad |x|\le t/10~{\rm or}~|x|\le2,\\
\w{x}^{1/2}\w{t-|x|}^{1-\eta}|\p^2_tw_2^c(t,x)|
&\ls\ln(2+t+|x|)\sum_{|a|\le2}\sup_{(s,y)\in[0,t]\times\overline{\R^2\setminus\cK_2}}\w{y}^{1/2}\cW_{1,1}(s,y)|Z^aF^0(s,y)|,\\
\w{t+|x|}^{1/2}\w{t-|x|}^{1/2}|w_r^c(t,x)|
&\ls\ln(2+t+|x|)\sup_{(s,y)\in[0,t]\times\overline{\R^2\setminus\cK_2}}\w{y}^{1/2}\cW_{3/2,1+\nu}(s,y)|G^r_2(s,y)|,\\
\sum_{|b|\le8}\w{x}^{1/2}\w{t-|x|}|\p^b\p w_r^c(t,x)|
&\ls\sum_{|a|\le9}\sup_{(s,y)\in[0,t]\times\overline{\R^2\setminus\cK_2}}\w{y}^{1/2}\cW_{1+\mu,1}(s,y)|Z^aG^r_2(s,y)|\\
&\quad+\sum_{|a|\le7}\sup_{y\in\cK}|\p^aG^r_2(0,y)|.
\end{split}
\end{equation}
Collecting \eqref{dpw:div:pf7}, \eqref{dpw:div:pf8} with $\p_tw_2^b=\p_tv_2^b+\chi_{[1,2]}(x)(\p^2_tw_2^c+\p_tw_r^c)$ derives
\begin{equation}\label{dpw:div:pf9}
\begin{split}
&\w{x}^{1/2}\w{t-|x|}^{1-\eta}|\p_tw_2^b(t,x)|\ls\sum_{|a|\le8}\sup_{y\in\cK}|\p^a(F^0,G^r_2)(0,y)|\\
&\qquad+\ln(2+t+|x|)\sum_{|a|\le13}\sup_{(s,y)\in[0,t]\times\overline{\R^2\setminus\cK_2}}\w{y}^{1/2}\cW_{1,1+\mu}(s,y)|Z^aF^0(s,y)|\\
&\qquad+\ln(2+t+|x|)\sup_{(s,y)\in[0,t]\times\overline{\R^2\setminus\cK_2}}\w{y}^{1/2}\cW_{3/2,1+\nu}(s,y)|G^r_2(s,y)|\\
&\qquad+\sum_{|a|\le9}\sup_{(s,y)\in[0,t]\times\overline{\R^2\setminus\cK_2}}\w{y}^{1/2}\cW_{1+\mu,1}(s,y)|Z^aG^r_2(s,y)|
\end{split}
\end{equation}
and for $|x|\le t/10$
\begin{equation}\label{dpw:div:pf10}
\begin{split}
&\w{t}|\p_tw_2^b(t,x)|\ls\sum_{|a|\le8}\sup_{y\in\cK}|\p^a(F^0,G^r_2)(0,y)|\\
&\qquad+\sum_{|a|\le13}\sup_{(s,y)\in[0,t]\times\overline{\R^2\setminus\cK_2}}\w{y}^{1/2}\cW_{1,1+\mu}(s,y)|Z^aF^0(s,y)|\\
&\qquad+\ln(2+t+|x|)\sup_{(s,y)\in[0,t]\times\overline{\R^2\setminus\cK_2}}\w{y}^{1/2}\cW_{3/2,1+\nu}(s,y)|G^r_2(s,y)|\\
&\qquad+\sum_{|a|\le9}\sup_{(s,y)\in[0,t]\times\overline{\R^2\setminus\cK_2}}\w{y}^{1/2}\cW_{1+\mu,1}(s,y)|Z^aG^r_2(s,y)|.
\end{split}
\end{equation}
Next, we turn to the estimate of $\p_xw_2^b$.
Let $\tilde w_2^c$ be the Cauchy problem of
\begin{equation*}
\Box\tilde w_2^c=\p_t(\chi_{[2,3]}(x)F^0(t,x)),
\quad (\tilde w_2^c,\p_t\tilde w_2^c)(0,x)=(0,0).
\end{equation*}
Therefore, $\tilde v_2^b=w_2^b-\chi_{[1,2]}(x)(\tilde w_2^c+w_r^c)$ is the solution of the IBVP
\begin{equation*}
\Box\tilde v_2^b=[\Delta,\chi_{[1,2]}(x)](\tilde w_2^c+w_r^c),\quad \tilde v_2^b|_{\p\cK}=0,\quad (\tilde v_2^b,\p_t\tilde v_2^b)(0,x)=(0,0).
\end{equation*}
Note that the estimate of $\tilde v_2^b$ is similar to that of $v_2^b$ in \eqref{dpw:div:pf7}. Then it holds that
\begin{equation}\label{dpw:div:pf11}
\w{x}^{1/2}\w{t-|x|}|\p\tilde v_2^b(t,x)|
\ls\sum_{|b|\le9}\sup_{s\in[0,t]}\w{s}\|\p^b(\tilde w_2^c+w_r^c)(s,y)\|_{L^\infty(|y|\le2)}.
\end{equation}
By utilizing \eqref{pw:ivp:div} to $\Box\p_i\tilde w_2^c=\p_t\p_i(\chi_{[2,3]}(x)F^0(t,x))$ with $(\tilde w_2^c,\p_t\tilde w_2^c)(0,x)=(0,0)$ and Corollary \ref{coro:dt:impr} to $\Box\p^b\tilde w_2^c=\p_t\p^b(\chi_{[2,3]}(x)F^0(t,x))$, respectively, one has
\begin{equation}\label{dpw:div:pf12}
\begin{split}
\sum_{|b|\le9}\w{t}|\p^b\tilde w_2^c(t,x)|
&\ls\sum_{|a|\le13}\sup_{(s,y)\in[0,t]\times\overline{\R^2\setminus\cK_2}}\w{y}^{1/2}\cW_{1,1+\mu}(s,y)|Z^aF^0(s,y)|\\
&\quad+\sum_{|a|\le9}\sup_{y\in\cK}|\p^aF^0(0,y)|,\qquad |x|\le t/10~{\rm or}~|x|\le2,\\
\w{x}^{1/2}\w{t-|x|}^{1-\eta}|\p_i\tilde w_2^c(t,x)|
&\ls\ln(2+t+|x|)\sum_{|a|\le2}\sup_{(s,y)\in[0,t]\times\overline{\R^2\setminus\cK_2}}\w{y}^{1/2}\cW_{1,1}(s,y)|Z^aF^0(s,y)|.
\end{split}
\end{equation}
It follows from \eqref{dpw:div:pf8}, \eqref{dpw:div:pf11}, \eqref{dpw:div:pf12} and $\p_iw_2^b=\p_i\tilde v_2^b+\chi_{[1,2]}(x)\p_i(\tilde w_2^c+w_r^c)+(\p_i\chi_{[1,2]}(x))(\tilde w_2^c+w_r^c)$ that
\begin{equation}\label{dpw:div:pf13}
\begin{split}
&\w{x}^{1/2}\w{t-|x|}^{1-\eta}|\p_iw_2^b(t,x)|\ls\sum_{|a|\le8}\sup_{y\in\cK}|\p^a(\p^{\le1}F^0,G^r_2)(0,y)|\\
&\qquad+\ln(2+t+|x|)\sum_{|a|\le13}\sup_{(s,y)\in[0,t]\times\overline{\R^2\setminus\cK_2}}\w{y}^{1/2}\cW_{1,1+\mu}(s,y)|Z^aF^0(s,y)|\\
&\qquad+\ln(2+t+|x|)\sup_{(s,y)\in[0,t]\times\overline{\R^2\setminus\cK_2}}\w{y}^{1/2}\cW_{3/2,1+\nu}(s,y)|G^r_2(s,y)|\\
&\qquad+\sum_{|a|\le9}\sup_{(s,y)\in[0,t]\times\overline{\R^2\setminus\cK_2}}\w{y}^{1/2}\cW_{1+\mu,1}(s,y)|Z^aG^r_2(s,y)|
\end{split}
\end{equation}
and for $|x|\le t/10$
\begin{equation}\label{dpw:div:pf14}
\begin{split}
&\w{t}|\p_iw_2^b(t,x)|\ls\sum_{|a|\le8}\sup_{y\in\cK}|\p^a(\p^{\le1}F^0,G^r_2)(0,y)|\\
&\qquad+\sum_{|a|\le13}\sup_{(s,y)\in[0,t]\times\overline{\R^2\setminus\cK_2}}\w{y}^{1/2}\cW_{1,1+\mu}(s,y)|Z^aF^0(s,y)|\\
&\qquad+\ln(2+t+|x|)\sup_{(s,y)\in[0,t]\times\overline{\R^2\setminus\cK_2}}\w{y}^{1/2}\cW_{3/2,1+\nu}(s,y)|G^r_2(s,y)|\\
&\qquad+\sum_{|a|\le9}\sup_{(s,y)\in[0,t]\times\overline{\R^2\setminus\cK_2}}\w{y}^{1/2}\cW_{1+\mu,1}(s,y)|Z^aG^r_2(s,y)|.
\end{split}
\end{equation}
Therefore, \eqref{dpw:div} is achieved by \eqref{dpw:div:pf1}, \eqref{dpw:div:pf5}, \eqref{dpw:div:pf9}-\eqref{dpw:div:pf10} and \eqref{dpw:div:pf13}-\eqref{dpw:div:pf14}.
\end{proof}

\section{Energy estimates and local energy decay estimates}\label{sect5}

\subsection{Bootstrap assumptions on the error solution}

Let $u$ and $w_{hom}$ be the solutions to \eqref{QWE} and \eqref{hom:eqn}, respectively, define the error solution
\begin{equation}\label{err:def}
v(t,x)=u(t,x)-w_{hom}(t,x).
\end{equation}
We make the following bootstrap assumptions for all $t\in [0,T_\ve]$ with $T_\ve=e^{\kappa/\ve^2}$,
\begin{align}
&\sum_{|a|\le N+1}|\p Z^av|\le\ve_1\w{x}^{-1/2}\w{t-|x|}^{\dl_0-1},\label{BA1}\\
&\sum_{|a|\le N+1}|\bar\p Z^av|\le\ve_1\w{x}^{-1/2}\w{t+|x|}^{\dl_0-1},\label{BA2}\\
&\sum_{|a|\le N+1}|Z^av|\le\ve_1\w{t+|x|}^{-1/2}\w{t-|x|}^{3\dl_0-1/2},\label{BA3}\\
&\sum_{|a|\le N+1}|\p Z^av|\le\ve_2^2\w{t}^{-1},\quad |x|\le{t}/10,\label{BA4}
\end{align}
where the small positive constants $\ve_1$ and $\ve_2$ will be determined later, $\ve\le\ve_1\le\ve_2\le1$ and $\dl_0=10^{-3}$.
In addition, for convenience, the following assumption is imposed in the proof of Theorem \ref{thm1}
\begin{equation}\label{BA:time}
\ve_1^2\ln(2+t)\le\ve_2^2\ln(2+t)\ls\frac{\kappa\ve_2^2}{\ve^2}\ls1.
\end{equation}
Note that due to $\supp(u_0,u_1)\subset\{x: |x|\le M_0\}$,
then the solution $u$ and $v$ are supported in $\{x\in\cK:|x|\le t+M_0\}$.

\subsection{Energy estimates}

\begin{lemma}\label{lem:energy:time}
Under the assumptions of Theorem \ref{thm1}, let $v$ be defined by \eqref{err:def} and suppose that \eqref{BA1}-\eqref{BA4} hold.
Then we have that for $t\in [0,T_\ve]$,
\begin{equation}\label{energy:time}
\sum_{i\le2N}\|\p\p_t^iv(t)\|_{L^2(\cK)}\ls\ve.
\end{equation}
\end{lemma}
\begin{proof}
By \eqref{eqn:high} together with $\Box w_{hom}\equiv0$ and \eqref{err:def}, one has
\begin{equation}\label{energy:t1}
\begin{split}
&\quad\Box Z^av=Z^aQ(\p u,\p^2u)\\
&=Q^{\alpha\beta\gamma}(\p^2_{\alpha\beta}Z^av\p_{\gamma}u+\p^2_{\alpha\beta}Z^aw_{hom}\p_{\gamma}u)
+Q^{\alpha\beta\gamma\delta}(\p^2_{\alpha\beta}Z^av\p_{\gamma}u\p_{\delta}u+\p^2_{\alpha\beta}Z^aw_{hom}\p_{\gamma}u\p_{\delta}u)\\
&\quad+\sum_{b+c\le a,b<a}Q^{\alpha\beta\gamma}_{abc}\p^2_{\alpha\beta}Z^bu\p_{\gamma}Z^cu
+\sum_{b+c+d\le a,b<a}Q^{\alpha\beta\gamma\delta}_{abcd}\p^2_{\alpha\beta}Z^bu\p_{\gamma}Z^cu\p_{\delta}Z^du,
\end{split}
\end{equation}
where the summations $\ds\sum_{\alpha,\beta,\gamma=0}^2$ and $\ds\sum_{\alpha,\beta,\gamma,\delta=0}^2$ in \eqref{energy:t1} and the following  are omitted.
Multiplying \eqref{energy:t1} by $e^q\p_tZ^av$ with the ghost weight $q=q(|x|-t)$ and $\ds q(s)=\int_{-\infty}^s\frac{d\lambda}{\w{\lambda}^{3/2}}$ introduced in \cite{Alinhac01a}, yields
\begin{equation}\label{energy:t2}
\begin{split}
&\quad\frac12\p_t[e^q(|\p_tZ^av|^2+|\nabla Z^av|^2)]-\dive(e^q\p_tZ^av\nabla Z^av)+\frac{e^q}{2\w{t-|x|}^{3/2}}|\bar\p Z^av|^2\\
&=Q^{\alpha\beta\gamma}e^q\p_tZ^av(\p^2_{\alpha\beta}Z^av\p_{\gamma}u+\p^2_{\alpha\beta}Z^aw_{hom}\p_{\gamma}u)\\
&\quad+Q^{\alpha\beta\gamma\delta}e^q\p_tZ^av(\p^2_{\alpha\beta}Z^av\p_{\gamma}u\p_{\delta}u+\p^2_{\alpha\beta}Z^aw_{hom}\p_{\gamma}u\p_{\delta}u)\\
&\quad+\sum_{\substack{b+c\le a,\\b<a}}Q^{\alpha\beta\gamma}_{abc}e^q\p_tZ^av\p^2_{\alpha\beta}Z^bu\p_{\gamma}Z^cu
+\sum_{\substack{b+c+d\le a,\\b<a}}Q^{\alpha\beta\gamma\delta}_{abcd}e^q\p_tZ^av\p^2_{\alpha\beta}Z^bu\p_{\gamma}Z^cu\p_{\delta}Z^du,
\end{split}
\end{equation}
where $|\bar\p f|=(|\bar\p_1f|^2+|\bar\p_2f|^2)^{1/2}$.
For the terms in the second and third lines of \eqref{energy:t2}, it follows from the symmetric conditions \eqref{sym:condition} and
direct computations that
\begin{equation}\label{energy:t3}
\begin{split}
&\quad\;Q^{\alpha\beta\gamma}e^q\p_tZ^av\p^2_{\alpha\beta}Z^av\p_{\gamma}u
+Q^{\alpha\beta\gamma\delta}e^q\p_tZ^av\p^2_{\alpha\beta}Z^av\p_{\gamma}u\p_{\delta}u\\
&=\p_{\alpha}(Q^{\alpha\beta\gamma}e^q\p_tZ^av\p_{\beta}Z^av\p_{\gamma}u)
-Q^{\alpha\beta\gamma}\p_tZ^av\p_{\beta}Z^av\p_{\alpha}(e^q\p_{\gamma}u)\\
&\quad-\frac12\p_t(Q^{\alpha\beta\gamma}e^q\p_{\alpha}Z^av\p_{\beta}Z^av\p_{\gamma}u)
+\frac12Q^{\alpha\beta\gamma}\p_{\alpha}Z^av\p_{\beta}Z^av\p_t(e^q\p_{\gamma}u)\\
&\quad+\p_{\alpha}(Q^{\alpha\beta\gamma\delta}e^q\p_tZ^av\p_{\beta}Z^av\p_{\gamma}u\p_{\delta}u)
-Q^{\alpha\beta\gamma\delta}\p_tZ^av\p_{\beta}Z^av\p_{\alpha}(e^q\p_{\gamma}u\p_{\delta}u)\\
&\quad-\frac12\p_t(Q^{\alpha\beta\gamma\delta}e^q\p_{\alpha}Z^av\p_{\beta}Z^av\p_{\gamma}u\p_{\delta}u)
+\frac12Q^{\alpha\beta\gamma\delta}\p_{\alpha}Z^av\p_{\beta}Z^av\p_t(e^q\p_{\gamma}u\p_{\delta}u).
\end{split}
\end{equation}
Integrating \eqref{energy:t2} and \eqref{energy:t3} over $[0,t]\times\cK$ with $Z^a=\p_t^i$ ($i\le2N$) and the boundary
conditions $\p_t^lu|_{[0,\infty)\times\p\cK}=0$ ($l\ge0$), implies
\begin{equation}\label{energy:t4}
\begin{split}
&\quad\|\p\p_t^iv(t)\|^2_{L^2(\cK)}+\int_0^t\int_{\cK}\frac{|\bar\p\p_t^iv|^2}{\w{s-|x|}^{3/2}}dxds\\
&\ls\|\p\p_t^iv(0)\|^2_{L^2(\cK)}(1+\|\p u(0)\|_{L^\infty(\cK)})+\|\p u(t)\|_{L^\infty(\cK)}\|\p\p_t^iv(t)\|^2_{L^2(\cK)}\\
&\quad+\sum_{b+c\le a,b<a}\Big|\int_0^t\int_{\cK}I_1^{abc}dxds\Big|
+\sum_{b+c+d\le a,b<a}\Big|\int_0^t\int_{\cK}I_2^{abcd}dxds\Big|+\int_0^t\int_{\cK}|I^a|dxds,
\end{split}
\end{equation}
where
\begin{equation}\label{energy:t5}
\begin{split}
&I_1^{abc}=Q^{\alpha\beta\gamma}_{abc}e^q\p_tZ^av\p^2_{\alpha\beta}Z^bu\p_{\gamma}Z^cu,\\
&I_2^{abcd}=Q^{\alpha\beta\gamma\delta}_{abcd}e^q\p_tZ^av\p^2_{\alpha\beta}Z^bu\p_{\gamma}Z^cu\p_{\delta}Z^du,\\
&I^a=-Q^{\alpha\beta\gamma}e^q\p_tZ^av\p_{\beta}Z^av\p_{\alpha}q\p_{\gamma}u
-Q^{\alpha\beta\gamma}e^q\p_tZ^av\p_{\beta}Z^av\p^2_{\alpha\gamma}u\\
&\quad+\frac12Q^{\alpha\beta\gamma}e^q\p_{\alpha}Z^av\p_{\beta}Z^av\p_tq\p_{\gamma}u
+\frac12Q^{\alpha\beta\gamma}e^q\p_{\alpha}Z^av\p_{\beta}Z^av\p^2_{0\gamma}u\\
&\quad-Q^{\alpha\beta\gamma\delta}\p_tZ^av\p_{\beta}Z^av\p_{\alpha}(e^q\p_{\gamma}u\p_{\delta}u)
+\frac12Q^{\alpha\beta\gamma\delta}\p_{\alpha}Z^av\p_{\beta}Z^av\p_t(e^q\p_{\gamma}u\p_{\delta}u)\\
&\quad+Q^{\alpha\beta\gamma}e^q\p_tZ^av\p^2_{\alpha\beta}Z^aw_{hom}\p_{\gamma}u
+Q^{\alpha\beta\gamma\delta}e^q\p_tZ^av\p^2_{\alpha\beta}Z^aw_{hom}\p_{\gamma}u\p_{\delta}u.
\end{split}
\end{equation}
Although $Z^a=\p_t^i$ has been selected, the estimates of $I_1^{abc},I_2^{abcd},I^a$ also hold for the vector fields $Z\in\{\p,\Omega\}$.
By the virtue of definition \eqref{err:def}, one has
\begin{equation}\label{energy:t6}
\begin{split}
&I_1^{abc}=Q^{\alpha\beta\gamma}_{abc}e^q\p_tZ^av\p^2_{\alpha\beta}Z^bw_{hom}\p_{\gamma}Z^cw_{hom}+I_{11}^{abc},\\
&I_{11}^{abc}=Q^{\alpha\beta\gamma}_{abc}e^q\p_tZ^av(\p^2_{\alpha\beta}Z^bv\p_{\gamma}Z^cv+\p^2_{\alpha\beta}Z^bv\p_{\gamma}Z^cw_{hom}
+\p^2_{\alpha\beta}Z^bw_{hom}\p_{\gamma}Z^cv).
\end{split}
\end{equation}
For the first term on the right hand side of the first line in \eqref{energy:t6},
it follows from  \eqref{null:structure}, \eqref{null:high} and \eqref{hom:pw1} that
\begin{equation*}
\begin{split}
|Q^{\alpha\beta\gamma}_{abc}\p^2_{\alpha\beta}Z^bw_{hom}\p_{\gamma}Z^cw_{hom}|
&\ls|\bar\p\p Z^bw_{hom}||\p Z^cw_{hom}|+|\p^2Z^bw_{hom}||\bar\p Z^cw_{hom}|\\
&\ls\ve^2\w{x}^{-1}\w{s+|x|}^{-1}\w{s-|x|}^{-1},
\end{split}
\end{equation*}
which yields
\begin{equation}\label{energy:t7}
\|Q^{\alpha\beta\gamma}_{abc}\p^2_{\alpha\beta}Z^bw_{hom}\p_{\gamma}Z^cw_{hom}\|_{L^2(\cK)}\ls\ve^2\w{s}^{-3/2}.
\end{equation}
When $|b|\le N$ in the first term $I_{11}^{abc}$, it can be deduced from \eqref{null:structure}, \eqref{null:high}, \eqref{BA1} and \eqref{BA2} in the region $|x|\ge t/10$ that
\begin{equation}\label{energy:t8}
\begin{split}
&|Q^{\alpha\beta\gamma}_{abc}\p^2_{\alpha\beta}Z^bv\p_{\gamma}Z^cv|
\ls|\bar\p\p Z^bv||\p Z^cv|+|\p^2Z^bv||\bar\p Z^cv|\\
&\ls\ve_1\w{s}^{\dl_0-3/2}|\p Z^cv|+\ve_1\w{s}^{-1/2}\w{s-|x|}^{\dl_0-1}|\bar\p Z^cv|.
\end{split}
\end{equation}
While in the region $|x|\le s/10$, \eqref{hom:pw2} and \eqref{BA4} directly lead to
\begin{equation}\label{energy:t9}
\begin{split}
&\|\p^2_{\alpha\beta}Z^bv\p_{\gamma}Z^cv\|_{L^2(\cK\cap\{|x|\le s/10\})}\\
&\ls\sum_{j\le2N}\|\p\p_t^jv(s)\|_{L^2(\cK)}(\ve\w{s}^{-1}\ln^{-2}(2+s)+\ve_2^2\w{s}^{-1}).
\end{split}
\end{equation}
When $|b|\ge N+1$ in the first term $I_{11}^{abc}$, one can see that for $|c|\le N-1$,
\begin{equation*}
\begin{split}
&Q^{\alpha\beta\gamma}_{abc}e^q\p_tZ^av\p^2_{\alpha\beta}Z^bv\p_{\gamma}Z^cv\\
&=\p_{\alpha}(Q^{\alpha\beta\gamma}_{abc}e^q\p_tZ^av\p_{\beta}Z^bv\p_{\gamma}Z^cv)
-\p_t(Q^{\alpha\beta\gamma}_{abc}e^q\p_{\alpha}Z^av\p_{\beta}Z^bv\p_{\gamma}Z^cv)\\
&\quad-Q^{\alpha\beta\gamma}_{abc}e^q\p_tZ^av\p_{\beta}Z^bv\p^2_{\alpha\gamma}Z^cv
-Q^{\alpha\beta\gamma}_{abc}e^q\p_{\alpha}q\p_tZ^av\p_{\beta}Z^bv\p_{\gamma}Z^cv\\
&\quad+Q^{\alpha\beta\gamma}_{abc}\p_{\alpha}Z^av\p_t(e^q\p_{\beta}Z^bv\p_{\gamma}Z^cv).
\end{split}
\end{equation*}
Note that the terms in the last two lines above can be treated as for \eqref{energy:t8} and \eqref{energy:t9} with $\bar\p q=0$.
Then by the Cauchy-Schwarz inequality, we arrive at
\begin{equation}\label{energy:t10}
\begin{split}
&\Big|\int_0^t\int_{\cK}Q^{\alpha\beta\gamma}_{abc}e^q\p_tZ^av\p^2_{\alpha\beta}Z^bv\p_{\gamma}Z^cvdxds\Big|
\ls\tilde\delta\sum_{j\le2N}\int_0^t\int_{\cK}\frac{|\bar\p\p_t^jv|^2}{\w{s-|x|}^{3/2}}dxds\\
&\qquad+\tilde\delta^{-1}\sum_{j\le2N}\int_0^t\|\p\p_t^jv\|^2_{L^2(\cK)}(\ve_2^2\w{s}^{-1}+\ve_1\w{s}^{-1}\ln^{-2}(2+s))ds,
\end{split}
\end{equation}
where $\tilde\delta>0$ is a small constant. In addition,
the estimates of the other terms in $I_{11}^{abc}$ are analogous.

Therefore, collecting \eqref{energy:t6}-\eqref{energy:t10} yields
\begin{equation}\label{energy:t11}
\begin{split}
&\sum_{b+c\le a,b<a}\Big|\int_0^t\int_{\cK}I_1^{abc}dxds|\\
&\ls\|(\p w_{hom},\p v)(0)\|_{L^\infty(\cK)}\|\p\p_t^iv(0)\|^2_{L^2(\cK)}+\|(\p w_{hom},\p v)(t)\|_{L^\infty(\cK)}\|\p\p_t^iv(t)\|^2_{L^2(\cK)}\\
&\quad+\ve^2\sum_{j\le2N}\int_0^t\|\p\p_t^jv\|_{L^2(\cK)}\w{s}^{-3/2}ds
+\tilde\delta\sum_{j\le2N}\int_0^t\int_{\cK}\frac{|\bar\p\p_t^jv|^2}{\w{s-|x|}^{3/2}}dxds\\
&\quad+\tilde\delta^{-1}\sum_{j\le2N}\int_0^t\|\p\p_t^jv\|^2_{L^2(\cK)}(\ve_1\w{s}^{-1}\ln^{-2}(2+s)+\ve_2^2\w{s}^{-1})ds.
\end{split}
\end{equation}
Analogously for the treatment on the term $I_2^{abcd}$ in \eqref{energy:t5}, it can be concluded from \eqref{hom:pw1} and \eqref{BA1} that
\begin{equation}\label{energy:t12}
\begin{split}
&\sum_{b+c+d\le a,b<a}\Big|\int_0^t\int_{\cK}I_2^{abcd}dxds\Big|\\
&\ls\ve^3\sum_{j\le2N}\int_0^t\|\p\p_t^jv\|_{L^2(\cK)}\w{s}^{-1}ds+\ve_2^2\sum_{j\le2N}\int_0^t\|\p\p_t^jv\|^2_{L^2(\cK)}\w{s}^{-1}ds.
\end{split}
\end{equation}
In addition, the estimates of $I_1^a$ and $I_2^a$ are similar to \eqref{energy:t11} and \eqref{energy:t12}.

Substituting \eqref{hom:pw1}, \eqref{BA1}, \eqref{energy:t11} and \eqref{energy:t12} into \eqref{energy:t4} with small $\ve_1$ and $\tilde\delta$ derives
\begin{equation*}
\begin{split}
E(t)^2\ls E(0)^2+\int_0^tE(s)(\ve^2\w{s}^{-3/2}+\ve^3\w{s}^{-1})ds
+\int_0^tE(s)^2(\ve_1\w{s}^{-1}\ln^{-2}(2+s)+\ve_2^2\w{s}^{-1})ds,
\end{split}
\end{equation*}
where $\ds E(t)=\sup_{s\in[0,t]}\sum_{i\le2N}\|\p\p_t^iv(s)\|_{L^2(\cK)}$.
Thus, one has
\begin{equation*}
\begin{split}
&E(t)\ls E(0)+\ve^2+\ve^3\ln(2+t)+\int_0^tE(s)(\ve_1\w{s}^{-1}\ln^{-2}(2+s)+\ve_2^2\w{s}^{-1})ds.
\end{split}
\end{equation*}
This, together with \eqref{BA:time} and Lemma \ref{lem-Gronwall}, yields \eqref{energy:time}.
\end{proof}

\begin{lemma}\label{lem:energy:dtdx}
Under the assumptions of Theorem \ref{thm1}, let $v$ be defined by \eqref{err:def} and suppose that \eqref{BA1}-\eqref{BA4} hold.
Then we have that for $t\in [0,T_\ve]$,
\begin{equation}\label{energy:dtdx}
\sum_{|a|\le2N}\|\p\p^av(t)\|_{L^2(\cK)}\ls\ve_1.
\end{equation}
\end{lemma}
\begin{proof}
Although the proof is somewhat analogous to that of Lemma 6.2 in \cite{HouYinYuan24}, we still provide the details for the sake of completeness.

Set $\ds E_j(t)=\sum_{k=0}^{2N-j}\|\p\p_t^kv(t)\|_{H^j(\cK)}$ with $0\le j\le2N$.
Then one can find that for $j\ge1$,
\begin{equation}\label{energy:dtdx1}
\begin{split}
E_j(t)&\ls\sum_{k=0}^{2N-j}[\|\p\p_t^kv(t)\|_{L^2(\cK)}
+\sum_{1\le|a|\le j}(\|\p_t\p_t^k\p_x^av(t)\|_{L^2(\cK)}+\|\p_x\p_t^k\p_x^av(t)\|_{L^2(\cK)})]\\
&\ls E_0(t)+E_{j-1}(t)+\sum_{k=0}^{2N-j}\sum_{2\le|a|\le j+1}\|\p_t^k\p_x^av(t)\|_{L^2(\cK)},
\end{split}
\end{equation}
where we have used the fact of
\begin{equation*}
\sum_{k=0}^{2N-j}\sum_{1\le|a|\le j}\|\p_t\p_t^k\p_x^av(t)\|_{L^2(\cK)}
\ls\sum_{k=0}^{2N-j}\sum_{|a'|\le j-1}\|\p_x\p_t^{k+1}\p_x^{a'}v(t)\|_{L^2(\cK)}
\ls E_{j-1}(t).
\end{equation*}
For the last term in the second line of \eqref{energy:dtdx1}, it can be deduced from the elliptic estimate \eqref{ellip} that
\begin{equation}\label{energy:dtdx2}
\begin{split}
\|\p_t^k\p_x^av\|_{L^2(\cK)}&\ls\|\Delta_x\p_t^kv\|_{H^{|a|-2}(\cK)}+\|\p_t^kv\|_{H^{|a|-1}(\cK_{R+1})}\\
&\ls\|\p_t^{k+2}v\|_{H^{|a|-2}(\cK)}+\|\p_t^kQ(\p u,\p^2u)\|_{H^{|a|-2}(\cK)}\\
&\quad+\|\p_t^kv\|_{L^2(\cK_{R+1})}+\sum_{1\le|b|\le|a|-1}\|\p_t^k\p_x^bv\|_{L^2(\cK_{R+1})},
\end{split}
\end{equation}
where $\Delta v=\p_t^2v-\Box v=\p_t^2v-Q(\p u,\p^2u)$ has been used.
In addition, by \eqref{hom:pw1}, \eqref{BA1} with $|a|\le j+1$, $k\le2N-j$, we have
\begin{equation}\label{energy:dtdx3}
\begin{split}
\|\p_t^kQ(\p u,\p^2u)\|_{H^{|a|-2}(\cK)}
&\ls\ve^2+\sum_{l\le k}\ve_1(\|\p_t^l\p v\|_{H^{|a|-1}(\cK)}+\|\p_t^{l+2}v\|_{H^{|a|-2}(\cK)})\\
&\ls\ve^2+\ve_1[E_j(t)+E_{j-1}(t)].
\end{split}
\end{equation}
Collecting \eqref{energy:dtdx1}-\eqref{energy:dtdx3} yields that for $j\ge1$,
\begin{equation*}
E_j(t)\ls\ve^2+E_0(t)+E_{j-1}(t)+\ve_1E_j(t)+\|v\|_{L^2(\cK_{R+1})}.
\end{equation*}
Together with \eqref{BA3}, \eqref{energy:time} and the smallness of $\ve_1$, this completes the proof of \eqref{energy:dtdx}.
\end{proof}

\begin{lemma}\label{lem:energy:Z}
Under the assumptions of Theorem \ref{thm1}, let $v$ be defined by \eqref{err:def} and suppose that \eqref{BA1}-\eqref{BA4} hold.
Then we it holds that for $t\in [0,T_\ve]$,
\begin{equation}\label{energyA}
\sum_{|a|\le2N-1}\|\p Z^av(t)\|_{L^2(\cK)}\ls\ve_1(1+t)^{1/2}.
\end{equation}
\end{lemma}
\begin{proof}
Analogously to Lemma \ref{lem:energy:time}, we can obtain that for any multi-index $a\in\N_0^4$ with $|a|\le2N-1$,
\begin{equation}\label{energyA1}
\begin{split}
&\sum_{|a|\le2N-1}\|\p Z^av(t)\|^2_{L^2(\cK)}\ls\sum_{|a|\le2N-1}\|\p Z^av(0)\|^2_{L^2(\cK)}
+\sum_{|a|\le2N-1}(|\cB_1^a|+|\cB_2^a|+|\cB_3^a|)\\
&\quad+\sum_{|a|\le2N-1}\int_0^t\|\p Z^av(s)\|_{L^2(\cK)}(\ve^2\w{s}^{-3/2}+\ve^3\w{s}^{-1})ds\\
&\quad+\sum_{|a|\le2N-1}\int_0^t\|\p Z^av(s)\|^2_{L^2(\cK)}(\ve_1\w{s}^{-1}\ln^{-2}(2+s)+\ve_2^2\w{s}^{-1})ds,\\
\end{split}
\end{equation}
where
\begin{equation*}\label{energyA1}
\begin{split}
&\cB^a_1=\int_0^t\int_{\p\cK}e^q(n_1(x),n_2(x))\cdot\nabla Z^av(s,x)\p_tZ^av(s,x)d\sigma ds,\\
&\cB^a_2=\int_0^t\int_{\p\cK}\sum_{j=1}^2\sum_{\beta,\gamma=0}^2
Q^{j\beta\gamma}n_j(x)e^q\p_tZ^av\p_{\beta}Z^av\p_{\gamma}ud\sigma ds,\\
&\cB^a_3=\int_0^t\int_{\p\cK}\sum_{j=1}^2\sum_{\beta,\gamma,\delta=0}^2
Q^{j\beta\gamma\delta}n_j(x)e^q\p_tZ^av\p_{\beta}Z^av\p_{\gamma}u\p_{\delta}ud\sigma ds,
\end{split}
\end{equation*}
$(n_1(x),n_2(x))$ is the unit outer normal of the boundary $\p\cK$ and $d\sigma$ is the curve measure on the boundary $\p\cK$.
According to $\p\cK\subset\overline{\cK_1}$ and the trace theorem, we arrive at
\begin{equation}\label{energyA2}
\begin{split}
|\cB^a_1|&\ls\sum_{|b|\le|a|}\int_0^t\|(1-\chi_{[1,2]}(x))\p Z^av\|^2_{H^1(\cK)}ds\\
&\ls\sum_{|b|\le|a|+1}\int_0^t\|\p \p^bv\|^2_{L^2(\cK_2)}ds\ls\ve_1^2(1+t),
\end{split}
\end{equation}
where \eqref{energy:dtdx} has been used.
Similarly, $\cB^a_2$ and $\cB^a_3$ can be also treated.
Substituting \eqref{energyA2} into \eqref{energyA1} yields
\begin{equation}\label{energyA3}
\begin{split}
\tilde E(t)^2&\ls\tilde E(0)^2+\ve_1^2(1+t)+\int_0^t\tilde E(s)(\ve^2\w{s}^{-3/2}+\ve^3\w{s}^{-1})ds\\
&\quad+\int_0^t\tilde E(s)^2(\ve_1\w{s}^{-1}\ln^{-2}(2+s)+\ve_2^2\w{s}^{-1})ds,
\end{split}
\end{equation}
where $\ds\tilde E(t)=\sum_{|a|\le2N-1}\|\p Z^av(t)\|_{L^2(\cK)}$.
If $\ds\tilde E(t)\ge\ve_1(1+t)^{1/2}$, otherwise \eqref{energyA} is directly obtained, \eqref{energyA3} can be reduced to
\begin{equation*}
\begin{split}
\tilde E(t)&\ls\tilde E(0)+\ve_1(1+t)^{1/2}+\int_0^t(\ve^2\w{s}^{-3/2}+\ve^3\w{s}^{-1})ds\\
&\quad+\int_0^t\tilde E(s)(\ve_1\w{s}^{-1}\ln^{-2}(2+s)+\ve_2^2\w{s}^{-1})ds\\
&\ls\ve_1(1+t)^{1/2}+\int_0^t\tilde E(s)(\ve_1\w{s}^{-1}\ln^{-2}(2+s)+\ve_2^2\w{s}^{-1})ds.
\end{split}
\end{equation*}
Thus, \eqref{energyA} is achieved by Lemma \ref{lem-Gronwall} and \eqref{BA:time}.
\end{proof}

\subsection{Local energy decay estimates}

It is noted that the improvement on the energy estimate in \eqref{energyA} mainly depends on the estimate of the boundary term in \eqref{energyA2}, which can be enhanced by the local energy decay estimates.
\begin{lemma}
Under the assumptions of Theorem \ref{thm1}, let $v$ be defined by \eqref{err:def} and suppose that \eqref{BA1}-\eqref{BA4} hold.
Then we have that for $t\in [0,T_\ve]$,
\begin{equation}\label{loc:energyA}
\|\p_x^{\le1}v(t)\|_{L^2(\cK_R)}\ls\ve_1^2\w{t}^{-1}.
\end{equation}
\end{lemma}
\begin{proof}
The proof of \eqref{loc:energyA} is analogous to that of \eqref{dt:loc} with the good structure in \eqref{eqn:null:div}.
As in \cite{HouYinYuan24}, we firstly eliminate the quadratic $Q_0$ term in \eqref{eqn:null:div} with $a=0$ by introducing the good known
\begin{equation}\label{loc:energyA1}
V_0=v-\sum_{\mu=0}^2\frac12C^{0,\mu}(\omega)u\p_{\mu}u
\end{equation}
and by using the identity for any functions $f,g$
\begin{equation}\label{Q0:ident}
\Box(fg)=f\Box g+g\Box f+2Q_0(f,g),
\end{equation}
which can be verified by direct computation.
Thus, one has that $V_0|_{[0,\infty)\times\p\cK}=0$ and
\begin{equation}\label{loc:energyA2}
\begin{split}
\Box V_0&=\p_\theta(\frac{C(\omega)u\p^2u}{|x|})+\sum_{\mu=0}^2\p_{\mu}(\frac{C(\omega)uZ\p u}{|x|})+\frac{C(\omega)uZ\p u}{|x|^2}
+\sum_{\mu=0}^2[\Delta,C^{0,\mu}(\omega)]u\p_{\mu}u\\
&\quad+\sum_{\mu=0}^2C(\omega)(\p_{\mu}u\Box u+u\Box\p_{\mu}u)
+\frac13\p_t\big[\sum_{\alpha,\beta,\gamma,\delta=0}^2Q^{\alpha\beta\gamma\delta}\omega_{\alpha}\omega_{\beta}\omega_{\gamma}\omega_{\delta}(\p_tu)^3\big]\\
&\quad+\sum_{\alpha,\beta,\gamma,\delta=0}^2Q^{\alpha\beta\gamma\delta}[\p^2_{\alpha\beta}u\p_{\gamma}u\p_{\delta}u
-\omega_{\alpha}\omega_{\beta}\omega_{\gamma}\omega_{\delta}\p_t^2u(\p_tu)^2].
\end{split}
\end{equation}
By substituting the expressions of $\Box u$ and $\Box\p_{\mu}u$ into \eqref{loc:energyA2}, we have
\begin{equation}\label{loc:energyA3}
\begin{split}
&\Box V_0=\p_tF^0+\sum_{\Xi=0,1,2,3,\theta}\p_{\Xi}G^{\Xi}+G_1^r+G_2^r,\\
&F^0(t,x)=\frac13\eta(t,x)\sum_{\alpha,\beta,\gamma,\delta=0}^2Q^{\alpha\beta\gamma\delta}\omega_{\alpha}\omega_{\beta}\omega_{\gamma}
\omega_{\delta}(\p_tu)^3\\
&G^{\alpha}(t,x)=\frac{C(\omega)uZ\p u}{|x|},\quad\alpha=0,1,2,\quad G^{\theta}(t,x):=\frac{C(\omega)u\p^2u}{|x|},\\
&G_1^r(t,x)=\frac13\p_t[(1-\eta(t,x))\sum_{\alpha,\beta,\gamma,\delta=0}^2Q^{\alpha\beta\gamma\delta}\omega_{\alpha}
\omega_{\beta}\omega_{\gamma}\omega_{\delta}(\p_tu)^3]\\
&\quad+\frac{C(\omega)uZ\p u}{|x|^2}
+\sum_{\alpha,\beta,\gamma,\mu=0}^2C(\omega)Q^{\alpha\beta\gamma}(\p_{\mu}u\p^2_{\alpha\beta}u\p_{\gamma}u
+u\p_{\mu}(\p^2_{\alpha\beta}u\p_{\gamma}u))\\
&\quad+\sum_{\alpha,\beta,\gamma,\delta=0}^2Q^{\alpha\beta\gamma\delta}[\p^2_{\alpha\beta}u\p_{\gamma}u\p_{\delta}u
-\omega_{\alpha}\omega_{\beta}\omega_{\gamma}\omega_{\delta}\p_t^2u(\p_tu)^2]
+\frac12\sum_{\mu=0}^2[\Delta,C^{0,\mu}(\omega)]u\p_{\mu}u,\\
&G_2^r(t,x)=\sum_{\alpha,\beta,\gamma,\delta,\mu=0}^2C(\omega)Q^{\alpha\beta\gamma\delta}
(\p_{\mu}u\p^2_{\alpha\beta}u\p_{\gamma}u\p_{\delta}u+u\p_{\mu}(\p^2_{\alpha\beta}u\p_{\gamma}u\p_{\delta}u)),
\end{split}
\end{equation}
where
\begin{equation}\label{cutoff:cone}
\eta(t,x)=\chi_{[\frac{9}{10},\frac{19}{20}]}(\frac{|x|}{t})
\end{equation}
is a conic cutoff function and is supported in $\{|x|\ge\frac{9t}{10}\}$.
For the first term in the last line of $G_1^r(t,x)$ in \eqref{loc:energyA3}, similarly to the proof of \eqref{null:structure}, one has
that for any $|a|\le1$,
\begin{equation}\label{cubic:null}
\begin{split}
&|Z^a(\p^2_{\alpha\beta}u\p_{\gamma}u\p_{\delta}u-\omega_{\alpha}\omega_{\beta}\omega_{\gamma}\omega_{\delta}\p_t^2u(\p_tu)^2)|\\
&\ls\sum_{b+c+d\le a}(|\bar\p\p Z^bu||\p Z^cu||\p Z^du|+|\p^2Z^bu||\bar\p Z^cu||\p Z^du|+\frac{1}{|x|}|\p^2Z^bu||\p Z^cu||\p Z^du|)\\
&\ls\ve_1^3\w{y}^{-3/2}\w{s+|y|}^{\delta_0-1}\w{s-|y|}^{2\delta_0-2},
\end{split}
\end{equation}
where \eqref{hom:pw1}, \eqref{BA1} and \eqref{BA2} have been used.
For the last term in $G_1^r(t,x)$, by $\Delta=\p_r^2+\frac{1}{r}\p_r+\frac{\Omega^2}{r^2}$ and $\p_r\omega=0$,
we have that for any function $f$,
\begin{equation}\label{commut}
|[\Delta,C^{0,\mu}(\omega)]f|\ls|x|^{-2}|Z^{\le1}f|.
\end{equation}
On the other hand, it can be concluded from \eqref{1st:null:condition}, \eqref{null:structure}, \eqref{hom:pw1}, \eqref{BA1}-\eqref{BA3} with $N\ge53$ that
\begin{equation}\label{loc:energyA4}
\begin{split}
\sum_{|a|\le5}|\p^aF^0(s,y)|
&\ls\ve_1^3\w{s+|y|}^{-3/2}\w{s-|y|}^{-2},\\
\sum_{|a|\le2}|Z^aG^{\Xi}(s,y)|
&\ls\ve_1^2\w{y}^{-3/2}\w{s+|y|}^{-1/2}\w{s-|y|}^{4\dl_0-3/2}, \quad\Xi=0,1,2,\theta,\\
\sum_{|a|\le1}|\p^aG_1^r(s,y)|
&\ls\ve_1^3\w{y}^{-3/2}\w{s+|y|}^{3\dl_0-3}+\ve_1^2\w{y}^{-5/2}\w{s+|y|}^{-1/2}\w{s-|y|}^{4\dl_0-3/2}\\
&\quad+\ve_1^3\w{y}^{-1}\w{s+|y|}^{\dl_0-3/2}\w{s-|y|}^{4\dl_0-3/2},\\
\sum_{|a|\le1}|\p^aG_2^r(s,y)|
&\ls\ve_1^4\w{y}^{-2}\w{s-|y|}^{6\dl_0-7/2}.
\end{split}
\end{equation}
Let $w_1^b$ and $w_2^b$ be the solutions of
\begin{equation}\label{loc:energyA5}
\left\{
\begin{aligned}
&\Box w_1^b=(1-\chi_{[2,3]}(x))[\p_tF^0+\sum_{\Xi=0,1,2,\theta}\p_{\Xi}G^{\Xi}+G_1^r+G_2^r]-\sum_{i=1}^2G^i(\p_i\chi_{[2,3]}(x)),\\
&\Box w_2^b=\p_t(\chi_{[2,3]}(x)F^0)+\sum_{\Xi=0,1,2,\theta}\p_{\Xi}(\chi_{[2,3]}(x)G^{\Xi})+\chi_{[2,3]}(x)(G_1^r+G_2^r),\\
&w_1^b|_{\p\cK}=w_2^b|_{\p\cK}=0,\\
&(w_1^b,\p_tw_1^b)(0,x)=(V_0,\p_tV_0)(0,x),\quad (w_2^b,\p_tw_2^b)(0,x)=(0,0),
\end{aligned}
\right.
\end{equation}
respectively.
It follows from \eqref{loc:energyA3} and \eqref{loc:energyA5} that $V_0=w_1^b+w_2^b$.
By \eqref{loc:decay} with $m=1$ and $\eta=1$, one has
\begin{equation}\label{loc:energyA6}
\w{t}\|\p_x^{\le1}w_1^b(t,x)\|_{L^2(\cK_R)}
\ls\ve^2+\sup_{s\in[0,t]}\w{s}\|\Box w_1^b(s)\|_{L^2(\cK_3)}\ls\ve_1^2,
\end{equation}
where \eqref{initial:data}, \eqref{loc:energyA1} are \eqref{loc:energyA4} are used.

Next, we treat the estimates of $w_2^b$. Let $w^c$, $w_{\Xi}^c$ $(\Xi=0,1,2,\theta)$, $w_{r,1}^c$,
$w_{r,2}^c$ be the solutions of the Cauchy problem
\begin{equation}\label{loc:energyA7}
\Box w^c=\chi_{[2,3]}(x)F^0,\quad \Box w_{\Xi}^c=\chi_{[2,3]}(x)G^{\Xi}\quad (\Xi=0,1,2,\theta),\quad \Box w_{r,l}^c=\chi_{[2,3]}(x)G_l^r,l=1,2,
\end{equation}
with zero initial data, respectively.
Then $\ds v_2^b=w_2^b-\chi_{[1,2]}(x)(\p_tw^c+\sum_{\Xi=0,1,2,\theta}\p_{\Xi}w_{\Xi}^c+\sum_{l=1}^2w_{r,l}^c)$ is the solution of the IBVP
\begin{equation*}
\left\{
\begin{aligned}
&\Box v_2^b=[\Delta,\chi_{[1,2]}]\big(\p_tw^c+\sum_{\Xi=0,1,2,\theta}\p_{\Xi}w_{\Xi}^c+\sum_{l=1}^2w_{r,l}^c\big),\\
&v_2^b|_{\p\cK}=0,\\
&(v_2^b,\p_tv_2^b)(0,x)=(0,\chi_{[2,3]}(x)(F^0(0,x)+G^0(0,x))).
\end{aligned}
\right.
\end{equation*}
Since $\Box v_2^b$ has compact support in the space variable, and the estimate of $v_2^b$ is analogous to that of $w_1^b$ in \eqref{loc:energyA6},
one then has
\begin{equation}\label{loc:energyA8}
\begin{split}
\w{t}\|\p_x^{\le1}v_2^b(t,x)\|_{L^2(\cK_R)}
&\ls\ve_1^2+\sup_{s\in[0,t]}\w{s}(\|\p_x^{\le1}\p_tw^c(s)\|_{L^\infty(|x|\le2)}\\
&\quad+\sum_{\Xi=0,1,2,\theta}\|\p_x^{\le1}\p w_{\Xi}^c(s)\|_{L^\infty(|x|\le2)}
+\sum_{l=1}^2\|\p_x^{\le1}w_{r,l}^c(s)\|_{L^\infty(|x|\le2)}).
\end{split}
\end{equation}
Applying \eqref{pw:ivp2}, \eqref{pw:ivp3}, \eqref{dpw:ivp} and Lemma \ref{lem:dt:impr} with $\mu=\nu=\dl_0$ to $w_{r,1}^c$, $w_{r,2}^c$, $w_{\Xi}^c$ and $w^c$ in \eqref{loc:energyA7}, respectively, yields
\begin{equation}\label{loc:energyA9}
\begin{split}
\w{t+|x|}^{1/2}\w{t-|x}^{1/2}|\p_x^{\le1}w_{r,1}^c(t,x)|&\ls\ve_1^2,\\
\w{t+|x|}^{1/2}\w{t-|x}^{1/2}|\p_x^{\le1}w_{r,2}^c(t,x)|&\ls\ve_1^4\ln(2+t),\\
\w{t-|x}^{1+\dl_0}|\p_x^{\le1}\p w_{\Xi}^c(t,x)|&\ls\ve_1^2,\quad\Xi=0,1,2,\theta,\\
\w{t}\|\p_x^{\le1}\p_tw^c(t,x)\|_{L^\infty(|x|\le t/10)}&\ls\ve_1^3,
\end{split}
\end{equation}
where \eqref{loc:energyA4} is also used.
Collecting \eqref{loc:energyA1}, \eqref{loc:energyA6}, \eqref{loc:energyA8} and \eqref{loc:energyA9} with $V_0=w_1^b+w_2^b$, $\ds w_2^b=v_2^b+\chi_{[1,2]}(x)(\p_tw^c+\sum_{\Xi=0,1,2,\theta}\p_{\Xi}w_{\Xi}^c+\sum_{l=1}^2w_{r,l}^c)$, \eqref{hom:pw1}, \eqref{BA3} and \eqref{BA:time} finishes the proof of \eqref{loc:energyA}.
\end{proof}

\begin{lemma}
Under the assumptions of Theorem \ref{thm1}, let $v$ be defined by \eqref{err:def} and suppose that \eqref{BA1}-\eqref{BA4} hold.
Then we have that for $t\in [0,T_\ve]$,
\begin{equation}\label{loc:energyB}
\sum_{|a|\le2N-5}\|\p^av(t)\|_{L^2(\cK_R)}\ls\ve_1^2\w{t}^{4\dl_0-1/2}.
\end{equation}
\end{lemma}
\begin{proof}
Choosing $Z^a=\p_t^i$ in \eqref{eqn:null} yields
\begin{equation*}
\begin{split}
\Box\p_t^iv&=\sum_{j+k\le i}\Big\{\sum_{\mu,\nu=0}^2|x|^{-1}(C(\omega)\p^2_{\mu\nu}\p_t^ju\Omega\p_t^ku
+C(\omega)\Omega\p_{\mu}\p_t^ju\p_{\nu}\p_t^ku)\\
&\quad+\sum_{\mu=0}^2\Big[C_{ijk}^{12,\mu}Q_{12}(\p_{\mu}\p_t^ju,\p_t^ku)
+C_{ijk}^{0,\mu}(\omega)Q_0(\p_{\mu}\p_t^ju,\p_t^ku)+C(\omega)\Box\p_t^ju\p_{\mu}\p_t^k u\Big]\Big\}\\
&\quad+\p_t^i\big[\sum_{\alpha,\beta,\gamma,\delta=0}^2Q^{\alpha\beta\gamma\delta}\p^2_{\alpha\beta}u\p_{\gamma}u\p_{\delta}u\big].
\end{split}
\end{equation*}
As in the definition \eqref{loc:energyA1}, introducing the good unknown 
\begin{equation}\label{Vi:def}
V_i=\p_t^iv-\frac12\sum_{\mu=0}^2\sum_{j+k\le i}C_{ijk}^{0,\mu}(\omega)\p_{\mu}\p_t^ju\p_t^ku.
\end{equation}
Then one has $V_i|_{[0,\infty)\times\p\cK}=0$ and
\begin{equation}\label{loc:energyB1}
\begin{split}
\Box V_i&=\sum_{j+k\le i}\Big\{\sum_{\mu,\nu=0}^2|x|^{-1}(C(\omega)\p^2_{\mu\nu}\p_t^ju\Omega\p_t^ku
+C(\omega)\Omega\p_{\mu}\p_t^ju\p_{\nu}\p_t^ku)\\
&\quad+\sum_{\mu=0}^2\Big[\frac12\sum_{\mu=0}^2[\Delta,C_{ijk}^{0,\mu}(\omega)]\p_{\mu}\p_t^ju\p_t^ku
+C_{ijk}^{12,\mu}Q_{12}(\p_{\mu}\p_t^ju,\p_t^ku)\\
&\qquad+C(\omega)\p_{\mu}\p_t^ju\Box\p_t^ku+C(\omega)\p_t^ku\Box\p_{\mu}\p_t^ju\Big]\Big\}
+\p_t^i\big[\sum_{\alpha,\beta,\gamma,\delta=0}^2Q^{\alpha\beta\gamma\delta}\p^2_{\alpha\beta}u\p_{\gamma}u\p_{\delta}u\big].
\end{split}
\end{equation}
Applying \eqref{dt:loc} to \eqref{loc:energyB1} with $F_1=\Box V_i$, $F_0\equiv0$, $i\le2N-7$, $\mu=\dl_0$ and \eqref{initial:data} yields
\begin{equation}\label{loc:energyB2}
\begin{split}
\w{t}\|\p_x^{\le1}\p_tV_i\|_{L^2(\cK_{R_1})}\ls\ve^2+\sum_{|a|\le1}\sup_{s\in[0,t]}\w{s}\|\p^aF_1(s)\|_{L^2(\cK_3)}\\
+\sum_{|a|\le2}\sup_{(s,y)\in[0,t]\times\overline{\R^2\setminus\cK_2}}\w{y}^{1/2}\cW_{1+\dl_0,1}(s,y)|Z^aF_1(s,y)|.
\end{split}
\end{equation}
Note that \eqref{Sobo:ineq} and \eqref{energyA} lead to
\begin{equation}\label{loc:energyB3}
\sum_{|a|\le2N-3}|\p Z^av|\ls\ve_1\w{x}^{-1/2}(1+t)^{1/2}.
\end{equation}
It follows from \eqref{Q12:null:structure}, \eqref{hom:pw1}, \eqref{BA1}, \eqref{BA3}, \eqref{commut}
and \eqref{loc:energyB3} that
\begin{equation}\label{loc:energyB4}
\begin{split}
\sum_{|a|\le2}|Z^aF_1(s,y)|&\ls\ve_1^2\w{y}^{-3/2}\w{s+|y|}^{-1/2}\w{s-|y|}^{3\dl_0-1/2}(1+s)^{1/2}\\
&\quad+\ve_1^3\w{y}^{-1}\w{s+|y|}^{-1/2}\w{s-|y|}^{4\dl_0-3/2}(1+s)^{1/2}.
\end{split}
\end{equation}
On the other hand, by \eqref{hom:pw1}, \eqref{BA3} and \eqref{energy:dtdx}, one has
\begin{equation}\label{loc:energyB5}
\sum_{|a|\le1}\|\p^aF_1(s)\|_{L^2(\cK_3)}\ls\ve_1^2\w{s}^{3\delta_0-1}.
\end{equation}
Therefore, collecting \eqref{hom:pw1}, \eqref{BA3}, \eqref{energy:dtdx}, \eqref{loc:energyA}, \eqref{loc:energyB2}, \eqref{loc:energyB4} and \eqref{loc:energyB5} with the definition \eqref{Vi:def} yields
\begin{equation}\label{loc:energyB6}
\sum_{i\le2N-7}\|\p^{\le1}\p_t\p_t^iv\|_{L^2(\cK_{R_1})}\ls\ve_1^2\w{t}^{4\dl_0-1/2}.
\end{equation}
The remaining proof is similar to that of Lemma 6.6 in \cite{HouYinYuan24}.
For $j=0,1,2,\cdots,2N-5$, denote $\ds E_j^{loc}(t)=\sum_{k\le j}\|\p_t^k\p_x^{2N-5-j}v\|_{L^2(\cK_{R+j})}$.
Then \eqref{loc:energyA} and \eqref{loc:energyB6} imply
\begin{equation}\label{loc:energyB7}
E_{2N-5}^{loc}(t)+E_{2N-6}^{loc}(t)\ls\ve_1^2\w{t}^{4\dl_0-1/2}.
\end{equation}
When $j\le2N-7$, we have $\p_x^{2N-5-j}=\p_x^2\p_x^{2N-7-j}$.
Thus, one can apply the elliptic estimate \eqref{ellip} to $(1-\chi_{[R+j,R+j+1]})\p_t^{k}v$
to obtain
\begin{equation}\label{loc:energyB8}
\begin{split}
E_j^{loc}(t)&\ls\sum_{k\le j}\|\p_x^{2N-5-j}(1-\chi_{[R+j,R+j+1]})\p_t^kv\|_{L^2(\cK)}\\
&\ls\sum_{k\le j}[\|\Delta(1-\chi_{[R+j,R+j+1]})\p_t^kv\|_{H^{2N-7-j}(\cK)}\\
&\quad+\sum_{k\le j}\|(1-\chi_{[R+j,R+j+1]})\p_t^kv\|_{H^{2N-6-j}(\cK)}]\\
&\ls\sum_{k\le j}\|(1-\chi_{[R+j,R+j+1]})\Delta\p_t^kv\|_{H^{2N-7-j}(\cK)}+\sum_{j+1\le l\le2N-5}E_l^{loc}(t)\\
&\ls\sum_{k\le j}[\|\Box\p_t^kv\|_{H^{2N-7-j}(\cK_{R+j+1})}+\|\p_t^{k+2}v\|_{H^{2N-7-j}(\cK_{R+j+1})}]\\
&\quad+\sum_{j+1\le l\le2N-5}E_l^{loc}(t),
\end{split}
\end{equation}
where the fact of $\Delta=\p_t^2-\Box$ has been used.
Analogously to \eqref{loc:energyB5}, it can be concluded from $u=v+w_{hom}$, \eqref{hom:pw1}, \eqref{BA1} and \eqref{energy:dtdx}
that for $j\le2N-7$,
\begin{equation}\label{loc:energyB9}
\sum_{k\le j}\|\Box\p_t^kv\|_{H^{2N-7-j}(\cK_{R+j+1})}\ls\sum_{k\le j}\|\p_t^kQ(\p u,\p^2u)\|_{H^{2N-7-j}(\cK_{R+j+1})}\ls\ve_1^2(1+t)^{\delta_0-1}.
\end{equation}
Thus, \eqref{loc:energyB8} and \eqref{loc:energyB9} leads to that for $j\le2N-7$,
\begin{equation*}
E_j^{loc}(t)\ls\ve_1^2(1+t)^{\delta_0-1}+\sum_{j+1\le l\le2N-5}E_l^{loc}(t).
\end{equation*}
This, together with \eqref{loc:energyB7}, yields \eqref{loc:energyB}.
\end{proof}

\subsection{Improved energy estimates}

With the local energy decay estimate \eqref{loc:energyB},
the energy estimates including all the vector fields in \eqref{energyA} will be improved.
\begin{lemma}
Under the assumptions of Theorem \ref{thm1}, let $v$ be defined by \eqref{err:def} and suppose that \eqref{BA1}-\eqref{BA4} hold.
Then we have that for $t\in [0,T_\ve]$,
\begin{equation}\label{energyB}
\sum_{|a|\le2N-7}\|\p Z^av(t)\|_{L^2(\cK)}\ls\ve_1(1+t)^{4\dl_0}.
\end{equation}
\end{lemma}
\begin{proof}
Analogously to Lemma \ref{lem:energy:Z}, one can obtain the following better estimates on the boundary terms $\cB^a_1,\cB^a_2,\cB^a_3$ for $|a|\le2N-7$,
\begin{equation}\label{energyB1}
\sum_{l=1}^3|\cB^a_l|\ls\sum_{|b|\le|a|+1}\int_0^t\|\p \p^bv\|^2_{L^2(\cK_2)}ds
\ls\ve_1^2\int_0^t(1+s)^{8\dl_0-1}ds\ls\ve_1^2(1+t)^{8\dl_0},
\end{equation}
where \eqref{loc:energyB} has been used.
Then \eqref{energyB} is achieved.
\end{proof}

In addition, the local energy decay estimate \eqref{loc:energyB} can be further improved as follows.
\begin{lemma}
Under the assumptions of Theorem \ref{thm1}, let $v$ be defined by \eqref{err:def} and suppose that \eqref{BA1}-\eqref{BA4} hold.
Then it holds that for $t\in [0,T_\ve]$,
\begin{equation}\label{loc:energyC}
\sum_{|a|\le2N-11}\|\p^av(t)\|_{L^2(\cK_R)}\ls\ve_1^2\w{t}^{8\dl_0-1}.
\end{equation}
\end{lemma}
\begin{proof}
The proof is analogous to that of \eqref{loc:energyB}.
At first, \eqref{Sobo:ineq} and \eqref{energyB} imply
\begin{equation}\label{loc:energyC1}
\sum_{|a|\le2N-9}|\p Z^av|\ls\ve_1\w{x}^{-1/2}(1+t)^{4\dl_0}.
\end{equation}
Then, with \eqref{loc:energyC1} instead of \eqref{loc:energyB3}, \eqref{loc:energyB4} can be improved to
\begin{equation}\label{loc:energyC2}
\begin{split}
\sum_{|a|\le2}|Z^aF_1(s,y)|&\ls\ve_1^2\w{y}^{-3/2}\w{s+|y|}^{-1/2}\w{s-|y|}^{3\dl_0-1/2}(1+s)^{4\dl_0}\\
&+\ve_1^3\w{y}^{-1}\w{s+|y|}^{-1/2}\w{s-|y|}^{4\dl_0-3/2}(1+s)^{4\dl_0}.
\end{split}
\end{equation}
Substituting \eqref{loc:energyB5} and \eqref{loc:energyC2} into \eqref{loc:energyB2} together with \eqref{hom:pw1}, \eqref{BA3}, \eqref{loc:energyB}, \eqref{Vi:def} and $i\le2N-13$, yields
\begin{equation}\label{loc:energyC3}
\sum_{i\le2N-13}\|\p^{\le1}\p_t\p_t^iv\|_{L^2(\cK_{R_1})}\ls\ve_1^2\w{t}^{8\dl_0-1}.
\end{equation}
Therefore, \eqref{loc:energyC} is achieved by \eqref{loc:energyC3} with an analogous argument as in the proof of \eqref{loc:energyB}.
\end{proof}

\begin{lemma}
Under the assumptions of Theorem \ref{thm1}, let $v$ be defined by \eqref{err:def} and suppose that \eqref{BA1}-\eqref{BA4} hold.
Then we have that for $t\in [0,T_\ve]$,
\begin{equation}\label{energyC}
\sum_{|a|\le2N-13}\|\p Z^av(t)\|_{L^2(\cK)}\ls\ve_1.
\end{equation}
\end{lemma}
\begin{proof}
Since the proof is similar to that of \eqref{energyB} with \eqref{loc:energyC} instead of \eqref{loc:energyB}, we omit the details here.
\end{proof}

\section{Improved pointwise estimates and proof of Theorem \ref{thm1}}\label{sect6}

In this section, we will improve the pointwise estimates \eqref{BA1}-\eqref{BA4} and close the bootstrap assumptions in Section \ref{sect5}.
Although these estimates are analogous to the ones in Section 7 of \cite{HouYinYuan24} with some minor modifications, we still give the details
for readers' convenience.

\subsection{Decay estimates on the good derivatives of solutions}

\begin{lemma}\label{lem:BA3:A}
Under the assumptions of Theorem \ref{thm1}, let $v$ be defined by \eqref{err:def} and suppose that \eqref{BA1}-\eqref{BA4} hold.
Then it holds that for $t\in [0,T_\ve]$,
\begin{equation}\label{BA3:impr:A}
\sum_{|a|\le2N-16}|Z^av|\ls\ve_1\w{t+|x|}^{6\dl_0}.
\end{equation}
\end{lemma}
\begin{proof}
Until now, we have not obtained the estimates of $Z^au$ with higher order derivatives by use of \eqref{eqn:null}.
To achieve this, we will employ the primitive equation \eqref{eqn:high}.
In order to preserve the zero boundary conditions, set
\begin{equation}\label{cutoff:Z}
\tilde Z=\chi_{[1/2,1]}(x) Z,\quad Z\in\{\p,\Omega\}.
\end{equation}
Applying \eqref{pw:ibvp} to $\Box\tilde Z^av=\Box(\tilde Z^a-Z^a)v+\Box Z^av$ with $F_0(t,x)\equiv0$, $\mu=\nu=\dl_0/2$ and $|a|\le2N-15$ yields
\begin{equation}\label{BA3:impr:A1}
\begin{split}
&\w{t+|x|}^{1/2}|\tilde Z^av|\ls\ve^2+\sum_{|b|\le4}\sup_{s\in[0,t]}\w{s}^{1/2+\dl_0}\|\p^b\Box(\tilde Z^a-Z^a)v(s)\|_{L^2(\cK_4)}\\
&\qquad+\sum_{|b|\le5}\sup_{(s,y)\in[0,t]\times\cK}\w{y}^{1/2}\cW_{1+\dl_0,1}(s,y)|\p^b\Box Z^av(s,y)|,
\end{split}
\end{equation}
where we have used \eqref{initial:data} and the fact of $\tilde Z=Z$ in the region $|x|\ge1$.
In addition, it can be deduced from \eqref{loc:energyB} that
\begin{equation}\label{BA3:impr:A2}
\sum_{|a|\le2N-15}\sum_{|b|\le4}\sup_{s\in[0,t]}\w{s}^{1/2+\dl_0}\|\p^b\Box(\tilde Z^a-Z^a)v(s)\|_{L^2(\cK_4)}
\ls\ve_1^2\w{t}^{5\dl_0}.
\end{equation}
For the second line of \eqref{BA3:impr:A1}, one can conclude from \eqref{eqn:high}, \eqref{hom:pw1}, \eqref{BA1} and \eqref{loc:energyC1} that
\begin{equation}\label{BA3:impr:A3}
\sum_{|a|\le2N-15}\sum_{|b|\le5}|Z^b\Box Z^av(s,y)|
\ls\ve_1^2\w{y}^{-1}\w{s-|y|}^{\dl_0-1}(1+s)^{4\dl_0}.
\end{equation}
Substituting \eqref{BA3:impr:A2} and \eqref{BA3:impr:A3} into \eqref{BA3:impr:A1} leads to
\begin{equation*}
\sum_{|a|\le2N-15}|\tilde Z^av|\ls\ve_1\w{t+|x|}^{6\dl_0}.
\end{equation*}
This, together with \eqref{loc:energyB} and the Sobolev embedding, completes the proof of \eqref{BA3:impr:A}.
\end{proof}

\begin{lemma}\label{lem:BA2:A}
Under the assumptions of Theorem \ref{thm1}, let $v$ be defined by \eqref{err:def} and suppose that \eqref{BA1}-\eqref{BA4} hold.
Then we have that for $t\in [0,T_\ve]$ and $|x|\ge1+t/2$,
\begin{equation}\label{BA2:impr:A}
\sum_{|a|\le2N-17}|\bar\p Z^av|\ls\ve_1\w{t+|x|}^{6\dl_0-1}.
\end{equation}
\end{lemma}
\begin{proof}
First, it follows from \eqref{BA3:impr:A3} that for $|a|\le2N-18$ and $|y|\ge1+s/2$,
\begin{equation}\label{BA2:impr:A1}
|\Box Z^av(s,y)|\ls\ve_1^2\w{s+|y|}^{4\dl_0-1}\w{s-|y|}^{\dl_0-1}.
\end{equation}
By plugging \eqref{initial:data}, \eqref{BA3:impr:A} and \eqref{BA2:impr:A1} into \eqref{hom:pw:pf2} with $w=Z^av$
and $|a|\le2N-17$, one has
that for $|x|\ge1+t/2$,
\begin{equation*}
\begin{split}
|\p_+(r^{1/2}Z^av)(t,x)|&\ls\ve^2\w{x}^{-1}+\ve_1^2(r+t)^{4\dl_0-1/2}\int_0^t(1+|r+t-2s|)^{\dl_0-1}ds\\
&\quad+\ve_1\int_0^t(r+t-s)^{-3/2}\w{r+t}^{6\dl_0}ds\\
&\ls\ve_1\w{t+|x|}^{6\dl_0-1/2},
\end{split}
\end{equation*}
which further implies
\begin{equation*}
\begin{split}
|\p_+Z^av|&\ls r^{-1/2}|\p_+(r^{1/2}Z^av)(t,x)|+r^{-1}|Z^av(t,x)|\\
&\ls\ve_1\w{t+|x|}^{6\dl_0-1}.
\end{split}
\end{equation*}
Together with \eqref{hom:pw:pf2} and \eqref{BA3:impr:A}, this finishes the proof of \eqref{BA2:impr:A}.
\end{proof}

\begin{lemma}
Under the assumptions of Theorem \ref{thm1}, let $v$ be defined by \eqref{err:def} and suppose that \eqref{BA1}-\eqref{BA4} hold.
Then one has
\begin{equation}\label{BA3:impr:B}
\begin{split}
\sum_{|a|\le2N-21}|Z^av|&\ls\ve_1\w{t+|x|}^{8\dl_0-1/2},\\
\sum_{|a|\le2N-21}|\p Z^av|&\ls\ve_1\w{x}^{-1/2}\w{t-|x|}^{-1/2}\w{t+|x|}^{8\dl_0}.
\end{split}
\end{equation}
\end{lemma}
\begin{proof}
With the estimate \eqref{BA3:impr:A}, we will use the equation \eqref{eqn:null} to derive \eqref{BA3:impr:B}.
First, in order to cancel out the quadratic $Q_0$ term in \eqref{eqn:null} while preserve the zero boundary conditions, introducing
the good unknown
\begin{equation}\label{BA3:impr:B1}
V^a=\tilde Z^av-\frac12\sum_{b+c\le a}\sum_{\mu=0}^2C_{abc}^{0,\mu}(\omega)\chi_{[1,2]}(x)\p_{\mu}Z^buZ^cu.
\end{equation}
This, together with \eqref{eqn:null} and \eqref{Q0:ident}, yields $V^a|_{[0,\infty)\times\p\cK}=0$ and
\begin{equation}\label{BA3:impr:B2}
\begin{split}
&\Box V^a=\Box(\tilde Z^a-Z^a)v+F_1^a(t,x),\\
F_1^a&(t,x)=\sum_{b+c\le a}\Big\{\sum_{\mu,\nu=0}^2|x|^{-1}(C(\omega)\p^2_{\mu\nu}Z^bu\Omega Z^cu
+C(\omega)\Omega\p_{\mu}Z^bu\p_{\nu}Z^cu)\\
&+\sum_{\mu=0}^2\Big[C_{abc}^{12,\mu}Q_{12}(\p_{\mu}Z^bu,Z^cu)+C(\omega)\Box Z^bu\p_{\mu}Z^cu\\
&+\frac12[\Delta,C_{abc}^{0,\mu}(\omega)\chi_{[1,2]}(x)](\p_{\mu}Z^buZ^cu)
+C(\omega)\chi_{[1,2]}(x)(\Box\p_{\mu} Z^buZ^cu+\p_{\mu} Z^bu\Box Z^cu)\\
&+\frac12C_{abc}^{0,\mu}(\omega)(1-\chi_{[1,2]}(x))Q_0(\p_{\mu}Z^bu,Z^cu)\Big]\Big\}
+Z^a\big[\sum_{\alpha,\beta,\gamma,\delta=0}^2Q^{\alpha\beta\gamma\delta}\p^2_{\alpha\beta}u\p_{\gamma}u\p_{\delta}u\big].
\end{split}
\end{equation}
For the first term on the fourth line of \eqref{BA3:impr:B2}, recall (2.17) of \cite{HouYinYuan24}
\begin{equation*}
|[\Delta,C_{abc}^{0,\mu}(\omega)\chi_{[1,2]}(x)](\p_{\mu}Z^buZ^cu)|\ls|x|^{-2}|Z^{\le1}(\p_{\mu}Z^buZ^cu)|.
\end{equation*}
Then it can be deduced from \eqref{Q12:null:structure}, \eqref{hom:pw1}, \eqref{BA1}, \eqref{BA3}, \eqref{loc:energyC1} and \eqref{BA3:impr:A} that
\begin{equation}\label{BA3:impr:B3}
\begin{split}
\sum_{|a|\le2N-21}\sum_{|b|\le5}|Z^bF_1^a(s,y)|&\ls\ve_1^2\w{y}^{-3/2}(1+s)^{7\dl_0-1},\quad |y|\le1+s/2,\\
\sum_{|a|\le2N-21}\sum_{|b|\le5}|Z^bF_1^a(s,y)|&\ls\ve_1^2\w{s+|y|}^{7\dl_0-3/2}\w{s-|y|}^{-1},\quad |y|\ge1+s/2.
\end{split}
\end{equation}
Applying \eqref{pw:ibvp} and \eqref{dpw:ibvp} to $\Box V^a$ in \eqref{BA3:impr:B2} with \eqref{BA3:impr:A2}, \eqref{BA3:impr:B3}, $F_0\equiv0$, $\mu=\nu=\dl_0/3$ and $|a|\le2N-21$ derives
\begin{equation*}
\begin{split}
\sum_{|a|\le2N-21}|V^a|&\ls\ve_1^2\w{t+|x|}^{8\dl_0-1/2},\\
\sum_{|a|\le2N-21}|\p V^a|&\ls\ve_1^2\w{x}^{-1/2}\w{t-|x|}^{-1/2}\w{t+|x|}^{8\dl_0}.
\end{split}
\end{equation*}
This, together with \eqref{hom:pw1}, \eqref{BA3}, \eqref{loc:energyC}, \eqref{loc:energyC1} and \eqref{BA3:impr:B1},
yields \eqref{BA3:impr:B}.
\end{proof}

By the second line of \eqref{BA3:impr:B}, the estimate of the first line in \eqref{BA3:impr:B} can be further improved as follows.
\begin{lemma}
Under the assumptions of Theorem \ref{thm1}, let $v$ be defined by \eqref{err:def} and suppose that \eqref{BA1}-\eqref{BA4} hold.
Then one has
\begin{equation}\label{BA3:impr:C}
\sum_{|a|\le2N-27}|Z^av|\ls(\ve_1^2+\ve_1\kappa)\w{t+|x|}^{-1/2}\w{t-|x|}^{-0.4}.
\end{equation}
\end{lemma}
\begin{proof}
Instead of \eqref{BA3:impr:B2}, we use the following equation
\begin{equation}\label{BA3:impr:C1}
\begin{split}
&\Box V^a=\Box(\tilde Z^a-Z^a)v+\p_tF_0^a(t,x)+\tilde F_1^a(t,x),\\
&F_0^a=\frac{1}{3}Z^a\big[\sum_{\alpha,\beta,\gamma,\delta=0}^2
Q^{\alpha\beta\gamma\delta}\omega_{\alpha}\omega_{\beta}\omega_{\gamma}\omega_{\delta}(\p_tu)^3\big]
\end{split}
\end{equation}
and
\begin{equation*}
\begin{split}
\tilde F_1^a&(t,x)=\sum_{b+c\le a}\Big\{\sum_{\mu,\nu=0}^2|x|^{-1}(C(\omega)\p^2_{\mu\nu}Z^bu\Omega Z^cu
+C(\omega)\Omega\p_{\mu}Z^bu\p_{\nu}Z^cu)\\
&+\sum_{\mu=0}^2\Big[C_{abc}^{12,\mu}Q_{12}(\p_{\mu}Z^bu,Z^cu)+C(\omega)\Box Z^bu\p_{\mu}Z^cu
+\frac12C_{abc}^{0,\mu}(\omega)(1-\chi_{[1,2]}(x))Q_0(\p_{\mu}Z^bu,Z^cu)\\
&+\frac12[\Delta,C_{abc}^{0,\mu}(\omega)\chi_{[1,2]}(x)](\p_{\mu}Z^buZ^cu)
+C(\omega)\chi_{[1,2]}(x)(\Box\p_{\mu} Z^buZ^cu+\p_{\mu} Z^bu\Box Z^cu)\Big]\Big\}\\
&+Z^a\big[\sum_{\alpha,\beta,\gamma,\delta=0}^2Q^{\alpha\beta\gamma\delta}(\p^2_{\alpha\beta}u\p_{\gamma}u\p_{\delta}u
-\omega_{\alpha}\omega_{\beta}\omega_{\gamma}\omega_{\delta}\p_t^2u(\p_tu)^2)\big].
\end{split}
\end{equation*}
On the other hand, \eqref{Sobo:ineq} and \eqref{energyC} imply
\begin{equation}\label{BA3:impr:C2}
\sum_{|a|\le2N-15}|\p Z^av(t,x)|\ls\ve_1\w{x}^{-1/2}.
\end{equation}
Together with \eqref{hom:pw1} and \eqref{BA1}, this yields
\begin{equation}\label{BA3:impr:C3}
\sum_{|a|\le2N-27}\sum_{|b|\le6}|Z^bF_0^a(s,y)|\ls\ve_1^3\w{y}^{-3/2}\w{s-|y|}^{2\dl_0-2}.
\end{equation}
Analogously to the estimates \eqref{cubic:null}, \eqref{loc:energyA4} and \eqref{BA3:impr:B3}, based on \eqref{null:structure}, \eqref{Q12:null:structure}, \eqref{null:high}, \eqref{BA2:impr:A}, \eqref{BA3:impr:B} and \eqref{BA3:impr:C2}, one arrives at
\begin{equation}\label{BA3:impr:C4}
\begin{split}
\sum_{|a|\le2N-27}\sum_{|b|\le5}|Z^b\tilde F_1^a(s,y)|&\ls\ve_1^2\w{y}^{-3/2}(1+s)^{11\dl_0-3/2},\quad|y|\le1+s/2,\\
\sum_{|a|\le2N-27}\sum_{|b|\le5}|Z^b\tilde F_1^a(s,y)|&\ls\ve_1^2\w{s+|y|}^{11\dl_0-2}\w{s-|y|}^{-1},\quad|y|\ge1+s/2.
\end{split}
\end{equation}
With \eqref{loc:energyC} instead of \eqref{loc:energyB}, the estimate \eqref{BA3:impr:A2} can be improved to
\begin{equation}\label{BA3:impr:C5}
\sum_{|a|\le2N-27}\sum_{|b|\le4}\sup_{s\in[0,t]}\w{s}^{1-8\dl_0}\|\p^b\Box(\tilde Z^a-Z^a)v(s)\|_{L^2(\cK_4)}
\ls\ve_1^2.
\end{equation}
By using \eqref{pw:ibvp} to $\Box V^a$ in \eqref{BA3:impr:C1} together with \eqref{BA:time}, \eqref{BA3:impr:C3}, \eqref{BA3:impr:C4}, \eqref{BA3:impr:C5}, $\mu=0.4$, $\nu=\dl_0$, $|a|\le2N-27$, $F_0=F_0^a$ and $F_1=\Box(\tilde Z^a-Z^a)v+\tilde F_1^a$, we can achieve
\begin{equation*}
\sum_{|a|\le2N-27}|V^a|\ls(\ve_1^2+\ve_1\kappa)\w{t+|x|}^{-1/2}\w{t-|x|}^{-0.4}.
\end{equation*}
Combining this with \eqref{hom:pw1}, \eqref{BA3}, \eqref{loc:energyC}, \eqref{BA3:impr:B} and \eqref{BA3:impr:B1}
yields \eqref{BA3:impr:C}.
\end{proof}

Based on \eqref{BA3:impr:B}, the estimate \eqref{BA2:impr:A} of the good derivatives can be also refined.
\begin{lemma}
Under the assumptions of Theorem \ref{thm1}, let $v$ be defined by \eqref{err:def} and suppose that \eqref{BA1}-\eqref{BA4} hold.
Then it holds that for $t\in [0,T_\ve]$ and $|x|\ge1+t/2$,
\begin{equation}\label{BA2:impr:B}
\sum_{|a|\le2N-23}|\bar\p Z^av|\ls\ve_1\w{t+|x|}^{9\delta_0-3/2}.
\end{equation}
\end{lemma}
\begin{proof}
First, it follows from \eqref{null:structure}, \eqref{eqn:high}, \eqref{null:high}, \eqref{hom:pw1}, \eqref{BA1}, \eqref{BA2}, \eqref{BA2:impr:A}, \eqref{BA3:impr:B} and \eqref{BA3:impr:C2} that for $|x|\ge1+t/2$ and $|a|\le2N-23$,
\begin{equation}\label{BA2:impr:B1}
\begin{split}
|Z^aQ(\p u,\p^2u)|&\ls\sum_{b+c\le a}(|\bar\p\p Z^bu||\p Z^cu|+|\p^2Z^bu||\bar\p Z^cu|)
+\sum_{b+c+d\le a}|\p^2Z^bu||\p Z^cu||\p Z^du|\\
&\ls\ve_1^2\w{t+|x|}^{6\delta_0-3/2}\w{t-|x|}^{\delta_0-1}+\ve_1^2\w{t+|x|}^{9\delta_0-2}\w{t-|x|}^{-1/2}\\
&\quad+\ve_1^3\w{t+|x|}^{-3/2}\w{t-|x|}^{2\delta_0-2}.
\end{split}
\end{equation}
Similarly to Lemma \ref{lem:BA2:A}, by using \eqref{hom:pw:pf2} to $\Box Z^av=Z^aQ(\p u,\p^2u)$ in \eqref{eqn:high} with $|a|\le2N-23$ and \eqref{BA2:impr:B1}, one has that for $|x|\ge1+t/2$,
\begin{equation*}
\begin{split}
|\p_+(r^{1/2}Z^av)(t,x)|&\ls\ve^2\w{x}^{-1}+\ve_1^2(r+t)^{8\delta_0-1}\int_0^t(1+|r+t-2s|)^{\delta_0-1}ds\\
&\quad+\ve_1\int_0^t(r+t-s)^{-3/2}\w{r+t}^{8\delta_0-1/2}ds\\
&\ls\ve_1\w{t+|x|}^{9\delta_0-1},
\end{split}
\end{equation*}
where \eqref{BA3:impr:B} is used again.
This, together with \eqref{hom:pw:pf2} and \eqref{BA3:impr:B}, leads to \eqref{BA2:impr:B}.
\end{proof}

\subsection{Improved pointwise estimates}

Before improving \eqref{BA1}, the local energy decay estimates \eqref{loc:energyC} will be enhanced.
\begin{lemma}
Under the assumptions of Theorem \ref{thm1}, let $v$ be defined by \eqref{err:def} and suppose that \eqref{BA1}-\eqref{BA4} hold.
Then one has
\begin{equation}\label{loc:energyD}
\sum_{|a|\le2N-25}\|\p^av\|_{L^2(\cK_R)}\ls\ve_1^2\w{t}^{-1}.
\end{equation}
\end{lemma}
\begin{proof}
Instead of the equation \eqref{loc:energyB1}, we have
\begin{equation}\label{loc:energyD1}
\Box V_i=F_0^i(t,x)+F_1^i(t,x),
\end{equation}
where
\begin{equation*}
\begin{split}
&F_0^i(t,x)=\frac{1}{3}\eta(t,x)\p_t^i\big[\sum_{\alpha,\beta,\gamma,\delta=0}^2
Q^{\alpha\beta\gamma\delta}\omega_{\alpha}\omega_{\beta}\omega_{\gamma}\omega_{\delta}(\p_tu)^3\big],\\
&F_1^i(t,x)=\frac{1}{3}(1-\eta(t,x))\p_t^i\big[\sum_{\alpha,\beta,\gamma,\delta=0}^2
Q^{\alpha\beta\gamma\delta}\omega_{\alpha}\omega_{\beta}\omega_{\gamma}\omega_{\delta}(\p_tu)^3\big]\\
&\quad+\p_t^i\big[\sum_{\alpha,\beta,\gamma,\delta=0}^2Q^{\alpha\beta\gamma\delta}(\p^2_{\alpha\beta}u\p_{\gamma}u\p_{\delta}u
-\omega_{\alpha}\omega_{\beta}\omega_{\gamma}\omega_{\delta}\p_t^2u(\p_tu)^2)\big]\\
&\quad+\sum_{j+k\le i}\Big\{\sum_{\mu,\nu=0}^2|x|^{-1}(C(\omega)\p^2_{\mu\nu}\p_t^ju\Omega\p_t^ku
+C(\omega)\Omega\p_{\mu}\p_t^ju\p_{\nu}\p_t^ku)\\
&\quad+\sum_{\mu=0}^2\Big[\frac12\sum_{\mu=0}^2[\Delta,C_{ijk}^{0,\mu}(\omega)](\p_{\mu}\p_t^ju\p_t^ku)
+C_{ijk}^{12,\mu}Q_{12}(\p_{\mu}\p_t^ju,\p_t^ku)\\
&\qquad+C(\omega)\p_{\mu}\p_t^ju\Box\p_t^ku+C(\omega)\p_t^ku\Box\p_{\mu}\p_t^ju\Big]\Big\}.
\end{split}
\end{equation*}
Analogously to  \eqref{BA3:impr:C3}, \eqref{BA3:impr:C4}, one has
\begin{equation}\label{loc:energyD2}
\sum_{i\le2N-28}\sum_{|a|\le5}|\p^aF_0^i(s,y)|\ls\ve_1^3\w{y}^{-3/2}\w{s-|y|}^{2\delta_0-2}
\end{equation}
and
\begin{equation}\label{loc:energyD3}
\begin{split}
\sum_{i\le2N-27}\sum_{|a|\le2}|Z^aF_1^i(s,y)|&\ls\ve_1^2\w{y}^{-3/2}\w{s}^{11\delta_0-3/2},\quad |y|\le1+s/2,~{\rm~or}~|y|\le3,\\
\sum_{i\le2N-27}\sum_{|a|\le2}|Z^aF_1^i(s,y)|&\ls\ve_1^2\w{s+|y|}^{11\delta_0-2}\w{s-|y|}^{-1},\quad |y|\ge1+s/2.
\end{split}
\end{equation}
Applying \eqref{dt:loc} to \eqref{loc:energyD1} together with $\mu=\delta_0$, the definition \eqref{Vi:def}, \eqref{loc:energyD2}
and \eqref{loc:energyD3} yields
\begin{equation*}
\sum_{i\le2N-27}\|\p^{\le1}\p_t\p_t^iv\|_{L^2(\cK_{R_1})}\ls\ve_1^2\w{t}^{-1}.
\end{equation*}
This, together with \eqref{loc:energyA} and the argument in the proof of \eqref{loc:energyB}, derives \eqref{loc:energyD}.
\end{proof}

\begin{lemma}
Under the assumptions of Theorem \ref{thm1}, let $v$ be defined by \eqref{err:def} and suppose that \eqref{BA1}-\eqref{BA4} hold.
Then it holds that
\begin{equation}\label{BA1:impr:A}
\begin{split}
\sum_{|a|\le2N-46}|\p Z^av|&\ls(\ve_1^2+\ve_1\kappa)\w{x}^{-1/2}\w{t-|x|}^{\dl_0-1},\\
\sum_{|a|\le2N-46}|\p Z^av|&\ls\ve_1^2\w{t}^{-1},\quad |x|\le t/10.
\end{split}
\end{equation}
\end{lemma}
\begin{proof}
Let $V^a$ be defined by \eqref{BA3:impr:B1}. It follows from \eqref{eqn:null:div} instead of \eqref{eqn:null} that
\begin{equation}\label{BA1:impr:A1}
\Box V^a=\Box(\tilde Z^a-Z^a)v+\p_t\tilde F_a^0(t,x)+\sum_{\Xi=0,1,2,\theta}\p_{\Xi}G_a^{\Xi}(t,x)+G_{a,1}^r(t,x)+G_{a,2}^r(t,x),
\end{equation}
where
\begin{equation}\label{BA1:impr:A2}
\begin{split}
G_{a,1}^r&(t,x)=\sum_{b+c\le a}\Big\{\frac{C(\omega)Z^cuZ\p Z^bu}{|x|^2}
+\sum_{\mu=0}^2\Big[C(\omega)(1-\chi_{[1,2]}(x))Q_0(\p_{\mu}Z^bu,Z^cu)\\
&+C(\omega)\p_{\mu}Z^cu\sum_{b'+c'\le b}\sum_{\alpha,\beta,\gamma=0}^2Q^{\alpha\beta\gamma}_{bb'c'}\p^2_{\alpha\beta}Z^{b'}u\p_{\gamma}Z^{c'}u
+\frac12[\Delta,C_{abc}^{0,\mu}(\omega)\chi_{[1,2]}(x)]\p_{\mu}Z^buZ^cu\\
&+C(\omega)\chi_{[1,2]}(x)Z^cu\sum_{b'+c'\le b}\sum_{\alpha,\beta,\gamma=0}^2Q^{\alpha\beta\gamma}_{bb'c'}\p_{\mu}(\p^2_{\alpha\beta}Z^{b'}u\p_{\gamma}Z^{c'}u)\\
&+C(\omega)\chi_{[1,2]}(x)\p_{\mu} Z^bu\sum_{b'+c'\le c}\sum_{\alpha,\beta,\gamma=0}^2
Q^{\alpha\beta\gamma}_{cb'c'}\p^2_{\alpha\beta}Z^{b'}u\p_{\gamma}Z^{c'}u\Big]\Big\}\\
&+\frac{1}{3}\p_tZ^a\big[(1-\eta(t,x))\sum_{\alpha,\beta,\gamma,\delta=0}^2
Q^{\alpha\beta\gamma\delta}\omega_{\alpha}\omega_{\beta}\omega_{\gamma}\omega_{\delta}(\p_tu)^3\big]\\
&+Z^a\big[\sum_{\alpha,\beta,\gamma,\delta=0}^2Q^{\alpha\beta\gamma\delta}(\p^2_{\alpha\beta}u\p_{\gamma}u\p_{\delta}u
-\omega_{\alpha}\omega_{\beta}\omega_{\gamma}\omega_{\delta}\p_t^2u(\p_tu)^2)\big]
\end{split}
\end{equation}
and
\begin{equation*}
\begin{split}
\tilde F_a^0&(t,x)=\frac{1}{3}Z^a\big[\eta(t,x)\sum_{\alpha,\beta,\gamma,\delta=0}^2
Q^{\alpha\beta\gamma\delta}\omega_{\alpha}\omega_{\beta}\omega_{\gamma}\omega_{\delta}(\p_tu)^3\big],\\
G_a^{\theta}&(t,x)=\sum_{b+c\le a}\frac{C(\omega)Z^cu\p^2Z^bu}{|x|},\\
G_a^{\alpha}&(t,x)=\sum_{b+c\le a}\frac{C(\omega)Z^cuZ\p Z^bu}{|x|},\quad\alpha=0,1,2,\\
G_{a,2}^r&(t,x)=\sum_{b+c\le a}\sum_{\alpha,\beta,\gamma,\delta,\mu=0}^2\Big\{C(\omega)\p_{\mu}Z^cu\sum_{b'+c'+d\le b}
Q^{\alpha\beta\gamma\delta}_{bb'c'd}\p^2_{\alpha\beta}Z^{b'}u\p_{\gamma}Z^{c'}u\p_{\delta}Z^du\\
&+C(\omega)\chi_{[1,2]}(x)Z^cu\sum_{b'+c'+d\le b}Q^{\alpha\beta\gamma\delta}_{bb'c'd}\p_{\mu}(\p^2_{\alpha\beta}Z^{b'}u\p_{\gamma}Z^{c'}u\p_{\delta}Z^du)\\
&+C(\omega)\chi_{[1,2]}(x)\p_{\mu} Z^bu\sum_{b'+c'+d\le c}
Q^{\alpha\beta\gamma\delta}_{cb'c'd}\p^2_{\alpha\beta}Z^{b'}u\p_{\gamma}Z^{c'}u\p_{\delta}Z^du\Big\}.
\end{split}
\end{equation*}
For the estimates of $\tilde F_a^0$, $G_a^{\Xi}$ ($\Xi=0,1,2,\theta$) and $G_{a,2}^r$, it follows from \eqref{hom:pw1}, \eqref{BA1}, \eqref{BA3}, \eqref{BA3:impr:B}, \eqref{BA3:impr:C} and \eqref{BA3:impr:C2} that
\begin{equation}\label{BA1:impr:A3}
\begin{split}
\sum_{|a|\le2N-36}\sum_{|b|\le13}|Z^b\tilde F_a^0(s,y)|&\ls\ve_1^3\w{y}^{-3/2}\w{s-|y|}^{2\dl_0-2},\\
\sum_{|a|\le2N-36}\sum_{|b|\le10}|Z^bG_a^{\Xi}(s,y)|&\ls\ve_1^2\w{y}^{-3/2}\w{s+|y|}^{8\dl_0-1/2}\w{s-|y|}^{3\dl_0-1},\\
\sum_{|a|\le2N-36}\sum_{|b|\le9}|Z^bG_{a,2}^r(s,y)|&\ls\ve_1^4\w{y}^{-3/2}\w{s+|y|}^{-1/2}\w{s-|y|}^{-2}.
\end{split}
\end{equation}
Next, we focus on the estimate of $G_{a,1}^r$ in \eqref{BA1:impr:A2}.
In the region $|y|\le1+s/2$, by using \eqref{hom:pw1}, \eqref{BA1}, \eqref{BA3} and \eqref{BA3:impr:B}, one obtains
\begin{equation}\label{BA1:impr:A4}
\begin{split}
\sum_{|a|\le2N-36}\sum_{|b|\le9}|Z^bG_{a,1}^r(s,y)|&\ls\ve_1^2\w{y}^{-5/2}\w{s+|y|}^{8\delta_0-1/2}\w{s-|y|}^{3\delta_0-1}\\
&\quad+\ve_1^3\w{y}^{-1}\w{s+|y|}^{8\delta_0-1/2}\w{s-|y|}^{4\delta_0-2}.
\end{split}
\end{equation}
In the region $|y|\ge1+s/2$, since the treatment on the last term in $G_{a,1}^r$ is similar to that of \eqref{cubic:null},
 then it follows from the first line of \eqref{BA1:impr:A4} together with \eqref{null:structure}, \eqref{hom:pw1}, \eqref{BA1}-\eqref{BA3}, \eqref{BA3:impr:B} and \eqref{BA2:impr:B} that
\begin{equation}\label{BA1:impr:A5}
\begin{split}
\sum_{|a|\le2N-36}\sum_{|b|\le9}|Z^bG_{a,1}^r(s,y)|&\ls\ve_1^2\w{y}^{-5/2}\w{s+|y|}^{8\delta_0-1/2}\w{s-|y|}^{3\delta_0-1}\\
&\quad+\ve_1^3\w{y}^{-1/2}\w{s+|y|}^{9\delta_0-2}\w{s-|y|}^{3\delta_0-1}.
\end{split}
\end{equation}
On the other hand, by using \eqref{loc:energyD}, the estimate \eqref{BA3:impr:C5} can be further improved to
\begin{equation}\label{BA1:impr:A6}
\sum_{|a|\le2N-36}\sum_{|b|\le8}\sup_{s\in[0,t]}\w{s}\|\p^b\Box(\tilde Z^a-Z^a)v(s)\|_{L^2(\cK_4)}\ls\ve_1^2.
\end{equation}
Collecting \eqref{BA1:impr:A1}-\eqref{BA1:impr:A6} with \eqref{dpw:div}, $m^{\#}=1$, $|a|\le2N-36$, $\eta=\mu=\nu=\delta_0/3$
and $G_1^r=\Box(\tilde Z^a-Z^a)v+G_{a,1}^r$ leads to
\begin{equation}\label{BA1:impr:A7}
\sum_{|a|\le2N-36}|\p Z^av(t,x)|\ls\ve_1\w{x}^{-1/2}\w{t-|x|}^{\delta_0-1}\w{t}^{12\delta_0},
\end{equation}
where we have used the definition \eqref{BA3:impr:B1}, \eqref{hom:pw1}, \eqref{BA3}, \eqref{BA:time}, \eqref{BA3:impr:C} and \eqref{loc:energyD}.
To obtain \eqref{BA1:impr:A}, one needs to improve $\|G_a^i(s)\|_{L^\infty(\cK_3)}$ ($i=1,2$). It follows from the first line of \eqref{BA1:impr:A4} with \eqref{BA1:impr:A7} instead of \eqref{BA3:impr:B}
\begin{equation}\label{BA1:impr:A8}
\begin{split}
\sum_{|a|\le2N-46}\sum_{|b|\le9}|Z^bG_{a,1}^r(s,y)|&\ls\ve_1^2\w{y}^{-5/2}\w{s+|y|}^{12\delta_0-1/2}\w{s-|y|}^{-1.4}\\
&\quad+\ve_1^3\w{y}^{-1}\w{s+|y|}^{8\delta_0-1/2}\w{s-|y|}^{4\delta_0-2},\\
\sum_{|a|\le2N-46}\sum_{i=1,2}\|G_a^i(s)\|_{L^\infty(\cK_3)}&\ls\ve_1^2\w{s}^{-2},
\end{split}
\end{equation}
where \eqref{loc:energyD} is used.
Therefore, analogously to \eqref{BA1:impr:A7}, applying \eqref{dpw:div} with $m^{\#}=0,1$, $|a|\le2N-46$, \eqref{BA1:impr:A1}-\eqref{BA1:impr:A3}, \eqref{BA1:impr:A5}, \eqref{BA1:impr:A6} and \eqref{BA1:impr:A8} yields \eqref{BA1:impr:A}.
This completes the proof of \eqref{BA1:impr:A}.
\end{proof}

\begin{lemma}
Under the assumptions of Theorem \ref{thm1}, let $v$ be defined by \eqref{err:def} and suppose that \eqref{BA1}-\eqref{BA4} hold.
Then we have
\begin{equation}\label{BA2:impr:C}
\sum_{|a|\le2N-46}|\bar\p Z^av|\ls(\ve_1^2+\ve_1\kappa)\w{x}^{-1/2}\w{t+|x|}^{\dl_0-1}
\end{equation}
and
\begin{equation}\label{BA3:impr:D}
\sum_{|a|\le2N-52}|Z^av|\ls(\ve_1^2+\ve_1\kappa)\w{t+|x|}^{-1/2}\w{t-|x|}^{3\dl_0-1/2}.
\end{equation}
\end{lemma}
\begin{proof}
The proof is easier than that of \eqref{BA3:impr:C} and \eqref{BA2:impr:B}.
First, we prove \eqref{BA2:impr:C} in the region $|x|\ge1+t/2$ since
\eqref{BA2:impr:C} can be achieved directly in the region $|x|\le1+t/2$  by \eqref{BA1:impr:A}.
Note that the estimate \eqref{BA2:impr:B1} may be improved with \eqref{BA2:impr:B}, \eqref{BA1:impr:A7} instead of \eqref{BA2:impr:A}, \eqref{BA3:impr:B}, respectively,
\begin{equation*}\label{BA3:impr:D1}
\begin{split}
\sum_{|a|\le2N-46}|Z^aQ(\p u,\p^2u)|&\ls\ve_1^2\w{t+|x|}^{13\delta_0-2}\w{t-|x|}^{\delta_0-1}\\
&+\ve_1^3\w{t+|x|}^{-3/2}\w{t-|x|}^{2\delta_0-2},\quad |x|\ge1+t/2.
\end{split}
\end{equation*}
Thus, \eqref{BA2:impr:C} is obtained by \eqref{BA3:impr:C} and the same argument in the proof of \eqref{BA2:impr:B}.

Due to
\begin{equation}\label{BA3:impr:D2}
\begin{split}
&\Box\tilde Z^av=\Box(\tilde Z^a-Z^a)v+\p_tF_0^a+\breve{F_1^a},\\
&\breve{F_1^a}=\sum_{\alpha,\beta,\gamma=0}^2\sum_{b+c\le a}Q^{\alpha\beta\gamma}_{abc}\p^2_{\alpha\beta}Z^bu\p_{\gamma}Z^cu\\
&\quad +Z^a\big[\sum_{\alpha,\beta,\gamma,\delta=0}^2Q^{\alpha\beta\gamma\delta}\p^2_{\alpha\beta}u\p_{\gamma}u\p_{\delta}u
-\omega_{\alpha}\omega_{\beta}\omega_{\gamma}\omega_{\delta}\p_t^2u(\p_tu)^2\big],
\end{split}
\end{equation}
where $F_0^a$ is defined in \eqref{BA3:impr:C1},
then analogously to \eqref{cubic:null}, it follows from \eqref{hom:pw1}, \eqref{BA1:impr:A} and \eqref{BA2:impr:C} that
\begin{equation}\label{BA3:impr:D3}
\sum_{|a|\le2N-52}\sum_{|b|\le5}|Z^b\breve{F_1^a}|\ls\ve_1^2\w{x}^{-1}\w{t+|x|}^{\dl_0-1}\w{t-|x|}^{\dl_0-1}.
\end{equation}
As in the proof of \eqref{BA3:impr:C}, combining \eqref{pw:ibvp} with $\Box Z^av=\Box Z^au$, $|a|\le2N-52$, $\mu=1/2-3\delta_0$, $\nu=\delta_0$, \eqref{BA:time}, \eqref{BA3:impr:C3}, \eqref{loc:energyD} and \eqref{BA3:impr:D2} derives \eqref{BA3:impr:D}.
\end{proof}

\subsection{Proof of Theorem \ref{thm1}}
\begin{proof}[Proof of Theorem \ref{thm1}]
Collecting \eqref{BA1:impr:A}, \eqref{BA2:impr:C} and \eqref{BA3:impr:D} yields that there is a constant $C_1\ge1$ such that
for $t\in [0,e^{\kappa/\ve^2})$,
\begin{align*}
\sum_{|a|\le2N-46}|\p Z^av|&\le C_1(\ve_1^2+\ve_1\kappa)\w{x}^{-1/2}\w{t-|x|}^{\delta_0-1},\\
\sum_{|a|\le2N-46}|\p Z^av|&\le C_1\ve_1^2\w{t}^{-1},\quad |x|\le t/10,\\
\sum_{|a|\le2N-46}|\bar\p Z^av|&\le C_1(\ve_1^2+\ve_1\kappa)\w{x}^{-1/2}\w{t+|x|}^{\delta_0-1},\\
\sum_{|a|\le2N-52}|Z^av|&\le C_1(\ve_1^2+\ve_1\kappa)\w{t+|x|}^{-1/2}\w{t-|x|}^{3\delta_0-1/2}.
\end{align*}
Choosing $\ve_1=\ve$, $\ve_2=\sqrt{2C_1}\ve_1$ and $\ve_0=\kappa_0=\frac{1}{4C_1}$. Then for $N\ge53$, \eqref{BA:time} is verified and \eqref{BA1}-\eqref{BA4} can be improved to
\begin{align*}
\sum_{|a|\le N+1}|\p Z^av|&\le \frac12\ve_1\w{x}^{-1/2}\w{t-|x|}^{\delta_0-1},\\
\sum_{|a|\le N+1}|\p Z^av|&\le \frac12\ve_2^2\w{t}^{-1},\quad |x|\le t/10,\\
\sum_{|a|\le N+1}|\bar\p Z^av|&\le \frac12\ve_1\w{x}^{-1/2}\w{t+|x|}^{\delta_0-1},\\
\sum_{|a|\le N+1}|Z^av|&\le \frac12\ve_1\w{t+|x|}^{-1/2}\w{t-|x|}^{3\delta_0-1/2}.
\end{align*}
This, together with the local existence of classical solution to the initial boundary value problem
of the hyperbolic equation  under the compatibility conditions of initial boundary values, yields that problem \eqref{QWE} has a solution $u\in\bigcap\limits_{j=0}^{2N+1}C^{j}([0,e^{\kappa/\ve^2}), H^{2N+1-j}(\cK))$.
In addition, \eqref{thm1:decay:a}, \eqref{thm1:decay:b}, \eqref{thm1:decay:c}, \eqref{thm1:decay:d}, \eqref{thm1:LEdecay:a} and \eqref{thm1:LEdecay:b} can be obtained by \eqref{BA1:impr:A}, \eqref{BA2:impr:C}, \eqref{BA3:impr:D} and \eqref{loc:energyD} with \eqref{hom:pw1}, \eqref{hom:pw:pf1}, respectively.
\end{proof}

\vskip 0.2 true cm

{\bf \color{blue}{Conflict of Interest Statement:}}

\vskip 0.2 true cm

{\bf The authors declare that there is no conflict of interest in relation to this article.}

\vskip 0.2 true cm

{\bf \color{blue}{Data availability statement:}}

\vskip 0.2 true cm
{\bf  Data sharing is not applicable to this article as no data sets are generated
during the current study.}

\end{document}